\documentclass[nointlimits,11pt,oneside]{amsart}
\usepackage[utf8]{inputenc}
\usepackage[T1]{fontenc}
\usepackage{lmodern}
\usepackage{amssymb, amsmath, amsthm, mathtools, cases,enumerate}
\usepackage{color}
\usepackage{hyperref,mathscinet}
\usepackage{natbib}
\usepackage{microtype}

\DeclareFontShape{OMX}{cmex}{m}{n}{
  <-7.5> cmex7
  <7.5-8.5> cmex8
  <8.5-9.5> cmex9
  <9.5-> cmex10
}{}
\SetSymbolFont{largesymbols}{normal}{OMX}{cmex}{m}{n}
\SetSymbolFont{largesymbols}{bold}  {OMX}{cmex}{m}{n}

\usepackage{hyphenat}

\usepackage[%
	a4paper,
	total={16cm,23cm},
	top=3cm,
	marginparsep=2pt]
{geometry}

\hypersetup{
	unicode=true,
	breaklinks=true,
	pdftitle={Optimal behavior of weighted Hardy operators on rearrangement-invariant spaces},
	pdfauthor={Zdeněk Mihula}
	}

\theoremstyle{plain}
\newtheorem{theorem}{Theorem}[section]
\newtheorem{corollary}[theorem]{Corollary}
\newtheorem{proposition}[theorem]{Proposition}

\theoremstyle{definition}
\newtheorem{remark}[theorem]{Remark}
\newtheorem{remarks}[theorem]{Remarks}
\newtheorem*{conventions}{Conventions and notation}

\newcommand{\R}{\mathbb{R}}
\newcommand{\N}{\mathbb{N}}

\newcommand{\myref}[2]{\hyperref[#2]{#1~\ref*{#2}}}

\newcommand{\rn}{\R^n}
\newcommand{\dx}{{\fam0 d}}
\renewcommand{\d}[1]{\,\dx #1}

\newcommand{\M}{\mathfrak{M}}
\newcommand{\Mpl}{\M^+}

\newcommand{\Dinf}[2][]{\vphantom{D}\smash{\underline{D}}_{#1}^{#2}}
\newcommand{\Dsup}[2][]{\vphantom{D}\smash{\overline{D}}_{#1}^{#2}}

\DeclarePairedDelimiter\ceil{\lceil}{\rceil}

\DeclareMathOperator*{\esssup}{ess\,sup}
\DeclareMathOperator{\id}{id}
\DeclareMathOperator{\Int}{Int}
\DeclareMathOperator{\K}{K}

\allowdisplaybreaks

\title[Optimal behavior of weighted Hardy operators on r.i.~spaces]{Optimal behavior of weighted Hardy operators on rearrangement\hyp{}invariant spaces}
\author{Zden\v ek Mihula}
\date{\today}
\address{Zden\v ek Mihula, Czech Technical University in Prague, Faculty of Electrical Engineering, Department of Mathematics, Technick\'a~2, 166~27 Praha~6, Czech Republic --- AND --- Charles University, Faculty of Mathematics and Physics, Department of Mathematical Analysis, Sokolovsk\'a~83, 186~75 Praha~8, Czech Republic}
\email{mihulzde@fel.cvut.cz}
\email{mihulaz@karlin.mff.cuni.cz}
\urladdr{\href{https://orcid.org/0000-0001-6962-7635}{0000-0001-6962-7635}}

\numberwithin{equation}{section}

\begin{document}
\setcitestyle{numbers}
\bibliographystyle{plainnat}


\thanks{This research was supported by the project OPVVV CAAS CZ.02.1.01/0.0/0.0/16\_019/0000778 and by the grant SVV-2020-260583.}

\begin{abstract}
The behavior of certain weighted Hardy-type operators on rearrangement\hyp{}invariant function spaces is thoroughly studied with emphasis being put on the optimality of the obtained results. First, the optimal rearrangement\hyp{}invariant function spaces\textemdash that is, the best possible function spaces within the class of rearrangement\hyp{}invariant function spaces\textemdash guaranteeing the boundedness of the operators from/to a given rearrangement\hyp{}invariant function space are described. Second, the optimal rearrangement\hyp{}invariant function norms being sometimes complicated, the question of whether and how they can be simplified to more manageable expressions, arguably more useful in practice, is addressed. Last, iterated weighted Hardy-type operators are also studied.

Besides aiming to provide a comprehensive treatment of the optimal behavior of the operators on rearrangement\hyp{}invariant function spaces in one place, the paper is motivated by its applicability in various fields of mathematical analysis, such as harmonic analysis, extrapolation theory or the theory of Sobolev-type spaces.
\end{abstract}

\maketitle

\section{Introduction}
When we face a complicated problem, it is only natural to seek for a way to simplify it. Problems (not only) in mathematical analysis often amount to establishing boundedness of operators between certain function spaces. When the function spaces in question are endowed with norms invariant with respect to certain rearrangements/transformations of functions, an often successful way to simplify such problems is make use of the rearrangements. We shall now make this vague idea more definite.

Arguably the most straightforwardly, the idea can be illustrated by the following well-known example, which traces back to the 1930s. Consider the question of establishing the boundedness of the Hardy--Littlewood maximal operator $M$, defined for measurable functions $f$ on $\rn$ as
\begin{equation*}
Mf(x) = \sup_{x \ni Q}\frac1{|Q|}\int_Q |f(y)|\d{y},\ x\in\rn,
\end{equation*}
in which the supremum extends over all cubes in $\rn$ whose edges are parallel to the coordinate axes, from a function space to another. The famous inequality
\begin{equation}\label{intro:M_upper}
(Mf)^*(t)\lesssim \frac1{t}\int_0^t f^*(s)\d{s} \quad \text{for every $t\in(0,\infty)$}
\end{equation}
by F.~Riesz (\cite{R:32}, $n=1$) and N.~Wiener (\cite{W:39}, $n\in\N$), in which $^*$ denotes the nonincreasing rearrangement and the multiplicative constant depends only on $n$, combined with the classical Hardy--Littlewood inequality (\cite[p.~278]{HLP:52}, see also \eqref{ch1:ri:HLg=chiE} below), tells us that, when the function spaces in question are such that equimeasurable functions (i.e., functions whose distribution functions coincide) have the same norm, it is sufficient to establish the boundedness of a considerably simpler operator, acting on functions of a single variable; namely of the operator $R$ defined for measurable functions $g$ on $(0,\infty)$ as
\begin{equation*}
Rg(t) = \frac1{t}\int_0^t|g(s)|\d{s},\ t\in(0,\infty).
\end{equation*}
Moreover, since the reverse inequality to \eqref{intro:M_upper} holds, too, which was proved by C.~Hertz (\cite{H:68}, $n=1$) and by C.~Bennett and R.~Sharpley (\cite{BS:79}, $n\in\N$), the boundedness of $M$ between such function spaces is actually equivalently reduced to the boundedness of the Hardy-type operator $R$. In particular, this reduction is valid for rearrangement\hyp{}invariant function spaces (see \myref{Section}{sec:prel} for precise definitions), which constitute a broad class of functions spaces, containing Lebesgue spaces, Orlicz spaces or Lorentz(--Zygmund) spaces to name a few.

Another important operator of harmonic analysis, hardly needing an introduction, is the Hilbert transform (e.g., \cite[Chapter~3, Definition~4.1]{BS}). An inequality for the nonincreasing rearrangement of the Hilbert transform (more generally, of certain singular integral operators with odd kernels) is known (\cite[Theorem~16.12]{BR:80}, \cite[Lemma~2.1]{S:80}, cf.~\cite[p.~55]{C:97}), and, since the inequality is in a suitable sense sharp (cf.~\cite[p.~29]{EOP:96}), it is easy to show that the boundedness of the Hilbert transform on rearrangement\hyp{}invariant function spaces is equivalent to the boundedness of a sum of two Hardy-type operators acting on functions of a single variable\textemdash namely $R+H$, where $H$ is defined for measurable functions $g$ on $(0,\infty)$ as
\begin{equation*}
Hg(t)=\int_t^\infty|g(s)|\frac1{s}\d{s},\ t\in(0,\infty).
\end{equation*}
Other important operators of harmonic analysis for which (in a suitable sense) sharp inequalities that involve weighted Hardy-type operators being of the form
\begin{align}
(0,\infty)\ni t\mapsto v(t)\int_0^t|g(s)|\d{s}, \label{intro:Rv}
\intertext{and}
(0,\infty)\ni t\mapsto \int_0^t|g(s)|v(s)\d{s}, \label{intro:Hv}
\end{align}
in which $v$ is a fixed positive nonincreasing function on $(0,\infty)$ and $g$ is a measurable function on $(0,\infty)$ on which they act, for their nonincreasing rearrangements are known are certain convolution operators (\cite{On:63,EOP:96}), of which the Riesz potential is a prototypical example, or the fractional maximal operator and its variants (\cite[Theorem~1.1]{CKOP:00}, \cite[Theorem~3.1]{EO:02}). The interested reader can find more information on boundedness of some classical operators of harmonic analysis on rearrangement\hyp{}invariant function spaces in \cite{EMMP:20}.

Much as sharp inequalities for the nonincreasing rearrangements of operators are desired, the number of operators for which such sharp inequalities are known is limited. Nevertheless, what is often at our disposal is at least an upper bound on the nonincreasing rearrangement of a given operator $T$. Obviously, the better we control the upper bound, the better we control the operator $T$. Weighted Hardy-type operators such as those defined by \eqref{intro:Rv} and \eqref{intro:Hv} (or their more complicated forms, which are to be introduced soon) for suitable functions $v$ often serve as such upper bounds. It is worth noting that such upper bounds for various maximal operators may actually involve a Hardy-type operator inside a supremum (see \cite{L:05} and references therein), but the supremum usually does not cause any trouble (see \cite[Lemma~4.10]{EMMP:20}, cf.~\cite[Theorem~3.9]{KP:06}). For example, consider the fractional maximal operator $M_\gamma$ defined for measurable functions $f$ on $\rn$, $\gamma\in(0,n)$, as
\begin{equation*}
M_\gamma f(x) = \sup_{x \ni Q}\frac1{|Q|^{1-\frac{\gamma}{n}}}\int_Q |f(y)|\d{y},\ x\in\rn,
\end{equation*}
in which the supremum extends over all cubes in $\rn$ whose edges are parallel to the coordinate axes. The following pointwise inequality, which is also in a suitable sense sharp, is valid (\cite[Example~1]{L:05}, \cite[Theorem~1.1]{CKOP:00}):
\begin{equation*}
(M_\gamma f)^*(t) \lesssim \sup_{t< s < \infty}s^{\frac{\gamma}{n} - 1}\int_0^s f^*(\tau) \d{\tau} \quad \text{for every $t\in(0,\infty)$},
\end{equation*}
in which the multiplicative constant depends only on $\gamma$ and $n$. Now, while the operators
\begin{equation*}
f\mapsto \sup_{(\cdot)< s < \infty}s^{\frac{\gamma}{n} - 1}\int_0^s f^*(\tau) \d{\tau} \quad \text{and} \quad f\mapsto (\cdot)^{\frac{\gamma}{n} - 1}\int_0^{(\cdot)} f^*(s) \d{s}
\end{equation*}
are not pointwise equivalent, it turns out that their rearrangement\hyp{}invariant function norms are equivalent; it follows from \cite[Lemma~4.10]{EMMP:20} (combined with the Hardy--Littlewood inequality) that the supremum operator is bounded from a rearrangement\hyp{}invariant function function space to a rearrangement\hyp{}invariant function space if and only if the integral Hardy-type operator defined by \eqref{intro:Rv} with $v(t)=t^{\gamma/n - 1}$ is.

Reductions of complicated questions (often involving functions of several variables) to simpler ones of whether certain Hardy-type operators (acting on functions of a single variable) are bounded, are sometimes also achieved with the right use of interpolation or by making use of some intrinsic properties of the problem in question (such as symmetrization principles and isoperimetric inequalities, which are of great use in studying Sobolev spaces). Such approaches have been notably successful in connection with various embeddings of Sobolev-type spaces built upon rearrangement\hyp{}invariant function spaces into rearrangement\hyp{}invariant function spaces. For a wide variety of such embeddings, either complete characterizations (\cite{KP:06, CPS:20, CPS:15, CP:16, CP:09, CKP:08, ACPL:18, M:21, BC:21}) or at least sufficient and/or necessary conditions (\cite{CiaMaz:16, CiaMaz:20, CPS:20, M:21b}) for their validity by means of boundedness of Hardy-type operators have been obtained. For example, consider the Sobolev-type embedding
\begin{equation}\label{intr:Sob_emb}
W^m X(\Omega) \hookrightarrow Y(\overline{\Omega},\nu),
\end{equation}
in which $W^mX(\Omega)$ is the $m$th order Sobolev space built upon a rearrangement\hyp{}invariant function space $X$ over a bounded Lipschitz domain $\Omega$ in $\rn$, $m<n$, $m\in\N$, and $Y$ is a rearrangement\hyp{}invariant function space over $\overline{\Omega}$ endowed with a $d$-upper Ahlfors measure $\nu$, that is, a finite Borel measure $\nu$ on $\overline{\Omega}$ satisfying
\begin{equation*}
\sup_{x\in\rn, r>0}\frac{\nu(B_r(x) \cap \overline{\Omega})}{r^d}<\infty
\end{equation*}
with $d\in(0, n]$, in which $B_r(x)$ is the open ball centered at $x$ with radius $r$. When $d\in[n-m,n]$, it was shown in \cite[Theorem~4.1]{CPS:20} that the embedding \eqref{intr:Sob_emb} is valid (in the sense of traces) if the Hardy-type operator defined for measurable functions $g$ on $(0,1)$ as
\begin{equation*}
(0,1)\ni t\mapsto \int_{t^\frac{n}{d}}^1 |g(s)| s^{-1+\frac{m}{n}}\d{s}
\end{equation*}
is bounded from $X(0,1)$ to $Y(0,1)$ (the representation spaces of $X$ and $Y$ over the interval $(0,1)$). Moreover, if there is a point $x_0\in\overline{\Omega}$ for which the exponent $d$ is sharp, that is,
\begin{equation*}
\inf_{r\in(0,1)}\frac{\nu(B_r(x_0) \cap \overline{\Omega})}{r^d}>0,
\end{equation*}
then the boundedness of the Hardy-type operator is also necessary for the validity of \eqref{intr:Sob_emb} (\cite[Theorem~4.3]{CPS:20}). Finally, when $d\in(0,n-m)$, it follows from \cite[Theorem~5.1]{CPS:20} that the embedding \eqref{intr:Sob_emb} is valid if two weighted Hardy-type operators are bounded from $X(0,1)$ to $Y(0,1)$\textemdash namely the same one as in the case $d\in[n-m,n]$ and the one defined for measurable functions $g$ on $(0,1)$ as
\begin{equation*}
(0,1)\ni t\mapsto t^{-\frac{m}{n-d}}\int_0^{t^\frac{n}{d}} |g(s)| s^{-1+\frac{m}{n-d}}\d{s}.
\end{equation*}
We note for the interested reader that, if this is the case, then the rearrangement\hyp{}invariant function space $Y$ in \eqref{intr:Sob_emb} can actually be replaced by a rearrangement\hyp{}invariant function space that is in some cases smaller than $Y$ (see \cite{T:preprint}).

In this paper, we introduce the Hardy-type operators $R_{u,v,\nu}$ and $H_{u,v,\nu}$, of which Hardy-type operators mentioned in the preceding paragraphs are special instances, defined for measurable functions $g$ on $(0,L)$, $L\in(0,\infty]$, as
\begin{align}
R_{u,v,\nu}g(t)&=v(t)\int_0^{\nu(t)}|g(s)|u(s)\d{s},\ t\in(0,L), \label{opRdef}
\intertext{and}
H_{u,v,\nu}g(t)&=u(t)\int_{\nu(t)}^L|g(s)|v(s)\d{s},\ t\in(0,L), \label{opHdef}
\end{align}
where $u, v$ are nonnegative nonincreasing functions on $(0,L)$ and $\nu$ is an increasing bijection of the interval $(0,L)$ onto itself, and thoroughly study their behavior on rearrangement\hyp{}invariant function spaces, putting emphasis on the optimality of our results. First, after fixing some notation and recalling some preliminary results in \myref{Section}{sec:prel}, we characterize the optimal domain and the optimal target rearrangement\hyp{}invariant function space for the operators when the rearrangement\hyp{}invariant function space on the other side is fixed in \myref{Section}{sec:optimal_norms}. More precisely, given one of the operators and a rearrangement\hyp{}invariant function space $X$, we describe the largest and the smallest rearrangement\hyp{}invariant function space $Y$ (in other words, the weakest and the strongest rearrangement\hyp{}invariant function norm $\|\cdot\|_Y$) such that the operator is bounded from $Y$ to $X$ and from $X$ to $Y$, respectively. As a simple corollary, we also obtain a description of the optimal rearrangement\hyp{}invariant function spaces for a sum of the two operators (each with possibly different functions $u,v,\nu$), though the description is less explicit than it could be if we studied directly the sum. Next, in \myref{Section}{sec:simplification}, we take a close look at how to simplify the description of the optimal rearrangement\hyp{}invariant function norms for the Hardy-type operators and whether it is possible, for the simpler and more manageable description we have at our disposal, the more useful it is. It turns out that this problem is more complex than it may appear at first glance, as is often the case in mathematics. Last, in \myref{Section}{sec:iteration}, we investigate the optimal behavior of iterated Hardy-type operators\textemdash namely $R_{u_1,v_1,\nu_1}\circ R_{u_2,v_2,\nu_2}$ and $H_{u_1,v_1,\nu_1}\circ H_{u_2,v_2,\nu_2}$\textemdash on rearrangement\hyp{}invariant function spaces, explaining why it is of interest to study such operators at the beginning of the section.

In considerably less general settings, the questions mentioned in the preceding paragraph were already studied in some of the papers cited in the fifth paragraph and in \cite{EMMP:20} (cf.~\cite{ST:16, OS:07}); however, not only are those results limited to some particular choices of the functions $u$, $v$ and $\nu$ (namely, $u\equiv1$ and $v,\nu$ being power functions for the most part, but see \cite{GL:99, HS:98}), but they are also scattered and often hidden somewhere between the lines with varying degrees of generality. The aim of this paper is thoroughly address the questions in a coherent unified way so that the results obtained here will not only encompass their already-known particular cases but also provide a general theory suitable for various future applications. For example, thanks to our results, we get under control the optimal behavior of upper bounds, which are even sharp in some cases, for nonincreasing rearrangements of not only various less-standard (nonfractional and fractional) maximal operators (\cite{EO:02, L:05}) but also operators with certain behavior of their operator norms, which play a role in the a.e.~convergence of the partial spherical Fourier integrals or in the  solvability of the Dirichlet problem for the Laplacian on planar domains (\cite{CaPr:09, CO-C:18, CMM:12} and references therein). Another possible application is related to traces of Sobolev functions. Although, as was already mentioned in the fifth paragraph, embeddings of Sobolev-type spaces into rearrangement\hyp{}invariant function spaces with respect to $d$-upper Ahlfors measures were thoroughly studied in \cite{CPS:20}, there are $d$-dimensional sets (say, in the sense of the Hausdorff dimension) $\Omega_d\subseteq\rn$, $d\in(0,n]$, that are ``unrecognizable'' by $d$-upper Ahlfors measures $\nu$ (i.e., it may happen that $\nu(\Omega_d)=0$ for every $d$-upper Ahlfors measure $\nu$). For instance, this is (almost surely) the case when $\Omega_d$ is a Brownian path in $\rn$, $n\geq2$, which has (almost surely) the Hausdorff dimension $2$ but is unrecognizable by $2$-upper Ahlfors measures; other measure functions (see~\cite{C:67}) than power functions have to be considered to rectify the situation (\cite{ET:61, CT:62}). Inevitably, if one is to generalize the results of \cite{CPS:20} to cover such exceptional sets, one will need to deal with general enough Hardy-type operators, whose optimal behavior on rearrangement\hyp{}invariant function spaces is extensively studied here.

General as the results proved in this paper are, we do usually impose some mild restrictions on the functions $u$, $v$, $\nu$ so that we can obtain interesting, strong results; however, the imposed assumptions on the functions are actually not too restrictive for the most part and often exclude only cases being in a way pathological. The assumptions also often reflect the very forms of the Hardy-type operators considered here. In particular, the operators do not involve kernels. Although Hardy-type operators with kernels are undoubtedly of great importance, too, and to investigate thoroughly their behavior on rearrangement\hyp{}invariant function spaces would be of interest (e.g., \cite{CPS:15, AAB-MC:preprint}), it goes beyond the scope of this paper.

Finally, what is missing in this paper are particular examples of optimal function norms and related things. Given the generality of our setting, to provide an exhaustive list of examples would inevitably involve carrying out a large number of, to some extent straightforward but lengthy and technical, computations, which would make this paper unreasonably long and convoluted. The interested reader is referred to \cite{CP:16, EMMP:20, M:21} for some particular examples and hints on how to carry out the needed computations.

\section{Preliminaries}\label{sec:prel}
\begin{conventions}\hphantom{}
\begin{itemize}
\item Throughout the entire paper $L\in(0,\infty]$.
\item We adhere to the convention that $\frac1{\infty}=0\cdot\infty=0$.
\item We write $P\lesssim Q$, where $P,Q$ are nonnegative quantities, when there is a positive constant $c$ independent of all appropriate quantities appearing in the expressions $P$ and $Q$ such that $P\leq c\cdot Q$. If not stated explicitly,  what ``the appropriate quantities appearing in the expressions $P$ and $Q$'' are should be obvious from the context. At the few places where it is not obvious, we will explicitly specify what the appropriate quantities are. We also write $P\gtrsim Q$ with the obvious meaning, and $P\approx Q$ when $P\lesssim Q$ and $P\gtrsim Q$ simultaneously.
\item When $A\subseteq(0,L)$ is a (Lebesgue) measurable set, $|A|$ stands for its Lebesgue measure.
\item When $u$ is a nonnegative measurable function defined on $(0,L)$, we denote by $U$ the function defined as $U(t)=\int_0^t u(s)\d{s}$, $t\in(0,L]$. We say that $u$ is \emph{nondegenerate} if there is $t_0\in(0 ,L)$ such that $0< U(t_0) < \infty$.
\end{itemize}
\end{conventions}

We set
\begin{align*}
\M(0, L)&= \{f\colon \text{$f$ is a measurable function on $(0,L)$ with values in $[-\infty,\infty]$}\},\\
\M_0(0,L)&= \{f \in \M(0, L)\colon f \ \text{is finite}\ \text{a.e.~on $(0,L)$}\},\\
\intertext{and}
\Mpl(0,L)&= \{f \in \M(0,L)\colon f \geq 0\ \text{a.e.~on $(0,L)$}\}.
\end{align*}

The \emph{nonincreasing rearrangement} $f^*\colon(0,\infty) \to [0, \infty]$ of a function $f\in \M(0,L)$  is defined as
\begin{equation*}
f^*(t)=\inf\{\lambda\in(0,\infty)\colon |\{s\in (0,L)\colon|f(s)|>\lambda\}|\leq t\},\ t\in(0,\infty).
\end{equation*}
Note that $f^*(t)=0$ for every $t\in[L,\infty)$. We say that functions $f,g\in\M(0,L)$ are \emph{equimeasurable}, and we write $f\sim g$, if $|\{s\in(0,L)\colon|f(s)|>\lambda\}|=|\{s\in (0,L)\colon|g(s)|>\lambda\}|$ for every $\lambda\in(0,\infty)$. We always have that $f\sim f^*$. The relation $\sim$ is transitive.

The \emph{maximal nonincreasing rearrangement} $f^{**} \colon (0,\infty) \to [0, \infty]$ of a function $f\in \M(0, L)$  is
defined as
\begin{equation*}
f^{**}(t)=\frac1t\int_0^ t f^{*}(s)\d{s},\ t\in(0,\infty).
\end{equation*}
The mapping $f \mapsto f^*$ is monotone in the sense that, for every $f,g\in\M(0,L)$,
\begin{equation*}
|f| \leq |g| \quad \text{a.e.~on $(0,L)$} \quad \Longrightarrow \quad f^* \leq g^* \quad  \text{on $(0, \infty)$};
\end{equation*}
consequently, the same implication remains true if ${}^*$ is replaced by ${}^{**}$. We have that
\begin{equation}\label{ch1:ri:twostarsdominateonestar}
f^*\leq f^{**}\quad\text{for every $f\in\M(0,L)$}.
\end{equation}
The operation $f\mapsto f^*$ is neither subadditive nor multiplicative. Although $f\mapsto f^*$ is not subadditive, the following pointwise inequality is valid for every $f,g\in\M_0(0,L)$ (\cite[Chapter~2, Proposition~1.7, (1.16)]{BS}):
\begin{equation}\label{ch1:ri:halfsubadditivityofonestar}
(f + g)^*(t) \leq f^*\Big(\frac{t}{2}\Big) + g^*\Big(\frac{t}{2}\Big) \quad \text{for every $t\in(0,L)$}.
\end{equation}
Furthermore, the lack of subadditivity of the operation of taking the nonincreasing rearrangement is, up to some extent, compensated by the following fact (\cite[Chapter~2,~(3.10)]{BS}):
\begin{equation}\label{ch1:ri:subadditivityofdoublestar}
\int_0^t(f+g)^*(s)\d{s}\leq\int_0^tf^*(s)\d{s}+\int_0^tg^*(s)\d{s}
\end{equation}
for every $t\in(0,\infty)$, $f,g\in\M_0(0,L)$. In other words, the operation $f\mapsto f^{**}$ is subadditive.

There are a large number of inequalities concerning rearrangements (e.g., \cite[Chapter~II, Section~2]{KPS:82}, \cite{HLP:52}). We state two of them, which shall prove particularly useful for us. The \emph{Hardy-Littlewood inequality} (\cite[Chapter~2, Theorem~2.2]{BS}) tells us that, for every $f,g\in\M(0,L)$,
\begin{equation}\label{ch1:ri:HL}
\int_0^L |f(t)||g(t)|\d{t}\leq\int_0^Lf^*(t)g^*(t)\d{t}.
\end{equation}
In particular, by taking $g=\chi_E$ in \eqref{ch1:ri:HL}, one obtain that
\begin{equation}\label{ch1:ri:HLg=chiE}
\int_E |f(t)|\d{t}\leq\int_0^{|E|}f^*(t)\d{t}
\end{equation}
for every measurable $E\subseteq (0,L)$. The \emph{Hardy lemma} (\cite[Chapter~2, Proposition~3.6]{BS}) states that, for every $f,g\in\Mpl(0,\infty)$ and every nonincreasing $h\in\Mpl(0,\infty)$,
\begin{equation}\label{ch1:ri:hardy-lemma}
\begin{aligned}
&\text{if}\quad\int_0^tf(s)\d{s}\leq \int_0^tg(s)\d{s}\quad\text{for all $t\in(0,\infty)$,}\\
&\text{then}\quad\int_0^\infty f(t)h(t)\d{t}\leq \int_0^\infty g(t)h(t)\d{t}.
\end{aligned}
\end{equation}

A functional $\varrho\colon\Mpl(0,L)\to[0,\infty]$ is called a \emph{rearrangement\hyp{}invariant function norm} (on $(0,L)$) if, for all $f$, $g$ and $\{f_k\}_{k=1}^\infty$ in $\Mpl(0,L)$, and every $\lambda\in[0,\infty)$:
\begin{itemize}
\item[(P1)] $\|f\|_{X(0,L)}=0$ if and only if $f=0$ a.e.~on $(0,L)$; $\|\lambda f\|_{X(0,L)}= \lambda \|f\|_{X(0,L)}$;  $\|f+g\|_{X(0,L)}\leq \|f\|_{X(0,L)} + \|g\|_{X(0,L)}$;
\item[(P2)] $\|f\|_{X(0,L)}\leq\|g\|_{X(0,L)}$ if $ f\leq g$ a.e.~on $(0,L)$;
\item[(P3)] $\|f_k\|_{X(0,L)} \nearrow \|f\|_{X(0,L)}$ if $f_k \nearrow f$ a.e.~on $(0,L)$;
\item[(P4)] $\|\chi_E\|_{X(0,L)}<\infty$ for every measurable $E\subseteq(0, L)$ of finite measure;
\item[(P5)] for every measurable $E\subseteq(0, L)$ of finite measure, there is a positive, finite constant $C_{E,X}$, possibly depending on $E$ and  $\|\cdot\|_{X(0,L)}$ but not on $f$, such that $\int_E f(t)\d{t} \leq C_{E,X} \|f\|_{X(0,L)}$;
\item[(P6)] $\|f\|_{X(0,L)} = \|g\|_{X(0,L)}$ whenever $f\sim g$.
\end{itemize}

The \textit{Hardy--Littlewood--P\'olya principle} (\cite[Chapter~2, Theorem~4.6]{BS}) asserts that, for every $f,g\in\M(0,L)$ and every rearrangement\hyp{}invariant function norm $\|\cdot\|_{X(0,L)}$,
\begin{equation}\label{ch1:ri:HLP}
\text{if $\int_0^tf^*(s)\d{s}\leq \int_0^tg^*(s)\d{s}$ for all $t\in(0,L)$, then $\|f\|_{X(0,L)}\leq \|g\|_{X(0,L)}$}.
\end{equation}

With every rearrangement\hyp{}invariant function norm $\|\cdot\|_{X(0,L)}$, we associate another functional $\|\cdot\|_{X'(0,L)}$ defined as
\begin{equation}\label{ch1:ri:normX'}
\|f\|_{X'(0,L)}= \sup_{\substack{g\in{\Mpl(0,L)}\\\|g\|_{X(0,L)}\leq1}}\int_0^Lf(t)g(t)\d{t},\ f\in\Mpl(0,L).
\end{equation}
The functional $\|\cdot\|_{X'(0,L)}$ is also a rearrangement\hyp{}invariant function norm (\cite[Chapter~2, Proposition~4.2]{BS}), and it is called the \emph{associate function norm} of $\|\cdot\|_{X(0,L)}$. Furthermore, we always have that (\cite[Chapter~1, Theorem~2.7]{BS})
\begin{equation}\label{ch1:ri:normX''}
\|f\|_{X(0,L)}= \sup_{\substack{g\in{\Mpl(0,L)}\\\|g\|_{X'(0,L)}\leq1}}\int_0^Lf(t)g(t)\d{t} \quad\text{for every $f\in\Mpl(0,L)$},
\end{equation}
that is,
\begin{equation}\label{ch1:ri:X''=X}
\|\cdot \|_{(X')'(0,L)} = \|\cdot \|_{X(0,L)}.
\end{equation}
Consequently, statements like ``Let $\|\cdot\|_{X(0,L)}$ be \emph{the} rearrangement\hyp{}invariant function norm whose associate function norm is \dots'' are well justified. The supremum in \eqref{ch1:ri:normX''} does not change when the functions involved are replaced with their nonincreasing rearrangements (\cite[Chapter~2, Proposition~4.2]{BS}), that is,
\begin{equation}\label{ch1:ri:normX''down}
\|f\|_{X(0,L)}= \sup_{\substack{g\in{\Mpl(0,L)}\\\|g\|_{X'(0,L)}\leq1}}\int_0^Lf^*(t)g^*(t)\d{t} \quad\text{for every $f\in\Mpl(0,L)$}.
\end{equation}

Given a rearrangement\hyp{}invariant function norm $\|\cdot\|_{X(0,L)}$, we extend it from $\Mpl(0,L)$ to $\M(0,L)$ by $\|f\|_{X(0,L)}=\|\,|f|\,\|_{X(0,L)}$. The extended functional $\|\cdot\|_{X(0,L)}$ restricted to the linear set $X(0,L)$ defined as
\begin{equation*}
X(0, L)=\{f\in\M(0,L)\colon \|f\|_{X(0,L)}<\infty\}
\end{equation*}
is a norm (provided that we identify any two functions from $\M(0,L)$ coinciding a.e.~on $(0,L)$, as usual). In fact, $X(0,L)$ endowed with the norm $\|\cdot\|_{X(0,L)}$ is a Banach space (\cite[Chapter~1, Theorem~1.6]{BS}). We say that $X(0,L)$ is a \emph{rearrangement\hyp{}invariant function space}. Note that $f\in\M(0,L)$ belongs to $X(0,L)$ if and only if $\|f\|_{X(0,L)}<\infty$. We always have that
\begin{equation}\label{ch1:ri:XembeddedinM0}
S(0,L)\subseteq X(0,L)\subseteq\M_0(0,L),
\end{equation}
where $S(0,L)$ denotes the set of all simple functions on $(0,L)$ (by a simple function, we mean a (finite) linear combination of characteristic functions of measurable sets having finite measure). Moreover, the second inclusion is continuous if the linear set $\M_0(0,L)$ is endowed with the (metrizable) topology of convergence in measure on sets of finite measure (\cite[Chapter~1, Theorem~1.4]{BS}).

The class of rearrangement\hyp{}invariant function spaces contains a large number of customary function spaces, such as Lebesgue spaces $L^p$ ($p\in[1,\infty]$), Lorentz spaces $L^{p,q}$ (e.g., \cite[pp.~216--220]{BS}), Orlicz spaces (e.g., \cite{RR:91}) or Lorentz--Zygmund spaces (e.g., \cite{BR:80, OP:99}), to name but a few. Here we provide definitions of only those rearrangement\hyp{}invariant function norms that we shall explicitly need. For $p\in[1,\infty]$, we define the Lebesgue function norm $\|\cdot\|_{L^p(0,L)}$ as
\begin{equation*}
\|f\|_{L^p(0,L)} = \begin{cases}
\int_0^L f(t)^p \d{t} \quad &\text{if $p\in[1,\infty)$},\\
\esssup_{t\in(0,L)} f(t) \quad &\text{if $p=\infty$},
\end{cases}
\end{equation*}
$f\in\Mpl(0,L)$. Given a measurable function $v\colon(0,L)\to(0,\infty)$ such that $V(t)<\infty$ for every $t\in(0,L)$, where $V(t)=\int_0^t v(s)\d{s}$, we define the functional $\|\cdot\|_{\Lambda^1_v(0,L)}$ as
\begin{equation*}
\|f\|_{\Lambda^1_v(0,L)} = \int_0^L f^*(s) v(s) \d{s},\ f\in\Mpl(0,L).
\end{equation*}
The functional is equivalent to a rearrangement\hyp{}invariant function norm if and only if (\cite[Theorem~2.3]{CGS:96}, see also \cite[Proposition~1]{S:93} with regard to local embedding of $\Lambda^1_v(0,L)$ in $L^1(0,L)$)
\begin{equation}\label{ch1:ri:when_is_Lambda1_ri_norm}
\frac{V(t)}{t}\lesssim\frac{V(s)}{s}\quad\text{for every $0<s< t<L$}.
\end{equation}
By the fact that it is equivalent to a rearrangement\hyp{}invariant function norm, we mean that there is a rearrangement\hyp{}invariant function norm $\varrho$ on $(0,L)$ such that $\|f\|_{\Lambda^1_v(0,L)} \approx \varrho(f)$ for every $f\in\Mpl(0,L)$; hence we can treat $\Lambda^1_v(0,L)$ as a rearrangement\hyp{}invariant function space whenever \eqref{ch1:ri:when_is_Lambda1_ri_norm} is satisfied. Let $\psi\colon (0,L) \to (0,\infty)$ be a \emph{quasiconcave function}, that is, a nondecreasing function such that the function $(0,L)\ni t\mapsto \frac{\psi(t)}{t}$ is nonincreasing. The functional $\|\cdot\|_{M_\psi(0,L)}$ defined as
\begin{equation*}
\|f\|_{M_\psi(0,L)} = \sup_{t\in(0,L)} \psi(t) f^{**}(t),\ t\in(0,L),
\end{equation*}
is a rearrangement\hyp{}invariant function norm (\cite[Proposition~7.10.2]{PKJF:13}).

The rearrangement\hyp{}invariant function space $X'(0,L)$ built upon the associate function norm $\|\cdot\|_{X'(0,L)}$ of a rearrangement\hyp{}invariant function norm $\|\cdot\|_{X(0,L)}$ is called the \emph{associate function space} of $X(0,L)$. Thanks to \eqref{ch1:ri:X''=X}, we have that $(X')'(0,L)=X(0,L)$. Furthermore, one has that
\begin{equation}\label{ch1:ri:holder}
\int_0^L |f(t)||g(t)|\d{t}\leq \|f\|_{X(0,L)}\|g\|_{X'(0,L)}\quad\text{for every $f,g\in\M(0,L)$}.
\end{equation}
Inequality \eqref{ch1:ri:holder} is a H\"older-type inequality, and we shall refer to it as the H\"older inequality.

Let $X(0,L)$ and $Y(0,L)$ be rearrangement\hyp{}invariant function spaces. We say that $X(0,L)$ is \emph{embedded in} $Y(0,L)$, and we write $X(0,L)\hookrightarrow Y(0,L)$, if there is a positive constant $C$ such that $\|f\|_{Y(0,L)}\leq C\|f\|_{X(0,L)}$ for every $f\in\M(0,L)$. If $X(0,L)\hookrightarrow Y(0,L)$ and $Y(0,L)\hookrightarrow X(0,L)$ simultaneously, we write that $X(0,L)=Y(0,L)$. We have that (\cite[Chapter~1, Theorem~1.8]{BS})
\begin{equation}\label{ch1:ri:inclusion_is_always_continuous}
X(0,L)\hookrightarrow Y(0,L)\quad\text{if and only if}\quad X(0,L)\subseteq Y(0,L).
\end{equation}
Furthermore,
\begin{equation}\label{ch1:ri:XtoYiffY'toX'}
X(0,L)\hookrightarrow Y(0,L)\quad\text{if and only if}\quad Y'(0,L)\hookrightarrow X'(0,L)
\end{equation}
with the same embedding constants.

If $\|\cdot\|_{X(0,L)}$ and $\|\cdot\|_{Y(0,L)}$ are rearrangement\hyp{}invariant function norms, then so are $\|\cdot\|_{X(0,L) \cap Y(0,L)}$ and $\|\cdot\|_{(X + Y)(0,L)}$ defined as
\begin{align*}
\|f\|_{X(0,L) \cap Y(0,L)} &= \max\{\|f\|_{X(0,L)}, \|f\|_{Y(0,L)}\},\ f\in\Mpl(0,L),
\intertext{and}
\|f\|_{(X + Y)(0,L)} &= \inf_{f=g+h}(\|g\|_{X(0,L)} + \|h\|_{Y(0,L)}),\ f\in\Mpl(0,L),
\end{align*}
where the infimum extends over all possible decompositions $f=g+h$, $g,h\in\Mpl(0,L)$. Furthermore, we have that (\cite[Theorem~3.1]{L:78}, also \cite[Lemma~1.12]{CNS:03})
\begin{equation}\label{ch1:ri:dual_sum_and_inter}
(X(0,L) \cap Y(0,L))' = (X' + Y')(0,L) \quad \text{and} \quad (X + Y)'(0,L) = X'(0,L) \cap Y'(0,L)
\end{equation}
with equality of norms. The \emph{$\K$-functional} between $X(0,L)$ and $Y(0,L)$ is, for every $f\in (X + Y)(0,L)$ and $t\in(0,\infty)$, defined as
\begin{equation*}
\K(f, t; X, Y) = \inf_{f=g+h}(\|g\|_{X(0,L)} + t\|h\|_{Y(0,L)}),
\end{equation*}
where the infimum extends over all possible decompositions $f=g+h$ with $g\in X(0,L)$ and $h\in Y(0,L)$. For every $f\in (X + Y)(0,L)\setminus\{0\}$, $\K(f, \cdot; X,Y)$ is a positive increasing concave function on $(0,\infty)$ (\cite[Chapter~5, Proposition~1.2]{BS}).

Let $X_0(0,L)$ and $X_1(0,L)$ be rearrangement\hyp{}invariant function spaces. We say that a rearrangement\hyp{}invariant function space $X(0,L)$ is an \emph{intermediate space} between $X_0(0,L)$ and $X_1(0,L)$ if $X_0(0,L) \cap X_1(0,L)\hookrightarrow X(0,L) \hookrightarrow (X_0 + X_1)(0,L)$. A linear operator $T$ defined on $(X_0 + X_1)(0,L)$ having values in $(X_0 + X_1)(0,L)$ is said to be \emph{admissible} for the couple $(X_0(0,L), X_1(0,L))$ if $T$ is bounded on both $X_0(0,L)$ and $X_1(0, L)$. An intermediate space $X(0,L)$ between $X_0(0,L)$ and $X_1(0,L)$ is an \emph{interpolation space} with respect to the couple $(X_0(0,L), X_1(0,L))$ if every admissible operator for the couple is bounded on $X(0,L)$. By \cite[Theorem~3]{C:66}, $X(0,L)$ is always an interpolation space with respect to the couple $(L^1(0,L), L^\infty(0,L))$.

We always have that (\cite[Chapter~2, Theorem~6.6]{BS})
\begin{equation*}
L^1(0,L)\cap L^\infty(0,L)\hookrightarrow X(0,L)\hookrightarrow L^1(0,L) + L^\infty(0,L).
\end{equation*}
In particular,
\begin{equation}\label{ch1:ri:smallestandlargestrispacefinitemeasure}
L^\infty(0,L)\hookrightarrow X(0,L)\hookrightarrow L^1(0,L)
\end{equation}
provided that $L<\infty$.

The \emph{dilation operator} is bounded on every rearrangement\hyp{}invariant function space $X(0,L)$. More precisely, we have that (\cite[Chapter~3, Proposition~5.11]{BS})
\begin{equation}\label{ch1:ri:dilation}
\|D_af\|_{X(0,L)}\leq\max\{1,a\}\|f\|_{X(0,L)}\quad\text{for every $f\in\M(0,L)$},
\end{equation}
where $D_af$ is defined on $(0,L)$ as
\begin{equation*}
D_af(t)=\begin{cases}

f(\frac{t}{a}), \quad &\text{if $L=\infty$},\\

f(\frac{t}{a})\chi_{(0, aL)}(t), \quad &\text{if $L<\infty$}.
\end{cases}
\end{equation*}

When $\|\cdot\|_{X(0,L)}$ is a rearrangement\hyp{}invariant function norm, we define its \emph{fundamental function} $\varphi_{X(0,L)}$ as
\begin{equation*}
\varphi_{X(0,L)}(t)=\|\chi_E\|_{X(0,L)},\ t\in[0, L),
\end{equation*}
where $E$ is any measurable subset of $(0,L)$ such that $|E|=t$. The fundamental function is well defined thanks to property (P6) of rearrangement\hyp{}invariant function norms and is a quasiconcave function. The fundamental functions of $\|\cdot\|_{X(0,L)}$ and $\|\cdot\|_{X'(0,L)}$ satisfy (\cite[Chapter~2, Theorem~5.2]{BS}) that
\begin{equation*}
\varphi_{X(0,L)}(t)\varphi_{X'(0,L)}(t)=t\quad\text{for every $t\in[0,L)$}.
\end{equation*}

\section{Hardy-type operators on r.i.~spaces}
We start with an easy but useful observation concerning the Hardy-type operators defined by \eqref{opRdef} and \eqref{opHdef}. Let $u,v\colon(0,L)\to(0,\infty)$ be measurable functions. Let $\nu\colon(0,L)\to(0,L)$ be an increasing bijection. The operators $R_{u,v,\nu}$ and $H_{u,v,\nu^{-1}}$, where $\nu^{-1}$ is the inverse function to $\nu$, are in a sense dual to each other. More precisely, by using the Fubini theorem, one can easily verify that
\begin{equation}\label{RaHdual}
\int_0^Lf(t)R_{u,v,\nu}g(t)\d{t}=\int_0^Lg(t)H_{u,v,\nu^{-1}}f(t)\d{t}\quad\text{for every $f,g\in\Mpl(0,L)$}.
\end{equation}

The validity of \eqref{RaHdual} has an unsurprising, well-known consequence, which we state here for future reference (see also \myref{Corollary}{cor:copsonrestrvsunrestr}).
\begin{proposition}\label{prop:RXtoYbddiffHY'toX'}
Let $\|\cdot\|_{X(0,L)}$, $\|\cdot\|_{Y(0,L)}$ be rearrangement\hyp{}invariant function norms. Let $u,v\colon(0,L)\to(0,\infty)$ be measurable. Let $\nu\colon(0,L)\to(0,L)$ be an increasing bijection. We have that
\begin{equation}\label{prop:RXtoYbddiffHY'toX':normRXtoZ=normHZ'toX'}
\sup_{\|f\|_{X(0,L)}\leq1}\|R_{u,v,\nu}f\|_{Y(0,L)}=\sup_{\|g\|_{Y'(0,L)}\leq1}\|H_{u,v,\nu^{-1}}g\|_{X'(0,L)}.
\end{equation}
In particular,
\begin{align}\label{prop:RXtoYbddiffHY'toX':RXtoZiffHZ'toX'}
\begin{split}
R_{u,v,\nu}\colon X(0,L)\to Y(0,L)\quad&\text{is bounded if and only if}\\
H_{u,v,\nu{-1}}\colon Y'(0,L)\to X'(0,L)\quad&\text{is bounded}.
\end{split}
\end{align}
\end{proposition}
\begin{proof}
We have that
\begin{align*}
\sup_{\|f\|_{X(0,L)}\leq1}\|R_{u,v,\nu}f\|_{Y(0,L)}&=\sup_{\|f\|_{X(0,L)}\leq1}\sup_{\|g\|_{Y'(0,L)}\leq1}\int_0^LR_{u,v,\nu}f(t)|g(t)|\d{t}\\
&=\sup_{\|f\|_{X(0,L)}\leq1}\sup_{\|g\|_{Y'(0,L)}\leq1}\int_0^L|f(t)|H_{u,v,\nu^{-1}}g(t)\d{t}\\
&=\sup_{\|g\|_{Y'(0,L)}\leq1}\|H_{u,v,\nu^{-1}}g\|_{X'(0,L)}
\end{align*}
thanks to \eqref{ch1:ri:normX''}, \eqref{RaHdual} and \eqref{ch1:ri:normX'}.
\end{proof}

\subsection{Optimal r.i.~function norms}\label{sec:optimal_norms}
In this section, we shall investigate optimal mapping properties of the operators $H_{u,v,\nu}$, $R_{u,v,\nu}$. Let $T$ be one of them. We say that a rearrangement\hyp{}invariant function space $Y(0,L)$ is \emph{the optimal target space} for the operator $T$ and a rearrangement\hyp{}invariant function space $X(0,L)$ if $T\colon X(0,L)\to Y(0,L)$ is bounded and $Y(0,L)\hookrightarrow Z(0,L)$ whenever $Z(0,L)$ is a rearrangement\hyp{}invariant function space such that $T\colon X(0,L)\to Z(0,L)$ is bounded (in other words, $\|\cdot\|_{Y(0,L)}$ is the strongest target rearrangement\hyp{}invariant function norm for $T$ and $\|\cdot\|_{X(0,L)}$). We say that a rearrangement\hyp{}invariant function space $X(0,L)$ is \emph{the optimal domain space} for the operator $T$ and a rearrangement\hyp{}invariant function space $Y(0,L)$ if $T\colon X(0,L)\to Y(0,L)$ is bounded and $Z(0,L)\hookrightarrow X(0,L)$ whenever $Z(0,L)$ is a rearrangement\hyp{}invariant function space such that $T\colon Z(0,L)\to Y(0,L)$ is bounded (in other words, $\|\cdot\|_{X(0,L)}$ is the weakest domain rearrangement\hyp{}invariant function norm for $T$ and $\|\cdot\|_{Y(0,L)}$).

We start by characterizing when the functional $\Mpl(0,L)\ni f\mapsto \|R_{u,v,\nu}(f^*)\|_{X(0,L)}$ is a rearrangement\hyp{}invariant function norm. It turns out that it also enables us to characterize optimal domain and target spaces for $R_{u,v,\nu}$ and $H_{u,v,\nu}$, respectively.
\begin{proposition}\label{prop:norminducedbyR}
Let $\|\cdot\|_{X(0,L)}$ be a rearrangement\hyp{}invariant function norm. Let $u\colon(0,L)\to(0,\infty)$ be a nondegenerate nonincreasing function. If $L<\infty$, assume that $u(L^-)>0$. Let $v\colon(0,L)\to(0,\infty)$ be measurable. Let $\nu\colon(0,L)\to(0,L)$ be an increasing bijection. Set
\begin{equation*}
\varrho(f)=\Bigg\|v(t)\int_0^{\nu(t)}f^*(s)u(s)\d{s}\Bigg\|_{X(0,L)},\ f\in\Mpl(0,L).
\end{equation*}
and
\begin{equation}\label{prop:Roptimaldomain:xidef}
\xi(t)=\begin{cases}
v(t)U(\nu(t)),\ t\in(0,L),\quad&\text{if $L<\infty$,}\\
v(t)U(\nu(t))\chi_{(0,1)}(t)+v(t)\chi_{(1,\infty)}(t),\ t\in(0, \infty),\quad&\text{if $L=\infty$.}
\end{cases}
\end{equation}
The functional $\varrho$ is a rearrangement\hyp{}invariant function norm if and only if $\xi\in X(0,L)$.

If $\xi\in X(0,L)$, then the rearrangement\hyp{}invariant function space induced by $\varrho$ is the optimal domain space for the operator $R_{u,v,\nu}$ and $X(0,L)$. If $\xi\not\in X(0,L)$, then there is no rearrangement\hyp{}invariant function space $Z(0,L)$ such that $R_{u,v,\nu}\colon Z(0,L)\to X(0,L)$ is bounded.
\end{proposition}
\begin{proof}
We shall show that $\varrho$ is a rearrangement\hyp{}invariant function norm provided that $\xi\in X(0,L)$. Before we do that, note that, since $u$ is positive and nonincreasing, its nondegeneracy implies that $0 < U(t) < \infty$ for every $t\in(0,L]\cap\R$.

\emph{Property \emph{(P1)}.} The positive homogeneity and positive definiteness of $\varrho$ can be readily verified. As for the subadditivity of $\varrho$, it follows from \eqref{ch1:ri:subadditivityofdoublestar} combined with Hardy's lemma \eqref{ch1:ri:hardy-lemma} that
\begin{equation*}
\int_0^L(f+g)^*(s)u(s)\chi_{(0,\nu(t))}(s)\d{s}\leq\int_0^Lf^*(s)u(s)\chi_{(0,\nu(t))}(s)\d{s}+\int_0^Lg^*(s)u(s)\chi_{(0,\nu(t))}(s)\d{s}
\end{equation*}
for every $f,g\in\Mpl(0,L)$ and $t\in(0,L)$ thanks to the fact that $u$ is nonincreasing. Since $\|\cdot\|_{X(0,L)}$ is subadditive, it follows that
\begin{equation*}
\varrho(f+g)\leq\varrho(f)+\varrho(g)\quad\text{for every $f,g\in\Mpl(0,L)$}.
\end{equation*}

\emph{Properties \emph{(P2)} and \emph{(P3)}.} Since $\|\cdot\|_{X(0,L)}$ has these properties, it can be readily verified that $\varrho$, too, has them.

\emph{Property \emph{(P4)}.} First, assume that $L<\infty$. Clearly, $\varrho(\chi_{(0,L)})<\infty$ if and only if $v(t)U(\nu(t))\in X(0,L)$. Since $\varrho$ has property (P2), $\varrho$ has property (P4) if and only if $v(t)U(\nu(t))\in X(0,L)$. Second, assume that $L=\infty$. Let $E\subseteq(0,\infty)$ be a set of finite positive measure. Clearly, $\varrho(\chi_E)<\infty$ if and only if $v(t)U(\nu(t))\chi_{(0,|E|)}(t)+v(t)\chi_{(|E|,\infty)}(t)\in X(0,\infty)$. If $|E|\leq1$, then 
\begin{align*}
&\|v(t)U(\nu(t))\chi_{(0,|E|)}(t)+v(t)\chi_{(|E|,\infty)}(t)\|_{X(0,\infty)}\\
&\leq\|v(t)U(\nu(t))\chi_{(0,1)}(t)\|_{X(0,\infty)}+\|v(t)\chi_{(|E|,1)}(t)\|_{X(0,\infty)}+\|v(t)\chi_{(1,\infty)}(t)\|_{X(0,\infty)}\\
&\leq\|v(t)U(\nu(t))\chi_{(0,1)}(t)\|_{X(0,\infty)}+\frac1{U(\nu(|E|))}\|U(\nu(t))v(t)\chi_{(|E|,1)}(t)\|_{X(0,\infty)}\\
&\quad+\|v(t)\chi_{(1,\infty)}(t)\|_{X(0,\infty)}\\
&\leq\Big(1+\frac1{U(\nu(|E|))}\Big)\|v(t)U(\nu(t))\chi_{(0,1)}(t)\|_{X(0,\infty)}+\|v(t)\chi_{(1,\infty)}(t)\|_{X(0,\infty)}.
\end{align*}
If $E\geq1$, we can obtain, in a similar way, that
\begin{align*}
&\|v(t)U(\nu(t))\chi_{(0,|E|)}(t)+v(t)\chi_{(|E|,\infty)}(t)\|_{X(0,\infty)}\\
&\leq\|v(t)U(\nu(t))\chi_{(0,1)}(t)\|_{X(0,\infty)}+(1+U(\nu(|E|)))\|v(t)\chi_{(1,\infty)}(t)\|_{X(0,\infty)}.
\end{align*}
Either way, we have that $\varrho(\chi_E)<\infty$ if and only if
\begin{equation*}
v(t)U(\nu(t))\chi_{(0,1)}(t)+v(t)\chi_{(1,\infty)}(t)\in X(0,\infty).
\end{equation*}

\emph{Property \emph{(P5)}.}
Let $E\subseteq(0,L)$ be a set of finite positive measure. Let $f\in\Mpl(0,L)$. Note that the function $(0,L)\ni t\mapsto\frac1{U(\nu(t))}\int_0^{\nu(t)}f^*(s)u(s)\d{s}$ is nonincreasing because it is the integral mean of a nonnegative nonincreasing function over the interval $(0,\nu(t))$ with respect to the measure $u(s)\d{s}$. Thanks to that and the monotonicity of $u$, we obtain that
\begin{align*}
\Big\|v(t)\int_0^{\nu(t)}f^*(s)u(s)\d{s}\Big\|_{X(0,L)}&\geq\Big\|v(t)\chi_{(0,\nu^{-1}(|E|))}(t)\int_0^{\nu(t)}f^*(s)u(s)\d{s}\Big\|_{X(0,L)}\\
&\geq\Big\|v(t)U(\nu(t))\chi_{(0,\nu^{-1}(|E|))}(t)\Big\|_{X(0,L)}\frac1{U(|E|)}\int_0^{|E|}f^*(s)u(s)\d{s}\\
&\geq\Big\|v(t)U(\nu(t))\chi_{(0,\nu^{-1}(|E|))}(t)\Big\|_{X(0,L)}\frac{u(|E|^-)}{U(|E|)}\int_0^{|E|}f^*(s)\d{s}\\
&\geq\Big\|v(t)U(\nu(t))\chi_{(0,\nu^{-1}(|E|))}(t)\Big\|_{X(0,L)}\frac{u(|E|^-)}{U(|E|)}\int_Ef(s)\d{s},
\end{align*}
where we used \eqref{ch1:ri:HLg=chiE} in the last inequality.

\emph{Property \emph{(P6)}.} Since $f^*=g^*$ when $f,g\in\Mpl(0,L)$ are equimeasurable, this is obvious.

Note that the necessity of $\xi\in X(0,L)$ for $\varrho$ to be a rearrangement\hyp{}invariant function norm was already proved in the paragraph devoted to property (P4).

Assume now that $\xi\in X(0,L)$ and denote by $Y(0,L)$ the rearrangement\hyp{}invariant function space induced by $\varrho$. Thanks to the Hardy--Littlewood inequality \eqref{ch1:ri:HL} and the monotonicity of $u$, we have that 
\begin{equation*}
\|R_{u,v,\nu}f\|_{X(0,L)}\leq\|R_{u,v,\nu}(f^*)\|_{X(0,L)}=\|f\|_{Y(0,L)}\quad\text{for every $f\in\Mpl(0,L)$}.
\end{equation*}
Hence $R_{u,v,\nu}\colon Y(0,L)\to X(0,L)$ is bounded. Next, if $Z(0,L)$ is a rearrangement\hyp{}invariant function space such that $R_{u,v,\nu}\colon Z(0,L)\to X(0,L)$ is bounded, then we have that
\begin{equation*}
\|f\|_{Y(0,L)}=\|R_{u,v,\nu}(f^*)\|_{X(0,L)}\lesssim\|f^*\|_{Z(0,L)}=\|f\|_{Z(0,L)}\quad\text{for every $f\in\Mpl(0,L)$},
\end{equation*}
and so $Z(0,L)\hookrightarrow Y(0,L)$. Finally, note that, if $R_{u,v,\nu}\colon Z(0,L)\to X(0,L)$ is bounded, then
\begin{equation*}
\|\xi\|_{X(0,L)}\approx\|R_{u,v,\nu}(\chi_{(0,a)})\|_{X(0,L)}\lesssim\|\chi_{(0,a)}\|_{Z(0,L)}<\infty,
\end{equation*}
where 
\begin{equation}\label{prop:Roptimaldomain:L_or_1}
a=\begin{cases}

L\quad&\text{if $L<\infty$},\\
1\quad&\text{if $L=\infty$};

\end{cases}
\end{equation}
hence $\xi\in X(0,L)$.
\end{proof}

\begin{remark}\label{rem:optimal_for_R_iff_for_H}
Thanks to \eqref{prop:RXtoYbddiffHY'toX':RXtoZiffHZ'toX'} and \eqref{ch1:ri:X''=X}, $Y(0,L)$ is the optimal target space for the operator $H_{u,v,\nu}$ and $X(0,L)$ if and only if $Y'(0,L)$ is the optimal domain space for the operator $R_{u,v,\nu^{-1}}$ and $X'(0,L)$. Therefore, \myref{Proposition}{prop:norminducedbyR} actually also characterizes optimal target spaces for the operator $H_{u,v,\nu}$.
\end{remark}

\begin{proposition}\label{prop:Hoptimalrange}
Let $\|\cdot\|_{X(0,L)}$ be a rearrangement\hyp{}invariant function norm. Let $\nu\colon(0,L)\to(0,L)$ be an increasing bijection. Let $u\colon(0,L)\to(0,\infty)$ be a nondegenerate nonincreasing function. If $L<\infty$, assume that $u(L^-)>0$. Let $v\colon(0,L)\to(0,\infty)$ be measurable. Assume that $\xi\in X'(0,L)$, where $\xi$ is defined by \eqref{prop:Roptimaldomain:xidef} with $\nu$ replaced by $\nu^{-1}$. Let $\|\cdot\|_{Y(0,L)}$ be the rearrangement\hyp{}invariant function norm whose associate function norm $\|\cdot\|_{Y'(0,L)}$ is defined as
\begin{equation*}
\|f\|_{Y'(0,L)}=\Bigg\|v(t)\int_0^{\nu^{-1}(t)}f^*(s)u(s)\d{s}\Bigg\|_{X'(0,L)},\ f\in\Mpl(0,L).
\end{equation*}
The rearrangement\hyp{}invariant function space $Y(0,L)$ is the optimal target space for the operator $H_{u,v,\nu}$ and $X(0,L)$. Moreover, if $\xi\not\in X'(0,L)$, then there is no rearrangement\hyp{}invariant function space $Z(0,L)$ such that $H_{u,v,\nu}\colon X(0,L)\to Z(0,L)$ is bounded.
\end{proposition}

We now turn our attention to $H_{u,v,\nu}$; it turns out that the situation becomes significantly more complicated. Notably the fact that, unlike with $R_{u,v,\nu}$, the integration is carried out over intervals away from $0$ often causes great difficulties. In particular, the functional $\Mpl(0,L)\ni f\mapsto \|H_{u,v,\nu}(f^*)\|_{X(0,L)}$ is hardly ever subadditive. Instead, in general, we need to consider a more complicated functional (see \myref{Proposition}{prop:norminducedbyHwhenRnonincreasing}, however). We will need to impose a mild condition on $\nu$. We write $\nu\in\Dinf{0}$, $\nu\in\Dinf{\infty}$, $\nu\in\Dsup{0}$ and $\nu\in\Dsup{\infty}$ if there is $\theta>1$ such that $\liminf_{t\to0^+}\frac{\nu(\theta t)}{\nu(t)}>1$, $\liminf_{t\to\infty}\frac{\nu(\theta t)}{\nu(t)}>1$, $\limsup_{t\to0^+}\frac{\nu(\theta t)}{\nu(t)}<\infty$ and $\limsup_{t\to\infty}\frac{\nu(\theta t)}{\nu(t)}<\infty$, respectively. When we need to emphasize the exact value of $\theta$, we will write $\nu\in\Dinf[\theta]{0}$ and so forth.
\begin{proposition}\label{prop:norminducedbyH}
Let $\|\cdot\|_{X(0,L)}$ be a rearrangement\hyp{}invariant function norm. Let $u,v\colon(0,L)\to(0,\infty)$ be nonincreasing. If $L<\infty$, assume that $v(L^-)>0$. Let $\nu\colon(0,L)\to(0,L)$ be an increasing bijection. If $L=\infty$, assume that $\nu\in\Dinf{\infty}$. Set
\begin{equation}\label{prop:norminducedbyH:normdef}
\varrho(f)=\sup_{h\sim f}\Big\|u(t)\int_{\nu(t)}^L h(s)v(s)\d{s}\Big\|_{X(0,L)},\ f\in\Mpl(0,L),
\end{equation}
where the supremum extends over all $h\in\Mpl(0,L)$ equimeasurable with $f$. The functional $\varrho$ is a rearrangement\hyp{}invariant function norm if and only if
\begin{equation}\label{prop:norminducedbyH:assum}
\begin{cases}
u(t)\int_{\nu(t)}^L v(s) \d{s}\in X(0,L)\quad&\text{if $L<\infty$},\\
\begin{gathered}[b]u(t)\chi_{(0,\nu^{-1}(1))}(t)\int_{\nu(t)}^1v(s)\d{s}\in X(0,\infty)\ \text{and}\\ \limsup_{\tau\to\infty}v(\tau)\|u\chi_{(0,\nu^{-1}(\tau))}\|_{X(0,\infty)}<\infty\end{gathered}\quad&\text{if $L=\infty$}.
\end{cases}
\end{equation}

If \eqref{prop:norminducedbyH:assum} is satisfied, then the rearrangement\hyp{}invariant function space induced by $\varrho$ is the optimal domain space for the operator $H_{u,v,\nu}$ and $X(0,L)$. If \eqref{prop:norminducedbyH:assum} is not satisfied, then there is no rearrangement\hyp{}invariant function space $Z(0,L)$ such that $H_{u,v,\nu}\colon Z(0,L)\to X(0,L)$ is bounded.
\end{proposition}
\begin{proof}
We shall show that $\varrho$ is a rearrangement\hyp{}invariant function norm provided that \eqref{prop:norminducedbyH:assum} is satisfied.

\emph{Property \emph{(P2)}.} Let $f,g\in\Mpl(0,L)$ be such that $f\leq g$ a.e. Consequently, $f^*\leq g^*$. Suppose that $\varrho(f) > \varrho(g)$. It implies that there is $\widetilde{f}\in\Mpl(0,L)$, $\widetilde{f}\sim f$, such that
\begin{equation}\label{prop:norminducedbyH:eq1}
\sup_{h\sim g}\Big\|u(t)\int_{\nu(t)}^L h(s)v(s)\d{s}\Big\|_{X(0,L)}<\Big\|u(t)\int_{\nu(t)}^L \widetilde{f}(s)v(s)\d{s}\Big\|_{X(0,L)}.
\end{equation}
When $L=\infty$, we may assume that $\lim_{t\to\infty}(\widetilde{f})^*(t)=\lim_{t\to\infty}f^*(t)=0$, for we would otherwise approximate $\widetilde f$ by functions $f_n=\widetilde{f}\chi_{(0,n)}$, $n\in\N$ (the monotone convergence theorem and property (P3) of $\|\cdot\|_{X(0,L)}$ would guarantee that the inequality above holds with $\widetilde{f}$ replaced by $f_n$ for $n$ large enough). Thanks to \cite[Chapter~2, Corollary~7.6]{BS} (also~\citep[Proposition~3]{R:70}), there is a measure-preserving transformation (in the sense of \cite[Chapter~2, Definition~7.1]{BS}) $\sigma\colon(0,L)\to(0,L)$ such that $\widetilde{f}=f^*\circ\sigma$. Since $\sigma$ is measure preserving, we have that $(g^*\circ\sigma)\sim g^*\sim g$ (\cite[Chapter~2, Proposition~7.2]{BS}). Consequently,
\begin{equation}\label{prop:norminducedbyH:eq2}
\begin{aligned}
\sup_{h\sim g}\Big\|u(t)\int_{\nu(t)}^L h(s)v(s)\d{s}\Big\|_{X(0,L)}&\geq\Big\|u(t)\int_{\nu(t)}^L g^*(\sigma(s))v(s)\d{s}\Big\|_{X(0,L)}\\
&\geq\Big\|u(t)\int_{\nu(t)}^L f^*(\sigma(s))v(s)\d{s}\Big\|_{X(0,L)}\\
&=\Big\|u(t)\int_{\nu(t)}^L \widetilde{f}(s)v(s)\d{s}\Big\|_{X(0,L)}.
\end{aligned}
\end{equation}
By combining \eqref{prop:norminducedbyH:eq1} and \eqref{prop:norminducedbyH:eq2}, we reach a contradiction. Hence $\varrho(f)\leq\varrho(g)$.

\emph{Property \emph{(P3)}.} Let $f, f_k\in\Mpl(0,L)$, $k\in \N$, be such that $f_k \nearrow f$ a.e. Thanks to property (P2) of $\varrho$, the limit $\lim_{k\to\infty}\varrho(f_k)$ exists and we clearly have that $\lim_{k\to\infty}\varrho(f_k)\leq\varrho(f)$. The fact that $\lim_{k\to\infty}\varrho(f_k)=\varrho(f)$ can be proved by contradiction in a similar way to the proof of (P2).

\emph{Property \emph{(P1)}.} The positive homogeneity and positive definiteness of $\varrho$ can be readily verified. As for the subadditivity of $\varrho$, let $f,g\in\Mpl(0,L)$ be simple functions. Let $h\in\Mpl(0,L)$ be such that $h\sim f+g$. Being equimeasurable with $f+g$, $h$ is a simple function having the same range as $f+g$. Furthermore, it is easy to see that $h$ can be decomposed as $h=h_1+h_2$, where $h_1,h_2\in\Mpl(0,L)$ are simple functions such that $h_1\sim f$ and $h_2\sim g$. Using the subadditivity of $\|\cdot\|_{X(0,L)}$, we obtain that
\begin{align*}
\Big\|u(t)\int_{\nu(t)}^L h(s)v(s)\d{s}\Big\|_{X(0,L)}&\leq\Big\|u(t)\int_{\nu(t)}^L h_1(s)v(s)\d{s}\Big\|_{X(0,L)}+\Big\|u(t)\int_{\nu(t)}^L h_2(s)v(s)\d{s}\Big\|_{X(0,L)}\\
&\leq\varrho(f)+\varrho(g).
\end{align*}
Hence $\varrho(f+g)\leq\varrho(f)+\varrho(g)$. When $f,g\in\Mpl(0,L)$  are general functions, we approximate each of them by a nondecreasing sequence of nonnegative, simple functions and use property (P3) of $\varrho$ to get $\varrho(f+g)\leq\varrho(f)+\varrho(g)$.

\emph{Property \emph{(P4)}.} Assume that $L<\infty$. Since $\varrho$ has property (P2), $\varrho$ has property (P4) if and only if $\varrho(\chi_{(0,L)})<\infty$. If $h\in\Mpl(0,L)$ is equimeasurable with $\chi_{(0,L)}$, then $h=1$ a.e.~on $(0,L)$; therefore,
\begin{equation*}
\varrho(\chi_{(0,L)})=\Big\|u(t)\int_{\nu(t)}^L v(s)\d{s}\Big\|_{X(0,L)}.
\end{equation*}
Hence $\varrho$ has property (P4) if and only if $u(t)\int_{\nu(t)}^L v(s)\d{s}\in X(0,L)$. Assume now that $L=\infty$. Fix $\theta>1$ such that $\nu\in\Dinf[\theta]{\infty}$. Let $E\subseteq(0,\infty)$ be of finite measure. Set $b=\max\Big\{1,\nu(1),\frac{\theta|E|}{M-1}\Big\}$, where $M=\inf_{t\in[1,\infty)}\frac{\nu(\theta t)}{\nu(t)}$. Note that $M>1$. Let $h\in\Mpl(0,\infty)$ be equimeasurable with $\chi_E$. It is easy to see that $h=\chi_F$ for some measurable $F\subseteq(0,\infty)$ such that $|F|=|E|$. Thanks to the (outer) regularity of the Lebesgue measure, there is an open set $G\supseteq F$ such that $|G|\leq \theta|F|$. Being an open set on the real line, $G\cap(b,\infty)$ can be expressed as $G\cap(b,\infty)=\bigcup_{k}(a_k,b_k)$, where $\{(a_k,b_k)\}_k$ is a countable system of mutually disjoint open intervals. We plainly have that $F\subseteq(0,b]\cup\big(G\cap(b,\infty)\big)$ and $a_k>b$. Furthermore, we have that $b_k-a_k\leq \theta|F|\leq (M-1)b < (M-1)a_k$, whence
\begin{equation}\label{prop:norminducedbyH:eq3}
\nu^{-1}(b_k)-\nu^{-1}(a_k)<(\theta-1)\nu^{-1}(a_k).
\end{equation}
We have that
\begin{align}
\Big\|u(t)\int_{\nu(t)}^\infty \chi_F(s)v(s)\d{s}\Big\|_{X(0,\infty)}&\leq\Big\|u(t)\int_{\nu(t)}^\infty \big(\chi_{(0,b]}(s)+\sum_k\chi_{(a_k,b_k)}(s)\big)v(s)\d{s}\Big\|_{X(0,\infty)}\notag\\
\begin{split}\label{prop:norminducedbyH:eq4}
&\leq\Big\|u(t)\chi_{(0,\nu^{-1}(b))}(t)\int_{\nu(t)}^bv(s)\d{s}\Big\|_{X(0,\infty)}\\
&\quad+\sum_k\Big\|u(t)\chi_{(0,\nu^{-1}(a_k))}(t)\int_{a_k}^{b_k}v(s)\d{s}\Big\|_{X(0,\infty)}\\
&\quad+\sum_k\Big\|u(t)\chi_{(\nu^{-1}(a_k),\nu^{-1}(b_k))}(t)\int_{\nu(t)}^{b_k}v(s)\d{s}\Big\|_{X(0,\infty)}.
\end{split}
\end{align}
Note that the assumption
\begin{equation}\label{prop:norminducedbyH:infinitecase:assump1}
\Big\|u(t)\chi_{(0,\nu^{-1}(1))}(t)\int_{\nu(t)}^1v(s)\d{s}\Big\|_{X(0,\infty)}<\infty
\end{equation}
together with the monotonicity of $u$ and $v$ implies that
\begin{equation}\label{prop:norminducedbyH:eq11}
\|u\chi_{(0,a)}\|_{X(0,\infty)}<\infty\quad\text{for every $a\in(0,\infty)$}.
\end{equation}
Indeed, since $u$ is nonincreasing, it is sufficient to show that $\|u\chi_{(0,\nu^{-1}(\frac1{2}))}\|_{X(0,\infty)}<\infty$, which follows from
\begin{align*}
\infty>\Big\|u(t)\chi_{(0,\nu^{-1}(1))}(t)\int_{\nu(t)}^1v(s)\d{s}\Big\|_{X(0,\infty)}&\geq\Big\|u(t)\chi_{(0,\nu^{-1}(\frac1{2}))}(t)\int_{\nu(t)}^1v(s)\d{s}\Big\|_{X(0,\infty)}\\
&\geq \frac{v(1)}{2}\Big\|u\chi_{(0,\nu^{-1}(\frac1{2}))}\Big\|_{X(0,\infty)}.
\end{align*}
Furthermore, note that \eqref{prop:norminducedbyH:eq11} guarantees that
\begin{equation}\label{prop:norminducedbyH:infinitecase:assump2}
\limsup_{\tau\to\infty}v(\tau)\|u\chi_{(0,\nu^{-1}(\tau))}\|_{X(0,\infty)}<\infty
\end{equation}
if and only if
\begin{equation}\label{prop:norminducedbyH:eq12}
\sup_{\tau\in[1,\infty)}v(\tau)\|u\chi_{(0,\nu^{-1}(\tau))}\|_{X(0,\infty)}<\infty.
\end{equation}
Now, as for the first term on the right-hand side of \eqref{prop:norminducedbyH:eq4}, we have that
\begin{align}\label{prop:norminducedbyH:eq7}
\begin{split}
\Big\|u(t)\chi_{(0,\nu^{-1}(b))}(t)\int_{\nu(t)}^bv(s)\d{s}\Big\|_{X(0,\infty)}&\leq\Big\|u(t)\chi_{(0,\nu^{-1}(1))}(t)\int_{\nu(t)}^1v(s)\d{s}\Big\|_{X(0,\infty)}\\
&\quad+\Big\|u(t)\chi_{(0,\nu^{-1}(1))}(t)\int_1^bv(s)\d{s}\Big\|_{X(0,\infty)}\\
&\quad+\Big\|u(t)\chi_{(\nu^{-1}(1),\nu^{-1}(b))}(t)\int_{\nu(t)}^bv(s)\d{s}\Big\|_{X(0,\infty)}\\
&\leq A<\infty,
\end{split}
\end{align}
where
\begin{align*}
A&=\Big\|u(t)\chi_{(0,\nu^{-1}(1))}(t)\int_{\nu(t)}^1v(s)\d{s}\Big\|_{X(0,\infty)}+v(1)(b-1)\|u\chi_{(0,\nu^{-1}(1))}\|_{X(0,\infty)}\\
&\quad+v(1)(b-1)\|u\chi_{(0,\nu^{-1}(b)-\nu^{-1}(1))}\|_{X(0,\infty)}.
\end{align*}
As for the second term on the right-hand side of \eqref{prop:norminducedbyH:eq4}, we have that
\begin{equation}\label{prop:norminducedbyH:eq5}
\begin{aligned}
\Big\|u(t)\chi_{(0,\nu^{-1}(a_k))}(t)\int_{a_k}^{b_k}v(s)\d{s}\Big\|_{X(0,\infty)}&\leq v(a_k)(b_k-a_k)\|u\chi_{(0,\nu^{-1}(a_k))}\|_{X(0,\infty)}\\
&\leq B(b_k-a_k),
\end{aligned}
\end{equation}
where $B$ is the supremum in \eqref{prop:norminducedbyH:eq12}, which is independent of $k$. Next,
\begin{align}
\Big\|u(t)\chi_{(\nu^{-1}(a_k),\nu^{-1}(b_k))}(t)\int_{\nu(t)}^{b_k}v(s)\d{s}\Big\|_{X(0,\infty)}&\leq\int_{a_k}^{b_k}v(s)\d{s}\|u\chi_{(\nu^{-1}(a_k),\nu^{-1}(b_k))}\|_{X(0,\infty)}\notag\\
\begin{split}\label{prop:norminducedbyH:eq6}
&\leq v(a_k)(b_k-a_k)\|u\chi_{(\nu^{-1}(a_k),\nu^{-1}(b_k))}\|_{X(0,\infty)}\\
&\leq v(a_k)(b_k-a_k)\|u\chi_{(0,(\theta-1)\nu^{-1}(a_k))}\|_{X(0,\infty)}\\
&\leq \ceil{\theta-1} v(a_k)(b_k-a_k)\|u\chi_{(0,\nu^{-1}(a_k))}\|_{X(0,\infty)}\\
&\leq \ceil{\theta-1} B(b_k-a_k),
\end{split}
\end{align}
where we used the monotonicity of $u$ and $v$ in the second inequality, \eqref{prop:norminducedbyH:eq3} in the third one, and the monotonicity of $u$ in the fourth one. By combining \eqref{prop:norminducedbyH:eq4} with \eqref{prop:norminducedbyH:eq7}, \eqref{prop:norminducedbyH:eq5} and \eqref{prop:norminducedbyH:eq6}, we obtain that
\begin{equation*}
\Big\|u(t)\int_{\nu(t)}^\infty h(s)v(s)\d{s}\Big\|_{X(0,\infty)}\leq A+\ceil{\theta}B\sum_k(b_k-a_k)\leq A+\ceil{\theta}\theta B|E|<\infty.
\end{equation*}
Hence $\varrho(\chi_E)<\infty$ provided that \eqref{prop:norminducedbyH:infinitecase:assump1} and \eqref{prop:norminducedbyH:infinitecase:assump2} are satisfied. The necessity of \eqref{prop:norminducedbyH:infinitecase:assump1} is obvious because we have that
\begin{equation*}
\Big\|u(t)\chi_{(0,\nu^{-1}(1))}(t)\int_{\nu(t)}^1v(s)\d{s}\Big\|_{X(0,\infty)}\leq\varrho(\chi_{(0,1)}).
\end{equation*}
As for the necessity of \eqref{prop:norminducedbyH:infinitecase:assump2}, suppose that $\limsup_{\tau\to\infty}v(\tau)\big\|u\chi_{(0,\nu^{-1}(\tau))}\big\|_{X(0,\infty)}=\infty$. It follows that there is a sequence $\tau_k\nearrow\infty$, $k\to\infty$, such that
\begin{equation*}
\lim_{k\to\infty}v(\tau_k)\|u\chi_{(0,\nu^{-1}(\tau_k))}\|_{X(0,\infty)}=\infty.
\end{equation*}
Since $\inf_{t\in[1,\infty)}\frac{\nu(\theta t)}{\nu(t)}>1$, we can find an $\varepsilon>0$ such that $\frac{\nu(\theta t)}{\nu(t)}\geq 1+\varepsilon$ for every $t\in[1,\infty)$. Moreover, we may clearly assume that $\tau_k\geq \nu(1)+1$ and $\frac{\tau_k}{\tau_k-1}\leq 1+\varepsilon$; hence
\begin{equation}\label{prop:norminducedbyH:eq10}
\nu^{-1}(\tau_k)-\nu^{-1}(\tau_k-1)\leq(\theta-1)\nu^{-1}(\tau_k-1)
\end{equation}
inasmuch as $\frac{\nu(\theta\nu^{-1}(\tau_k-1))}{\nu(\nu^{-1}(\tau_k-1))}\geq1+\varepsilon$. Using \eqref{prop:norminducedbyH:eq10} and the fact that $u$ is nonincreasing, we obtain that
\begin{align*}
\|u\chi_{(0, \nu^{-1}(\tau_k))}\|_{X(0,\infty)}&\leq\|u\chi_{(0, \nu^{-1}(\tau_k-1))}\|_{X(0,\infty)}+\|u\chi_{(\nu^{-1}(\tau_k-1),\nu^{-1}(\tau_k))}\|_{X(0,\infty)}\\
&\leq\|u\chi_{(0, \nu^{-1}(\tau_k-1))}\|_{X(0,\infty)}+\|u\chi_{(0,\nu^{-1}(\tau_k)-\nu^{-1}(\tau_k-1))}\|_{X(0,\infty)}\\
&\leq\|u\chi_{(0, \nu^{-1}(\tau_k-1))}\|_{X(0,\infty)}+\|u\chi_{(0, (\theta - 1)\nu^{-1}(\tau_k-1))}\|_{X(0,\infty)}\\
&=\ceil{\theta}\|u\chi_{(0, \nu^{-1}(\tau_k-1))}\|_{X(0,\infty)}.
\end{align*}
Therefore,
\begin{align*}
\varrho(\chi_{(0,1)})&\geq\Big\|u(t)\chi_{(0, \nu^{-1}(\tau_k-1))}(t)\int_{\nu(t)}^\infty\chi_{(\tau_k-1,\tau_k)}(s)v(s)\d{s}\Big\|_{X(0,\infty)}\\
&\geq v(\tau_k)\|u\chi_{(0, \nu^{-1}(\tau_k-1))}\|_{X(0,\infty)}\geq\frac1{\ceil{\theta}}v(\tau_k)\|u\chi_{(0, \nu^{-1}(\tau_k))}\|_{X(0,\infty)},
\end{align*}
which tends to $\infty$ as $k\to\infty$. Hence $\varrho(\chi_{(0,1)})=\infty$, and so $\varrho$ does not have property (P4).

\emph{Property \emph{(P5)}.} Assume that $L<\infty$. Note that \eqref{prop:norminducedbyH:assum} together with $v(L^-)>0$ implies that $\|u\|_{X(0,L)}<\infty$. Let $f\in\Mpl(0,L)$. Since $f^*$ is nonincreasing, we have that $\int_0^Lf^*(s)\d{s}\leq2\int_0^\frac{L}{2}f^*(s)\d{s}$. Since the function $(0,L)\ni t\mapsto f^*(L-t)$ is equimeasurable with $f$, we have that 
\begin{align*}
\varrho(f)&\geq\Big\|u(t)\int_{\nu(t)}^L f^*(L-s)v(s)\d{s}\Big\|_{X(0,L)}\\
&\geq v(L^-)\Big\|u(t)\chi_{(0,\nu^{-1}(\frac{L}{2}))}(t)\int_{\nu(t)}^L f^*(L-s)\d{s}\Big\|_{X(0,L)}\\
&=v(L^-)\Big\|u(t)\chi_{(0,\nu^{-1}(\frac{L}{2}))}(t)\int_0^{L-\nu(t)} f^*(s)\d{s}\Big\|_{X(0,L)}\\
&\geq v(L^-)\|u\chi_{(0,\nu^{-1}(\frac{L}{2}))}\|_{X(0,L)}\int_0^\frac{L}{2} f^*(s)\d{s}\\
&\geq\frac{v(L^-)}{2}\|u\chi_{(0,\nu^{-1}(\frac{L}{2}))}\|_{X(0,L)}\int_0^L f^*(s)\d{s}\\
&\geq\frac{v(L^-)}{2}\|u\chi_{(0,\nu^{-1}(\frac{L}{2}))}\|_{X(0,L)}\int_0^L f(s)\d{s},
\end{align*}
where we used \eqref{ch1:ri:HLg=chiE} in the last inequality. Since $\frac{v(L^-)}{2}\|u\chi_{(0,\nu^{-1}(\frac{L}{2}))}\|_{X(0,L)}\in(0,\infty)$ does not depend on $f$, property (P5) follows. Assume now that $L=\infty$. Recall that \eqref{prop:norminducedbyH:eq11} is satisfied provided that \eqref{prop:norminducedbyH:infinitecase:assump1} is satisfied. Let $f\in\Mpl(0,\infty)$ and $E\subseteq(0,\infty)$ be of finite measure. The function $(0,\infty)\ni t\mapsto f^*(t-|E|)\chi_{(|E|,\infty)}(t)$ is equimeasurable with $f$. By arguing similarly to the case $L<\infty$, we obtain that
\begin{equation*}
\varrho(f)\geq v(2|E|)\|u\chi_{(0,\nu^{-1}(|E|))}\|_{X(0,\infty)}\int_Ef(s)\d{s},
\end{equation*}
whence property (P5) follows.

\emph{Property \emph{(P6)}.} Since the relation $\sim$ is transitive, it plainly follows that $\varrho$ has property (P6).

Note that the necessity of \eqref{prop:norminducedbyH:assum} for $\varrho$ to be a rearrangement\hyp{}invariant function norm was already proved in the paragraph devoted to property (P4).

Assume now that \eqref{prop:norminducedbyH:assum} is satisfied and denote the rearrangement\hyp{}invariant function space induced by $\varrho$ by $Y(0,L)$. We plainly have that
\begin{equation*}
\|H_{u,v,\nu}f\|_{X(0,L)}\leq\|f\|_{Y(0,L)}\quad\text{for every $f\in\Mpl(0,L)$},
\end{equation*}
and so $H_{u,v,\nu}\colon Y(0,L)\to X(0,L)$ is bounded. Next, let $Z(0,L)$ be a rearrangement\hyp{}invariant function space such that $H_{u,v,\nu}\colon Z(0,L)\to X(0,L)$ is bounded. For every $f\in\Mpl(0,\infty)$ and each $h\in\Mpl(0,L)$ equimeasurable with $f$, we have that
\begin{equation*}
\|H_{u,v,\nu}h\|_{X(0,L)}\lesssim\|h\|_{Z(0,L)}=\|f\|_{Z(0,L)}.
\end{equation*}
Therefore,
\begin{equation*}
\|f\|_{Y(0,L)}\lesssim\|f\|_{Z(0,L)}\quad\text{for every $f\in\Mpl(0,L)$}.
\end{equation*}
Hence $Z(0,L)\hookrightarrow Y(0,L)$. Finally, we claim that \eqref{prop:norminducedbyH:assum} needs to be satisfied if there is any rearrangement\hyp{}invariant function space $Z(0,L)$ such that $H_{u,v,\nu}\colon Z(0,L)\to X(0,L)$ is bounded. If $L<\infty$, we plainly have that
\begin{equation*}
\Big\|u(t)\int_{\nu(t)}^L v(s)\d{s}\Big\|_{X(0,L)}=\|H_{u,v,\nu}\chi_{(0,L)}\|_{X(0,L)}\lesssim\|\chi_{(0,L)}\|_{Z(0,L)}<\infty.
\end{equation*}
If $L=\infty$, we can argue as in the paragraph devoted to property (P4) to show that, if \eqref{prop:norminducedbyH:assum} is not satisfied, then
\begin{equation*}
\sup_{h\sim\chi_{(0,1)}}\|H_{u,v,\nu}h\|_{X(0,\infty)}=\infty,
\end{equation*}
whence, thanks to the boundedness of $H_{u,v,\nu}\colon Z(0,\infty)\to X(0,\infty)$,
\begin{equation*}
\infty=\sup_{h\sim\chi_{(0,1)}}\|H_{u,v,\nu}h\|_{X(0,\infty)}\lesssim\|\chi_{(0,1)}\|_{Z(0,\infty)}<\infty,
\end{equation*}
which would be a contradiction.
\end{proof}

\begin{remarks}\hphantom{}
\begin{itemize}
\item The assumption $\nu\in\Dinf{\infty}$ is not overly restrictive. For example, it is satisfied whenever $\nu$ is equivalent to $t\mapsto t^\alpha b(t)$ near $\infty$ for some $\alpha>0$ and a slowly-varying function $b$ (cf.~\cite[Proposition~2.2]{GOT:05}). On the other hand, $\nu(t)=\log^\alpha(t)$ near $\infty$, where $\alpha>0$, is a typical example of a function not satisfying the assumption. The same remark (with the obvious modifications) is true for the assumption $\nu\in\Dinf{0}$, which will appear in \myref{Proposition}{prop:Hiteration}.

\item When $u\equiv1$, \eqref{prop:norminducedbyH:infinitecase:assump2} is equivalent to
\begin{equation*}
\limsup_{t\to\infty}v(\nu(t))\varphi_{X(0,\infty)}(t)<\infty.
\end{equation*}

\item The functional \eqref{prop:norminducedbyH:normdef} is quite complicated; however, we shall see in \myref{Section}{sec:simplification} that it can often be significantly simplified.
\end{itemize}
\end{remarks}

Since, owing to \eqref{prop:RXtoYbddiffHY'toX':RXtoZiffHZ'toX'} and \eqref{ch1:ri:X''=X}, $Y(0,L)$ is the optimal domain space for the operator $H_{u,v,\nu}$ and $X(0,L)$ if and only if $Y'(0,L)$ is the optimal target space for the operator $R_{u,v,\nu^{-1}}$ and $X'(0,L)$, \myref{Proposition}{prop:norminducedbyH} actually also characterizes optimal target spaces for the operator $R_{u,v,\nu}$.
\begin{proposition}\label{prop:Roptimalrange}
Let $\|\cdot\|_{X(0,L)}$ be a rearrangement\hyp{}invariant function norm. Let $\nu\colon(0,L)\to(0,L)$ be an increasing bijection. If $L=\infty$, assume that $\nu^{-1}\in\Dinf{\infty}$. Let $u,v\colon(0,L)\to(0,\infty)$ be nonincreasing. If $L<\infty$, assume that $v(L^-)>0$. Assume that
\begin{equation}\label{prop:Roptimalrange:assum}
\begin{cases}
u(t)\int_{\nu^{-1}(t)}^L v(s)\d{s}\in X'(0,L)\quad&\text{if $L<\infty$},\\
\begin{gathered}[b]u(t)\chi_{(0,\nu(1))}(t)\int_{\nu^{-1}(t)}^1v(s)\d{s}\in X'(0,\infty)\ \text{and}\\ \limsup_{\tau\to\infty}v(\tau)\|u\chi_{(0,\nu(\tau))}\|_{X'(0,\infty)}<\infty\end{gathered}\quad&\text{if $L=\infty$}.
\end{cases}
\end{equation}
Let $\|\cdot\|_{Y(0,L)}$ be the rearrangement\hyp{}invariant function norm whose associate function norm $\|\cdot\|_{Y'(0,L)}$ is defined as 
\begin{equation}\label{prop:Roptimalrange:normdef}
\|f\|_{Y'(0,L)}=\sup_{h\sim f}\Big\|u(t)\int_{\nu^{-1}(t)}^L h(s)v(s)\d{s}\Big\|_{X'(0,L)},\ f\in\Mpl(0,L),
\end{equation}
where the supremum extends over all $h\in\Mpl(0,L)$ equimeasurable with $f$. The rearrangement\hyp{}invariant function space $Y(0,L)$ is the optimal target space for the operator $R_{u,v,\nu}$ and $X(0,L)$. Moreover, if \eqref{prop:Roptimalrange:assum} is not satisfied, then there is no rearrangement\hyp{}invariant function space $Z(0,L)$ such that $R_{u,v,\nu}\colon X(0,L)\to Z(0,L)$ is bounded.
\end{proposition}

Although the functional $\Mpl(0,L)\ni f\mapsto \|H_{u,v,\nu}(f^*)\|_{X(0,L)}$ is usually not a rearrangement\hyp{}invariant function norm, a fact that, in turn, complicates the description of the optimal domain and optimal target spaces for the operators $H_{u,v,\nu}$ and $R_{u,v,\nu^{-1}}$, respectively, it is a rearrangement\hyp{}invariant function norm when $u$, $v$ and $\nu$ are related to each other in such a way that the function $R_{u,v,\nu^{-1}}(g^*)$ is nonincreasing for every $g\in\Mpl(0,L)$. This fact is the content of the following proposition, in which we omit its obvious consequence for optimal spaces.
\begin{proposition}\label{prop:norminducedbyHwhenRnonincreasing}
Let $\|\cdot\|_{X(0,L)}$ be a rearrangement\hyp{}invariant function norm. Let $u\colon(0,L)\to(0,\infty)$ be a nondegenerate nonincreasing function. Let $\nu\colon(0,L)\to(0,L)$ be an increasing bijection. Let $v\colon(0,L)\to(0,\infty)$ be defined by
\begin{equation*}
\frac1{v(t)}=\int_0^{\nu^{-1}(t)} u(s)\d{s} \quad \text{for every $t\in(0,L)$}.
\end{equation*}
Set
\begin{equation*}
\varrho(f)=\Big\|u(t)\int_{\nu(t)}^L f^*(s)v(s)\d{s}\Big\|_{X(0,L)},\ f\in\Mpl(0,L).
\end{equation*}
The functional $\varrho$ is a rearrangement\hyp{}invariant function norm if and only if
\begin{equation}\label{prop:norminducedbyHwhenRnonincreasing:assum}
\Big\|u(t)\chi_{(0,\nu^{-1}(a))}(t)\int_{\nu(t)}^a\frac1{U(\nu^{-1}(s))}\d{s}\Big\|_{X(0,L)} < \infty,
\end{equation}
where $a$ is defined by \eqref{prop:Roptimaldomain:L_or_1}.
\end{proposition}
\begin{proof}
We only sketch the proof (see also \myref{Theorem}{prop:norminducedbyT} and \myref{Remark}{rem:norminducedbyT:varphi_nonincreasing}(i)), which is significantly easier than that of \myref{Proposition}{prop:norminducedbyH}. The functional $\varrho$ plainly possesses properties (P2), (P3) and (P6). It is easy to see that $\varrho$ has property (P4) if and only if \eqref{prop:norminducedbyHwhenRnonincreasing:assum} is satisfied (to this end, note that \eqref{prop:norminducedbyHwhenRnonincreasing:assum} implies \eqref{prop:norminducedbyH:eq11}). As for property (P1), only the subadditivity needs a comment. The key observation is that $(0,L)\ni t\mapsto R_{u,v,\nu^{-1}}(h^*)(t)$ is nonincreasing for every $h\in\Mpl(0,L)$ inasmuch as it is the integral mean of a nonnegative nonincreasing function over the interval $(0,\nu^{-1}(t))$ with respect to the measure $u(s)\d{s}$. Hence, thanks to \eqref{ch1:ri:normX''down}, \eqref{RaHdual} and \eqref{ch1:ri:subadditivityofdoublestar} combined with the Hardy lemma \eqref{ch1:ri:hardy-lemma}, we have that
\begin{align*}
\varrho(f+g)&=\|H_{u,v,\nu}((f+g)^*)\|_{X(0,L)}=\sup_{\substack{h\in{\Mpl(0,L)}\\\|h\|_{X'(0,L)}\leq1}}\int_0^L H_{u,v,\nu}((f+g)^*)(t)h^*(t)\d{t} \\
&= \sup_{\substack{h\in{\Mpl(0,L)}\\\|h\|_{X'(0,L)}\leq1}}\int_0^L (f+g)^*(t)R_{u,v,\nu^{-1}}(h^*)(t)\d{t} \\
&\leq \sup_{\substack{h\in{\Mpl(0,L)}\\\|h\|_{X'(0,L)}\leq1}}\int_0^L f^*(t)R_{u,v,\nu^{-1}}(h^*)(t)\d{t} \\
& \quad + \sup_{\substack{h\in{\Mpl(0,L)}\\\|h\|_{X'(0,L)}\leq1}}\int_0^L g^*(t)R_{u,v,\nu^{-1}}(h^*)(t)\d{t} \\
&= \sup_{\substack{h\in{\Mpl(0,L)}\\\|h\|_{X'(0,L)}\leq1}}\int_0^L H_{u,v,\nu}(f^*)(t)h^*(t)\d{t} \\
& \quad + \sup_{\substack{h\in{\Mpl(0,L)}\\\|h\|_{X'(0,L)}\leq1}}\int_0^L H_{u,v,\nu}(g^*)(t)h^*(t)\d{t} \\
&= \varrho(f) + \varrho(g)
\end{align*}
for every $f,g\in\Mpl(0,L)$. Finally, as for the validity of property (P5), owing to \eqref{ch1:ri:normX''down}, \eqref{RaHdual}, the monotonicity of $R_{u,v,\nu^{-1}}(g^*)$ and the Hardy--Littlewood inequality \eqref{ch1:ri:HL}, we have that
\begin{align*}
\varrho(f)&\geq \varrho(f\chi_E)=\sup_{\substack{g\in{\Mpl(0,L)}\\\|g\|_{X'(0,L)}\leq1}}\int_0^L (f\chi_E)^*(t)R_{u,v,\nu^{-1}}(g^*)(t)\d{t} \\
&= \sup_{\substack{g\in{\Mpl(0,L)}\\\|g\|_{X'(0,L)}\leq1}}\int_0^{|E|} (f\chi_E)^*(t)R_{u,v,\nu^{-1}}(g^*)(t)\d{t} \\
&\geq \int_0^{|E|} (f\chi_E)^*(t)\d{t} \sup_{\substack{g\in{\Mpl(0,L)}\\\|g\|_{X'(0,L)}\leq1}}R_{u,v,\nu^{-1}}(g^*)(|E|) \\
&\geq R_{u,v,\nu^{-1}}\Big(\frac{\chi_{(0,|E|)}}{\|\chi_{(0,|E|)}\|_{X'(0,L)}}\Big)(|E|) \int_0^{|E|} (f\chi_E)^*(t)\d{t}\\
&\geq\frac{v(|E|)}{\varphi_{X'(0,L)}(|E|)}U(\nu^{-1}(|E|)) \int_E f(t)\d{t}
\end{align*}
for every $f\in\Mpl(0,L)$ and $E\subseteq(0,L)$ having finite measure.
\end{proof}

\begin{remark}
Note that $(R_{u_1, v_1, \nu_1} + H_{u_2, v_2, \nu_2})\colon X(0,L) \to Y(0,L)$ is bounded if and only if both $R_{u_1, v_1, \nu_1}$ and $H_{u_2, v_2, \nu_2}$ are bounded from $X(0,L)$ to $Y(0,L)$. Furthermore, it is easy to see that $(Y_1 + Y_2)(0,L)$ is the optimal target space for $R_{u_1, v_1, \nu_1} + H_{u_2, v_2, \nu_2}$ and $X(0,L)$, where $Y_1(0,L)$ and $Y_2(0,L)$ are the optimal target spaces for $R_{u_1, v_1, \nu_1}$ and $H_{u_2, v_2, \nu_2}$, respectively, and $X_1(0,L) \cap X_2(0,L)$ is the optimal domain space for $R_{u_1, v_1, \nu_1} + H_{u_2, v_2, \nu_2}$ and $Y(0,L)$, where $X_1(0,L)$ and $X_2(0,L)$ are the optimal domain spaces for $R_{u_1, v_1, \nu_1}$ and $H_{u_2, v_2, \nu_2}$.
\end{remark}

We conclude this subsection with an important result ensuring that, to verify the boundedness of $H_{u, v, \nu}$ between a pair of rearrangement\hyp{}invariant function spaces, it is sufficient to verify it on nonincreasing functions. While it is an easy consequence of Hardy--Littlewood inequality \eqref{ch1:ri:HL} that this is the case for the operator $R_{u, v, \nu}$, provided that $u$ is nonincreasing, the validity of such a result for the operator $H_{u, v, \nu}$ is far from obvious because this time the integration is not carried out over a right-neighborhood of $0$. Such a result was first obtained by Cianchi, Pick and Slav\'\i kov\'a in \cite[Theorem~9.5]{CPS:15} for $u\equiv1$, $\nu=\id$ and $L<\infty$. Later, Pe\v sa generalized their result to cover also the case $L=\infty$ in \cite[Theorem~3.10]{Pe:20}. In \cite{CM:20}, we needed such a result for $\nu(t)=t^\alpha$, $\alpha>0$, and $u\not\equiv1$, and, while we felt certain that their proofs would carry over to the needed setting, we still had to carefully check them because there is plenty of fine analysis involved. The following proposition extends the result to the generality considered in this paper. It turns out that their proofs can easily be adapted for our setting. Our proof is actually in a way simpler because they considered operators with kernels.
\begin{proposition}\label{prop:copsonrestrvsunrestr}
Let $\|\cdot\|_{X(0,L)}$ and $\|\cdot\|_{Y(0,L)}$ be rearrangement\hyp{}invariant function norms. Let $u,v\colon(0,L)\to(0, \infty)$ be nonincreasing. Let $\nu\colon(0,L)\to(0,L)$ be an increasing bijection. Assume that $\nu^{-1}\in\Dsup[\theta]{0}$ for some $\theta>1$. If $L=\infty$, assume that $\nu^{-1}\in\Dsup[\theta]{\infty}$. The following two statements are equivalent.
\begin{enumerate}[(i)]
\item There is a positive constant $C$ such that
\begin{equation}\label{prop:copsonrestrvsunrestr:unrestr}
\Big\|u(t)\int_{\nu(t)}^L|f(s)|v(s)\d{s}\Big\|_{Y(0,L)}\leq C\|f\|_{X(0,L)}
\end{equation}
for every $f\in\M(0,L)$.
\item There is a positive constant $C$ such that
\begin{equation}\label{prop:copsonrestrvsunrestr:restr}
\Big\|u(t)\int_{\nu(t)}^Lf^*(s)v(s)\d{s}\Big\|_{Y(0,L)}\leq C\|f\|_{X(0,L)}
\end{equation}
for every $f\in\M(0,L)$.
\end{enumerate}
Moreover, if \eqref{prop:copsonrestrvsunrestr:restr} holds with a constant $C$, then \eqref{prop:copsonrestrvsunrestr:unrestr} holds with the constant $C\frac{\theta}{\theta-1}\sup_{t\in(0,L)}\frac{\nu^{-1}(t)}{\nu^{-1}(\frac{t}{\theta})}$.
\end{proposition}
\begin{proof}
Since (i) plainly implies (ii), we only need to prove that (ii) implies (i). Since the quantities in \eqref{prop:copsonrestrvsunrestr:unrestr} and \eqref{prop:copsonrestrvsunrestr:restr} do not change when the function $v$ is redefined on a countable set, we may assume that $v$ is left continuous. Note that $H_{u,v,\nu}f$ is nonincreasing for every $f\in\M(0,L)$. Hence, thanks to \eqref{ch1:ri:normX''down} and \eqref{RaHdual}, in order to prove that (ii) implies (i), we need to show that
\begin{equation}\label{prop:copsonrestrvsunrestr:eq8}
\sup_{\substack{f\in{\Mpl(0,L)}\\\|f\|_{X(0,L)}\leq1}}\sup_{\substack{g\in{\Mpl(0,L)}\\\|g\|_{Y'(0,L)}\leq1}}\int_0^Lf(s)R_{u,v,\nu^{-1}}(g^*)(s)\d{s}\lesssim\sup_{\substack{f\in{\Mpl(0,L)}\\\|f\|_{X(0,L)}\leq1}}\sup_{\substack{g\in{\Mpl(0,L)}\\\|g\|_{Y'(0,L)}\leq1}}\int_0^Lf^*(s)R_{u,v,\nu^{-1}}(g^*)(s)\d{s}.
\end{equation}

We define the operator $G$ as 
\begin{equation*}
Gg(t)=\sup_{\tau\in[t,L)}R_{u,v,\nu^{-1}}(g^*)(\tau),\ t\in(0,L),
\end{equation*}
for every $g\in\Mpl(0,L)$. Note that $Gg$ is nonincreasing for every $g\in\Mpl(0,L)$. Fix $g\in\Mpl(0,L)$ such that $|\{t\in(0,L)\colon g(t)>0\}|<\infty$, and set
\begin{equation*}
E=\big\{t\in(0,L)\colon R_{u,v,\nu^{-1}}(g^*)(t)<Gg(t)\big\}.
\end{equation*}
It can be shown that there is a countable system $\{(a_k,b_k)\}_{k\in\mathcal{I}}$ of mutually disjoint, bounded intervals in $(0,L)$ such that
\begin{align}
E&=\bigcup_{k\in\mathcal{I}}(a_k,b_k);\label{prop:copsonrestrvsunrestr:eq4}\\
Gg(t)&=R_{u,v,\nu^{-1}}(g^*)(t)\quad\text{if $t\in(0,L)\setminus E$;}\label{prop:copsonrestrvsunrestr:eq5}\\
Gg(t)&=R_{u,v,\nu^{-1}}(g^*)(b_k)\quad\text{if $t\in(a_k,b_k)$ for  $k\in\mathcal{I}$}.\label{prop:copsonrestrvsunrestr:eq1}
\end{align}
This was proved in \cite[Proposition~9.3]{CPS:15} for $L<\infty$ and in \cite[Lemma~3.9]{Pe:20} for $L=\infty$. Their proofs are for $u\equiv1$ and $\nu=\id$, but the fact that $g^*u$ is nonincreasing and $R_{u,v,\nu^{-1}}(g^*)$ is upper semicontinuous remains valid in our situation, and so it can be readily seen that their proofs carry over verbatim to our setting.

Note that $M=\sup_{t\in(0,L)}\frac{\nu^{-1}(t)}{\nu^{-1}(\frac{t}{\theta})}<\infty$. Set $\sigma=\frac{\theta}{\theta-1}\in(1,\infty)$. Since $v$ and $g^*u$ are nonincreasing, we have that, for every $k\in\mathcal{I}$,
\begin{align}
(b_k - a_k)R_{u,v,\nu^{-1}}(g^*)(b_k)&=\sigma\int_{\frac{a_k+(\sigma-1)b_k}{\sigma}}^{b_k}R_{u,v,\nu^{-1}}(g^*)(b_k)\d{t}\notag\\
\begin{split}\label{prop:copsonrestrvsunrestr:eq6}
&=\sigma\int_{\frac{a_k+(\sigma-1)b_k}{\sigma}}^{b_k}\frac{v(b_k)}{\nu^{-1}(b_k)}\nu^{-1}(b_k)\int_0^{\nu^{-1}(b_k)}g^*(s)u(s)\d{s}\,\d{t}\\
&\leq\sigma\int_{\frac{a_k+(\sigma-1)b_k}{\sigma}}^{b_k}\frac{v(t)}{\nu^{-1}(t)}\nu^{-1}(b_k)\int_0^{\nu^{-1}(t)}g^*(s)u(s)\d{s}\,\d{t}\\
&\leq\sigma\frac{\nu^{-1}(b_k)}{\nu^{-1}(\frac{a_k+(\sigma-1)b_k}{\sigma})}\int_{\frac{a_k+(\sigma-1)b_k}{\sigma}}^{b_k}R_{u,v,\nu^{-1}}(g^*)(t)\d{t}\\
&\leq\sigma\frac{\nu^{-1}(b_k)}{\nu^{-1}(\frac{b_k}{\theta})}\int_{\frac{a_k+(\sigma-1)b_k}{\sigma}}^{b_k}R_{u,v,\nu^{-1}}(g^*)(t)\d{t}\\
&\leq\sigma M\int_{\frac{a_k+(\sigma-1)b_k}{\sigma}}^{b_k}R_{u,v,\nu^{-1}}(g^*)(t)\d{t}\\
&\leq\sigma M\int_{a_k}^{b_k}R_{u,v,\nu^{-1}}(g^*)(t)\d{t},
\end{split}
\end{align}
where we used the fact that $v$ and $(g^*u)^{**}$ are nonincreasing in the first inequality.

Consider the averaging operator $A$ defined as
\begin{equation*}
Af=f^*\chi_{(0,L)\setminus E}+\sum_{k\in\mathcal I}\Big(\frac1{b_k-a_k}\int_{a_k}^{b_k}f^*(s)\d{s}\Big)\chi_{(a_k, b_k)},\ f\in\Mpl(0,L).
\end{equation*}
Note that $Af$ is a nonincreasing function for every $f\in\Mpl(0,L)$. Furthermore, it is known (\cite[Chapter~2, Theorem~4.8]{BS}) that
\begin{equation}\label{prop:copsonrestrvsunrestr:eq2}
\|Af\|_{X(0,L)}\leq\|f\|_{X(0,L)}\quad\text{for every $f\in\Mpl(0,L)$}.
\end{equation}
We have, for every $f\in\Mpl(0,L)$, that
\begin{align}\label{prop:copsonrestrvsunrestr:eq3}
\int_0^Lf(t)R_{u,v,\nu^{-1}}(g^*)(t)\d{t}&\leq\int_0^Lf(t)Gg(t)\d{t}\leq\int_0^Lf^*(t)Gg(t)\d{t} \notag\\
\begin{split}
&=\int_{(0,L)\setminus E}f^*(t)R_{u,v,\nu^{-1}}(g^*)(t)\d{t}+\sum_{k\in\mathcal{I}}\int_{a_k}^{b_k}f^*(t)R_{u,v,\nu^{-1}}(g^*)(b_k)\d{t}\\
&\leq\int_0^Lf^*(t)\chi_{(0,L)\setminus E}(t)R_{u,v,\nu^{-1}}(g^*)(t)\d{t}\\
&\quad+ \sigma M\sum_{k\in\mathcal{I}}\Big(\frac1{b_k-a_k}\int_{a_k}^{b_k}f^*(t)\d{t}\Big)\Big(\int_{a_k}^{b_k}R_{u,v,\nu^{-1}}(g^*)(t)\d{t}\Big)\\
&\leq \sigma M\int_0^LAf(t)R_{u,v,\nu^{-1}}(g^*)(t)\d{t},
\end{split}
\end{align}
owing to the Hardy--Littlewood inequality \eqref{ch1:ri:HL}, \eqref{prop:copsonrestrvsunrestr:eq4}, \eqref{prop:copsonrestrvsunrestr:eq5}, \eqref{prop:copsonrestrvsunrestr:eq1}, and \eqref{prop:copsonrestrvsunrestr:eq6}.
If $L=\infty$ and $g\in\Mpl(0,\infty)$ is positive on a set of infinite measure, we consider $g\chi_{(0,n)}\nearrow g$, $n\to\infty$, and obtain \eqref{prop:copsonrestrvsunrestr:eq3} even for such functions $g$ thanks to the monotone convergence theorem; hence we have proved that
\begin{equation}\label{prop:copsonrestrvsunrestr:eq7}
\int_0^Lf(t)R_{u,v,\nu^{-1}}(g^*)(t)\d{t}\leq \sigma M\int_0^LAf(t)R_{u,v,\nu^{-1}}(g^*)(t)\d{t}\quad\text{for every $f,g\in\Mpl(0,L)$}.
\end{equation}
By combining \eqref{prop:copsonrestrvsunrestr:eq2} and \eqref{prop:copsonrestrvsunrestr:eq7}, we obtain that
\begin{equation*}
\int_0^Lf(t)R_{u,v,\nu^{-1}}(g^*)(t)\d{t}\leq \sigma M\sup_{\substack{h\in{\Mpl(0,L)}\\\|h\|_{X(0,L)}\leq1}}\int_0^Lh^*(t)R_{u,v,\nu^{-1}}(g^*)(t)\d{t}
\end{equation*}
for every $f\in\Mpl(0,L)$, $\|f\|_{X(0,L)}\leq1$, and $g\in\Mpl(0,L)$. Note that here we used the fact that $Af$ is nonincreasing for every $f\in\Mpl(0,L)$. By taking the supremum over all $f,g\in\Mpl(0,L)$ from the closed unit balls of $X(0,L)$ and $Y'(0,L)$, respectively, we obtain \eqref{prop:copsonrestrvsunrestr:eq8} with the multiplicative constant equal to $\sigma M$.
\end{proof}

\myref{Proposition}{prop:copsonrestrvsunrestr} together with \myref{Proposition}{prop:RXtoYbddiffHY'toX'} has the following important corollary, in which the first equality is just a consequence of the Hardy--Littlewood inequality \eqref{ch1:ri:HL} combined with the obvious inequality
\begin{equation*}
\sup_{\|f\|_{X(0,L)}\leq1}\|R_{u,v,\nu}(f^*)\|_{Y(0,L)}\leq\sup_{\|f\|_{X(0,L)}\leq1}\|R_{u,v,\nu}f\|_{Y(0,L)}.
\end{equation*}
\begin{corollary}\label{cor:copsonrestrvsunrestr}
Let $\|\cdot\|_{X(0,L)}$ and $\|\cdot\|_{Y(0,L)}$ be rearrangement\hyp{}invariant function norms. Let $u,v\colon(0,L)\to(0, \infty)$ be nonincreasing. Let $\nu\colon(0,L)\to(0,L)$ be an increasing bijection. Assume that $\nu\in\Dsup{0}$. If $L=\infty$, assume that $\nu\in\Dsup{\infty}$. We have that
\begin{equation*}
\begin{aligned}
\sup_{\|f\|_{X(0,L)}\leq1}\|R_{u,v,\nu}(f^*)\|_{Y(0,L)}&=\sup_{\|f\|_{X(0,L)}\leq1}\|R_{u,v,\nu}f\|_{Y(0,L)}\\
&=\sup_{\|g\|_{Y'(0,L)}\leq1}\|H_{u,v,\nu^{-1}}g\|_{X'(0,L)}\\
&\approx\sup_{\|g\|_{Y'(0,L)}\leq1}\|H_{u,v,\nu^{-1}}(g^*)\|_{X'(0,L)}.
\end{aligned}
\end{equation*}
\end{corollary}

\subsection{Simplification of optimal r.i.~norms}\label{sec:simplification}
In this subsection, we will study deeper properties of the important but complicated functional \eqref{prop:norminducedbyH:normdef}. In particular, we shall see that the functional is actually often equivalent to a significantly more manageable functional (cf.~\cite[Theorem~4.2]{EMMP:20}). To this end, we need to introduce a supremum operator. For a fixed function $\varphi\colon(0,L)\to(0,\infty)$, we define the operator $T_\varphi$ as
\begin{equation}\label{opTdef}
T_\varphi f(t)=\frac1{\varphi(t)}\esssup_{s\in[t,L)}\varphi(s)f^*(s),\ t\in(0,L),\ f\in\M(0,L).
\end{equation}
Note that $T_\varphi f(t)=\frac1{\varphi(t)}\sup_{s\in[t,L)}\varphi(s)f^*(s)$ for every $t\in(0, L)$ provided that $\varphi$ is nondecreasing and/or right-continuous. If $\varphi$ is nonincreasing, we have that $T_\varphi f(t)= \frac{f^*(t)}{\varphi(t)}\varphi(t^+)$ for every $t\in(0,L)$, and so $T_\varphi f=f^*$ possibly up to a countably many points.
\begin{proposition}\label{prop:norminducedbyT}
Let $\|\cdot\|_{X(0,L)}$ be a rearrangement\hyp{}invariant function norm. Let $\nu\colon(0,L)\to(0,L)$ be an increasing bijection. If $L=\infty$, assume that $\nu\in\Dinf{\infty}$. Let $u\colon(0,L)\to(0,\infty)$ be nonincreasing. Let $v\colon(0,L)\to(0,\infty)$ be defined by
\begin{equation*}
\frac1{v(t)}=\int_0^{\nu^{-1}(t)}\xi(s)\d{s} \quad \text{for every $t\in(0,L)$},
\end{equation*}
where $\xi\colon(0,L)\to(0,\infty)$ is a measurable function. If $L<\infty$, assume that $v(L^-)>0$. Assume that
\begin{equation*}
\Big\|u(t)\chi_{(0,\nu^{-1}(a))}(t)\int_{\nu(t)}^a v(s)\d{s}\Big\|_{X(0,L)}<\infty,
\end{equation*}
where $a$ is defined by \eqref{prop:Roptimaldomain:L_or_1}, and that the operator $T_\varphi$ defined by \eqref{opTdef} with $\varphi=\frac{u}{\xi}$ is bounded on $X'(0,L)$. Let $\varrho$ be the functional defined by \eqref{prop:norminducedbyH:normdef} and set
\begin{equation*}
\widetilde{\varrho}(f)=\sup_{\substack{g\in{\Mpl(0,L)}\\\|g\|_{X'(0,L)}\leq1}}\int_0^Lf^*(s)v(t)\int_0^{\nu^{-1}(t)}T_\varphi g(s)u(s)\d{s}\d{t},\ f\in\Mpl(0,L).
\end{equation*}
The functionals $\varrho$ and $\widetilde{\varrho}$ are rearrangement\hyp{}invariant function norms. Furthermore, we have that
\begin{align}\label{prop:norminducedbyT:equivalencewithHf*}
\begin{split}
\|H_{u,v,\nu}(f^*)\|_{X(0,L)}&\leq\sup_{h\sim f}\Big\|u(t)\int_{\nu(t)}^L h(s)v(s)\d{s}\Big\|_{X(0,L)}\leq\widetilde{\varrho}(f)\\
&\leq\|T_\varphi\|_{X'(0,L)}\|H_{u,v,\nu}(f^*)\|_{X(0,L)}
\end{split}
\end{align}
for every $f\in\Mpl(0,L)$, where $\|T_\varphi\|_{X'(0,L)}$ stands for the operator norm of $T_\varphi$ on $X'(0,L)$. In particular, the rearrangement\hyp{}invariant function norms $\varrho$ and $\widetilde{\varrho}$ are equivalent.
\end{proposition}
\begin{proof}
Since $f\sim f^*$ for every $f\in\Mpl(0,L)$, the first inequality in \eqref{prop:norminducedbyT:equivalencewithHf*} plainly holds. As for the second inequality, note that the function $(0,L)\ni t\mapsto R_{u,v,\nu^{-1}}(T_\varphi g)(t)$ is nonincreasing for every $g\in\Mpl(0,L)$ because it is the integral mean of the nonincreasing function $(0,L)\ni s\mapsto \esssup_{\tau\in[s,L)}\varphi(\tau)g^*(\tau)$ over the interval $(0,\nu^{-1}(t))$ with respect to the measure $\xi(s)\d{s}$. Consequently, for every $f\in\Mpl(0,L)$ and every $h\in\Mpl(0,L)$ equimeasurable with $f$, we have that
\begin{align}
\|H_{u,v,\nu}h\|_{X(0,L)}&=\sup_{\substack{g\in{\Mpl(0,L)}\\\|g\|_{X'(0,L)}\leq1}}\int_0^Lh(t)R_{u,v,\nu^{-1}}(g^*)(t)\d{t} \notag\\
\begin{split}\label{prop:norminducedbyT:eq1}
&\leq\sup_{\substack{g\in{\Mpl(0,L)}\\\|g\|_{X'(0,L)}\leq1}}\int_0^Lh(t)R_{u,v,\nu^{-1}}(T_\varphi g)(t)\d{t}\\
&\leq\sup_{\substack{g\in{\Mpl(0,L)}\\\|g\|_{X'(0,L)}\leq1}}\int_0^Lh^*(t)R_{u,v,\nu^{-1}}(T_\varphi g)(t)\d{t}\\
&=\sup_{\substack{g\in{\Mpl(0,L)}\\\|g\|_{X'(0,L)}\leq1}}\int_0^Lf^*(t)R_{u,v,\nu^{-1}}(T_\varphi g)(t)\d{t}\\
&=\widetilde{\varrho}(f),
\end{split}
\end{align}
where we used \eqref{ch1:ri:normX''down} (note that the function $H_{u,v,\nu}h$ is nonincreasing for every $h\in\Mpl(0,L)$) together with \eqref{RaHdual} in the first equality, the pointwise estimate $g^*(t)\leq T_\varphi g(t)$ for a.e.~$t\in(0,L)$ in the first inequality, the Hardy--Littlewood inequality \eqref{ch1:ri:HL} in the second inequality, and the equimeasurability of $f$ and $h$ in the last inequality. Hence the second inequality in \eqref{prop:norminducedbyT:equivalencewithHf*} follows from \eqref{prop:norminducedbyT:eq1}. As for the third inequality in \eqref{prop:norminducedbyT:equivalencewithHf*}, we have that
\begin{align*}
\sup_{\substack{g\in{\Mpl(0,L)}\\\|g\|_{X'(0,L)}\leq1}}\int_0^Lf^*(t)R_{u,v,\nu^{-1}}(T_\varphi g)(t)\d{t}&=\sup_{\substack{g\in{\Mpl(0,L)}\\\|g\|_{X'(0,L)}\leq1}}\int_0^LT_\varphi g(t)H_{u,v,\nu}(f^*)(t)\d{t}\\
&\leq\|H_{u,v,\nu}(f^*)\|_{X(0,L)}\sup_{\substack{g\in{\Mpl(0,L)}\\\|g\|_{X'(0,L)}\leq1}}\|T_\varphi g\|_{X'(0,L)}\\
&=\|T_\varphi\|_{X'(0,L)}\|H_{u,v,\nu}(f^*)\|_{X(0,L)}
\end{align*}
for every $f\in\Mpl(0,L)$ thanks to \eqref{RaHdual} and the H\"older inequality \eqref{ch1:ri:holder}.

Second, we shall prove that the functional $\varrho$, defined by \eqref{prop:norminducedbyH:normdef}, is a rearrangement\hyp{}invariant function norm. If $L<\infty$, this follows immediately from \myref{Proposition}{prop:norminducedbyH}. If $L=\infty$, owing to \myref{Proposition}{prop:norminducedbyH} again, we only need to verify that \eqref{prop:norminducedbyH:infinitecase:assump2} is satisfied. To this end, it follows from the proof of property (P4) of $\varrho$ that, if \eqref{prop:norminducedbyH:infinitecase:assump2} did not hold, then we would have
\begin{equation*}
\sup_{h\sim \chi_{(0,1)}}\Big\|u(t)\int_{\nu(t)}^\infty h(s)v(s)\d{s}\Big\|_{X(0,\infty)}=\infty.
\end{equation*}
However, thanks to \eqref{prop:norminducedbyT:equivalencewithHf*}, we have that
\begin{align*}
\sup_{h\sim \chi_{(0,1)}}\Big\|u(t)\int_{\nu(t)}^L h(s)v(s)\d{s}\Big\|_{X(0,\infty)}&\approx\|H_{u,v,\nu}\chi_{(0,1)}\|_{X(0,L)}\\
&=\Big\|u(t)\chi_{(0,\nu^{-1}(1))}(t)\int_{\nu(t)}^1v(s)\d{s}\Big\|_{X(0,\infty)}<\infty.
\end{align*}
Therefore, \eqref{prop:norminducedbyH:infinitecase:assump2} is satisfied.

Finally, now that we know that the functionals $\varrho$ and $\widetilde{\varrho}$ are equivalent and the former is a rearrangement\hyp{}invariant function norm, it readily follows that $\widetilde{\varrho}$, too, is a rearangement-invariant function norm once we observe that $\widetilde{\varrho}$ is subadditive. The subadditivity follows from
\begin{align*}
\widetilde{\varrho}(f+g)&=\sup_{\substack{h\in{\Mpl(0,L)}\\\|h\|_{X'(0,L)}\leq1}}\int_0^L(f+g)^*(t)R_{u,v,\nu^{-1}}(T_\varphi h)(t)\d{t}\\
&\leq\sup_{\substack{h\in{\Mpl(0,L)}\\\|h\|_{X'(0,L)}\leq1}}\int_0^Lf^*(t)R_{u,v,\nu^{-1}}(T_\varphi h)(t)\d{t}\\
&\quad+\sup_{\substack{h\in{\Mpl(0,L)}\\\|h\|_{X'(0,L)}\leq1}}\int_0^Lg^*(t)R_{u,v,\nu^{-1}}(T_\varphi h)(t)\d{t}\\
&=\widetilde{\varrho}(f)+\widetilde{\varrho}(g)\quad\text{for every $f,g\in\Mpl(0,L)$},
\end{align*}
where we used \eqref{ch1:ri:subadditivityofdoublestar} together with the Hardy lemma \eqref{ch1:ri:hardy-lemma} (recall that the function $R_{u,v,\nu^{-1}}(T_\varphi h)$ is nonincreasing for every $h\in\Mpl(0,L)$).
\end{proof}

\begin{remarks}\label{rem:norminducedbyT:varphi_nonincreasing}\hphantom{}
\begin{enumerate}[(i)]
\item If $\varphi=\frac{u}{\xi}$ is (equivalent to) a nonincreasing function, $T_\varphi f(t)$ is (equivalent to) $f^*(t)$ for a.e.~$t\in(0,L)$; hence $T_\varphi$ is bounded on every rearrangement\hyp{}invariant function space in this case. Furthermore, when $\varphi=\frac{u}{\xi}$ is nonincreasing, the norm of $T_\varphi$ on every rearrangement\hyp{}invariant function space is equal to $1$; therefore all the inequalities in \eqref{prop:norminducedbyT:equivalencewithHf*} are actually equalities (cf.~\myref{Proposition}{prop:norminducedbyHwhenRnonincreasing}) in this case.

\item The boundedness of $T_\varphi$ on a large number of rearrangement\hyp{}invariant function spaces is characterized by \cite[Theorem~3.2]{GOP:06}.
\end{enumerate}
\end{remarks}

By combining \myref{Proposition}{prop:norminducedbyT} and \myref{Proposition}{prop:Roptimalrange}, we obtain the following proposition, which tells us that the optimal target space for the operator $R_{u,v,\nu}$ and a rearrangement\hyp{}invariant function space $X(0,L)$ has a much more manageable description than that given by \myref{Proposition}{prop:Roptimalrange} provided that the supremum operator $T_\varphi$ defined by \eqref{opTdef} with an appropriate function $\varphi$ is bounded on $X(0,L)$.
\begin{proposition}\label{prop:RoptimalrangeTbound}
Let $\|\cdot\|_{X(0,L)}$ be a rearrangement\hyp{}invariant function norm. Let $\nu\colon(0,L)\to(0,L)$ be an increasing bijection. If $L=\infty$, assume that $\nu^{-1}\in\Dinf{\infty}$. Let $u\colon(0,L)\to(0,\infty)$ be nonincreasing. Let $v\colon(0,L)\to(0,\infty)$ be defined by
\begin{equation*}
\frac1{v(t)}=\int_0^{\nu(t)}\xi(s)\d{s} \quad \text{for every $t\in(0,L)$},
\end{equation*}
where $\xi\colon(0,L)\to(0,\infty)$ is a measurable function. If $L<\infty$, assume that $v(L^-)>0$. Furthermore, assume that
\begin{equation*}
\Big\|u(t)\chi_{(0,\nu(a))}(t)\int_{\nu^{-1}(t)}^a v(s)\d{s}\Big\|_{X'(0,L)} < \infty,
\end{equation*}
where $a$ is defined by \eqref{prop:Roptimaldomain:L_or_1}. Finally, assume that the operator $T_\varphi$ is bounded on $X(0,L)$, where $\varphi=\frac{u}{\xi}$. Let $\|\cdot\|_{Y(0,L)}$ be the rearrangement\hyp{}invariant function norm whose associate function norm $\|\cdot\|_{Y'(0,L)}$ is defined as 
\begin{equation*}
\|f\|_{Y'(0,L)}=\sup_{\substack{g\in{\Mpl(0,L)}\\\|g\|_{X(0,L)}\leq1}}\int_0^Lf^*(s)v(t)\int_0^{\nu(t)}T_\varphi g(s)u(s)\d{s}\d{t},\ f\in\Mpl(0,L).
\end{equation*}
The rearrangement\hyp{}invariant function space $Y(0,L)$ is the optimal target space for the operator $R_{u,v,\nu}$ and $X(0,L)$. Moreover,
\begin{equation*}
\|H_{u,v,\nu^{-1}}(f^*)\|_{X'(0,L)}\leq \|f\|_{Y'(0,L)}\leq\|T_\varphi\|_{X(0,L)}\|H_{u,v,\nu^{-1}}(f^*)\|_{X'(0,L)}
\end{equation*}
for every $f\in\Mpl(0,L)$, where $\|T_\varphi\|_{X(0,L)}$ stands for the operator norm of $T_\varphi$ on $X(0,L)$.
\end{proposition}

\begin{remark}
Owing to \myref{Remark}{rem:optimal_for_R_iff_for_H}, \myref{Proposition}{prop:RoptimalrangeTbound} can also be used to get a simpler description of optimal domain spaces for the operator $H_{u,v,\nu}$.
\end{remark}

We already know that a sufficient condition for simplification of the complicated function norm \eqref{prop:Roptimalrange:normdef} is boundedness of a certain supremum operator. We shall soon see that the connection between the supremum operator and the question of whether the supremum in the function norm can be `simplified' is actually much tighter than it may look at first glance. Furthermore, not only is the boundedness of the supremum operator often also necessary for simplifying \eqref{prop:Roptimalrange:normdef}, but it also goes hand in hand with the notion of being an optimal function space and a certain interpolation property of the rearrangement\hyp{}invariant function space on which the supremum operator acts.

As the following theorem shows, there is a connection between a rearrangement\hyp{}invariant function space $X(0,L)$ being an interpolation space with respect to a certain pair of endpoint spaces and the boundedness of $T_\varphi$ on the associate space of $X(0,L)$ (cf.~\citep[Theorem~3.12]{KP:06}). We say that a measurable a.e.~positive function on $(0,L)$ satisfies \emph{the averaging condition \eqref{averaging_condition}} (cf.~\cite[Lemma~2.3]{S:79}) if
\begin{equation}\label{averaging_condition}
\esssup_{t\in(0, L)} \frac1{t w(t)} \int_0^t w(s) \d{s} < \infty,
\end{equation}
in which $w$ temporarily denotes the function in question. The value of the essential supremum will be referred to as \emph{the averaging constant} of the function.
\begin{theorem}\label{thm:Tbounded_iff_Xinterpolation}
Let $\|\cdot\|_{X(0,L)}$ be a rearrangement\hyp{}invariant function norm. Let $\varphi\colon (0,L) \to (0, \infty)$ be a measurable function that is equivalent to a continuous nondecreasing function. Set $\xi = \frac1{\varphi}$. Assume that $\xi$ satisfies the averaging condition~\eqref{averaging_condition}. Set $\psi(t)=\frac{t}{\int_0^t \xi(s) \d{s}}$, $t\in(0,L)$. Consider the following three statements.
\begin{enumerate}[(i)]
\item The operator $T_\varphi$, defined by \eqref{opTdef}, is bounded on $X'(0,L)$.
\item $X(0,L) \in \Int\big(\Lambda^1_\xi(0,L), L^\infty(0,L)\big)$.
\item $X'(0,L) \in \Int\big(L^1(0,L), M_\psi(0,L)\big)$.
\end{enumerate}
If $L<\infty$, then the three statements are equivalent to each other. If $L=\infty$, then (i) implies (ii), and (iii) implies (i).
\end{theorem}
\begin{proof}
We start off by noting that we may without loss of generality assume that $\varphi$ is continuous and nondecreasing. Furthermore, $(\Lambda^1_\xi)'(0,L)=M_\psi(0,L)$ (\cite[Theorem~10.4.1]{PKJF:13}) and
\begin{equation}\label{thm:Tbounded_iff_Xinterpolation:eq6}
\psi\approx\varphi \quad \text {on $(0,L)$}
\end{equation}
thanks to the fact that $\xi$ satisfies the averaging condition~\eqref{averaging_condition} and is (equivalent to) a nonincreasing function. We shall show that (i) implies (ii), whether $L$ is finite or infinite. First, we observe that $X(0,L)$ is an intermediate space between $\Lambda^1_\xi(0,L)$ and $L^\infty(0,L)$. Set $\Xi_L=\int_0^L \xi(s)\d{s}\in(0,\infty]$ (note that $\Xi_L<\infty$ if $L<\infty$). Let $\Xi^{-1}\colon (0, \Xi_L) \to (0,L)$ be the increasing bijection that is inverse to the function $(0,L) \ni t \mapsto \int_0^t \xi(s)\d{s}$. By \cite[Lemma~6.8]{S:72}, we have that
\begin{equation}\label{thm:Tbounded_iff_Xinterpolation:eq3}
\K(f,t; \Lambda^1_\xi, L^\infty) \approx \int_0^{\Xi^{-1}(t)} f^*(s) \xi(s) \d{s} \quad \text{for every $f\in\Mpl(0,L)$ and $t\in(0,\Xi_L)$}.
\end{equation}
Let $a$ be defined by \eqref{prop:Roptimaldomain:L_or_1}. The embedding $X(0,L) \hookrightarrow (\Lambda^1_\xi + L^\infty)(0,L)$ follows from
\begin{align*}
\| f \|_{(\Lambda^1_\xi + L^\infty)(0,L)} &= \K(f,1; \Lambda^1_\xi, L^\infty) \leq \max\Big\{1,\frac1{\int_0^a \xi(s)\d{s}}\Big\} \K(f,\int_0^a \xi(s)\d{s}; \Lambda^1_\xi, L^\infty) \\
&\approx \int_0^a f^*(t) \xi(t) \d{t} \lesssim \frac1{\psi(a)} \int_0^a f^*(t) T_\varphi \chi_{(0,a)}(t) \d{t} \\
&\lesssim \|f\|_{X(0,L)}\|T_\varphi \chi_{(0,a)}\|_{X'(0,L)} \lesssim \|f\|_{X(0,L)},
\end{align*}
where we used H\"older's inequality \eqref{ch1:ri:holder} in the last but one inequality and the boundedness of $T_\varphi$ on $X'(0,L)$ in the last one. We now turn our attention to the embedding $\Lambda^1_\xi(0,L)\cap L^\infty(0,L) \hookrightarrow X(0,L)$. If $L<\infty$, the embedding is plainly true owing to \eqref{ch1:ri:smallestandlargestrispacefinitemeasure}. If $L=\infty$, it is sufficient to observe that, for every $f\in X'(0,\infty)$,  $f^*=g+h$ for some functions $g\in L^1(0,\infty)$ and $h\in M_\psi(0,\infty)$ thanks to \eqref{ch1:ri:XtoYiffY'toX'}, the fact that $(\Lambda^1_\xi(0,\infty)\cap L^\infty(0,\infty))' = L^1(0,\infty) + M_\psi(0,\infty)$ by \eqref{ch1:ri:dual_sum_and_inter}, and \eqref{ch1:ri:inclusion_is_always_continuous}. Set $g=f^*\chi_{(0,1)}$ and $h=f^*\chi_{(1,\infty)}$. Clearly, $g\in L^1(0,\infty)$ thanks to property (P5) of $\|\cdot\|_{X'(0,L)}$. Furthermore
\begin{align*}
\|h\|_{M_\psi(0,\infty)} &\approx \sup_{t\in(0,\infty)}\psi(t) (f^*\chi_{(1,\infty)})^*(t) \lesssim \sup_{t\in(0,\infty)}\psi(t+1) f^*(t+1) = \frac{\psi(1)}{\psi(1)}\sup_{t\in [1,\infty)}\psi(t) f^*(t) \\
&= \psi(1) T_\psi f(1) \approx T_\varphi f(1) < \infty,
\end{align*}
where we used the fact that $\xi$ satisfies the averaging condition~\eqref{averaging_condition} in the first equivalence (cf.~\cite[Lemma~2.1]{MO:19}) and \eqref{thm:Tbounded_iff_Xinterpolation:eq6} in the last one. Note that $T_\varphi f(1)$ is finite owing to \eqref{ch1:ri:XembeddedinM0} inasmuch as $T_\varphi f\in X'(0,\infty)$ and it is a nonincreasing function. Next, now that we know that $X(0,L)$ is an intermediate space between $\Lambda^1_\xi(0,L)$ and $L^\infty(0,L)$, in order to prove that (i) implies (ii), it remains to show that every admissible operator $S$ for the couple $\big(\Lambda^1_\xi(0,L), L^\infty(0,L)\big)$ is bounded on $X(0,L)$. Let $S$ be such an operator. Since $S$ is linear and bounded on both $\Lambda^1_\xi(0,L)$ and $L^\infty(0,L)$, it follows that (see~\cite[Chapter~5, Theorem~1.11]{BS})
\begin{equation}\label{thm:Tbounded_iff_Xinterpolation:eq4}
\K(Sf,t; \Lambda^1_\xi, L^\infty) \lesssim \K(f,t; \Lambda^1_\xi, L^\infty) \quad \text{for every $f\in(\Lambda^1_\xi + L^\infty)(0,L)$ and $t\in(0,L)$}.
\end{equation}
By combining \eqref{thm:Tbounded_iff_Xinterpolation:eq3} and \eqref{thm:Tbounded_iff_Xinterpolation:eq4}, we obtain that
\begin{equation}\label{thm:Tbounded_iff_Xinterpolation:eq5}
\int_0^t (Sf)^*(s) \xi(s) \d{s} \lesssim \int_0^t f^*(s) \xi(s) \d{s} \quad \text{for every $f\in(\Lambda^1_\xi + L^\infty)(0,L)$ and $t\in(0,L)$}.
\end{equation}
Since the function $(0,L)\ni t \mapsto \sup_{t\leq s < L}\varphi(s)g^*(s)$ is nonincreasing for every $g\in\Mpl(0,L)$, the Hardy lemma \eqref{ch1:ri:hardy-lemma} together with \eqref{thm:Tbounded_iff_Xinterpolation:eq5} implies that
\begin{equation*}
\int_0^L (Sf)^*(t) T_\varphi g(t) \d{t} \lesssim \int_0^L f^*(t) T_\varphi g(t) \d{t} \quad \text{for every $f\in(\Lambda^1_\xi + L^\infty)(0,L)$ and $g\in\Mpl(0,L)$}.
\end{equation*}
Therefore
\begin{align*}
\|Sf\|_{X(0,L)} &= \sup_{\substack{g\in{\Mpl(0,L)}\\\|g\|_{X'(0,L)}\leq1}}\int_0^{L}(Sf)^*(t)g^*(t)\d{t} \leq \sup_{\substack{g\in{\Mpl(0,L)}\\\|g\|_{X'(0,L)}\leq1}}\int_0^{L}(Sf)^*(t)T_\varphi g(t)\d{t} \\
&\lesssim \sup_{\substack{g\in{\Mpl(0,L)}\\\|g\|_{X'(0,L)}\leq1}}\int_0^{L}f^*(t)T_\varphi g(t)\d{t} \leq \|f\|_{X(0,L)} \sup_{\substack{g\in{\Mpl(0,L)}\\\|g\|_{X'(0,L)}\leq1}}\|T_\varphi g\|_{X'(0,L)} \\
&\lesssim \|f\|_{X(0,L)}
\end{align*}
for every $f\in X(0,L)$, where we used \eqref{ch1:ri:normX''down} in the equality, H\"older's inequality \eqref{ch1:ri:holder} in the last but one inequality and the boundedness of $T_\varphi$ on $X'(0,L)$ in the last one. Hence $S$ is bounded on $X(0,L)$.

We shall now prove that (iii) implies (i), whether $L$ is finite or infinite. Since $\xi$ is nonincreasing and satisfies the averaging condition~\eqref{averaging_condition}, it follows from \citep[Theorem~3.2]{GOP:06} (cf.~\cite[Lemma~3.1]{MO:19}) that $T_\varphi$ is bounded on $L^1(0,L)$. Furthermore, $T_\varphi$ is also bounded on $M_\psi(0,L)$, for
\begin{align*}
\|T_\varphi f\|_{M_\psi(0,L)} &= \sup_{t\in(0,L)} (T_\varphi f)^{**}(t)\psi(t) = \sup_{t\in(0,L)} \frac1{\int_0^t \xi(s) \d{s}}\int_0^t \xi(s) \sup_{s\leq \tau < L} \varphi(\tau)f(\tau) \d{s} \\
&\leq \sup_{t\in(0,L)}\varphi(t)f^*(t) \lesssim \|f\|_{M_\psi(0,L)},
\end{align*}
where we used \eqref{thm:Tbounded_iff_Xinterpolation:eq6} and \eqref{ch1:ri:twostarsdominateonestar} in the last inequality. Fix $f\in (L^1 + M_\psi)(0,L)$. We claim that
\begin{equation}\label{thm:Tbounded_iff_Xinterpolation:eq1}
\K(T_\varphi f, t; L^1, M_\psi)\lesssim \K(f, t; L^1, M_\psi) \quad \text{for every $t\in(0,\infty)$}
\end{equation}
with a multiplicative constant independent of $f$. Let $f=g+h$ with $g\in L^1(0, L)$ and $h\in M_\psi(0, L)$ be a decomposition of $f$. Note that the fact that $\xi$ is nonincreasing and satisfies the averaging condition~\eqref{averaging_condition} implies that
\begin{equation*}
\varphi(s)\lesssim \varphi\Big(\frac{s}{2}\Big) \quad \text{for every $s\in(0,L)$}.
\end{equation*}
Thanks to this and \eqref{ch1:ri:halfsubadditivityofonestar}, we have that
\begin{equation}\label{thm:Tbounded_iff_Xinterpolation:eq2}
\begin{aligned}
T_\varphi f(s)&\leq \frac1{\varphi(s)}\Big(\sup_{s\leq\tau < L}\varphi(\tau) g^*\Big(\frac{\tau}{2}\Big) + \sup_{s\leq\tau < L}\varphi(\tau) h^*\Big(\frac{\tau}{2}\Big)\Big) \\
&\lesssim T_\varphi g\Big(\frac{s}{2}\Big) + T_\varphi h\Big(\frac{s}{2}\Big)
\end{aligned}
\end{equation}
for every $s\in(0, L)$. By combining \eqref{thm:Tbounded_iff_Xinterpolation:eq2} and the boundedness of the dilation operator $D_2$ (see~\eqref{ch1:ri:dilation}) with the fact that $T_\varphi$ is bounded on both $L^1(0,L)$ and $M_\psi(0,L)$, we obtain that (cf.~\cite[p.~497]{BK:91})
\begin{align*}
\K(T_\varphi f, t; L^1, M_\psi) &\lesssim \K\Big(T_\varphi g\Big(\frac{\cdot}{2}\Big), t; L^1, M_\psi\Big) + \K\Big(T_\varphi h\Big(\frac{\cdot}{2}\Big), t; L^1, M_\psi\Big) \\
&\lesssim \|T_\varphi g\|_{L^1(0, L)} + t\|T_\varphi h\|_{M_\psi(0, L)} \lesssim \|g\|_{L^1(0, L)} + t\|h\|_{M_\psi(0, L)}
\end{align*}
for every $t\in(0,\infty)$, in which the multiplicative constants are independent of $f,g,h$ and $t$. Hence \eqref{thm:Tbounded_iff_Xinterpolation:eq1} is true. Now, since we have \eqref{thm:Tbounded_iff_Xinterpolation:eq1} at our disposal, there is a linear operator $S$ bounded on both $L^1(0,L)$ and $M_\psi(0,L)$ with norms that can be bounded from above by a constant independent of $f$ such that $Sf=T_\varphi f$ by virtue of \cite[Theorem~2]{CN:85}. Owing to (iii), $S$ is also bounded on $X'(0,L)$; moreover, its norm on $X'(0,L)$ can be bounded from above by a constant independent of $f$ (\cite[Chapter~3, Proposition~1.11]{BS}). Therefore
\begin{equation*}
\|T_\varphi f\|_{X'(0,L)} = \|Sf\|_{X'(0,L)} \lesssim \|f\|_{X'(0,L)},
\end{equation*}
in which the multiplicative constant is independent of $f$; hence $T_\varphi$ is bounded on $X'(0,L)$.

Finally, if $L<\infty$, then (ii) is equivalent to (iii); hence the three statements are equivalent to each other in this case. Indeed, since $(\Lambda^1_\xi+L^\infty)(0,L)=\Lambda^1_\xi(0,L)$ and $(L^1 + M_\psi)(0,L)=L^1(0,L)$ owing to \eqref{ch1:ri:smallestandlargestrispacefinitemeasure} and both $\Lambda^1_\xi(0,L)$ and $L^1(0,L)$ have absolutely continuous norm (in the sense of \cite[Chapter~1, Definition~3.1]{BS}), the equivalence of (ii) and (iii) follows from \cite[Corollary~3.6]{M:89}.
\end{proof}

A great deal of our effort has been devoted to describing optimal rearrangement\hyp{}invariant function spaces. A natural, somewhat related question is, can every rearrangement\hyp{}invariant function space be an optimal space? Suppose that $Z(0,L)$ is the optimal domain space for $H_{u, v, \nu}$ and $X(0,L)$, and denote by $W(0,L)$ the optimal target space for $H_{u, v, \nu}$ and $Z(0,L)$. Owing to the optimality of $W(0,L)$, we immediately see that $W(0,L)\hookrightarrow X(0,L)$. What is not obvious, however, is whether the opposite embedding, too, (is)/(can be) true. This leads us to the following theorem.
\begin{theorem}\label{thm:char_of_optimal_spaces} 
Let $\|\cdot\|_{X(0,L)}$ be a rearrangement\hyp{}invariant function norm. Let $\nu\colon(0,L)\to(0,L)$ be an increasing bijection. If $L=\infty$, assume that $\nu\in\Dinf{\infty}$. Let $u,v\colon(0,L)\to(0,\infty)$ be nonincreasing functions. Assume that $u$ is nondegenerate. If $L<\infty$, assume that $u(L^-)>0$ and $v(L^-)>0$. Let $\varrho$ be the functional defined by \eqref{prop:norminducedbyH:normdef}. The following three statements are equivalent.
\begin{enumerate}[(i)]
\item The space $X(0,L)$ is the optimal target space for the operator $H_{u, v, \nu}$ and some rearrangement\hyp{}invariant function space.
\item The space $X'(0,L)$ is the optimal domain space for the operator $R_{u, v, \nu^{-1}}$ and some rearrangement\hyp{}invariant function space.
\item We have that
\begin{equation}\label{thm:char_of_optimal_spaces:normonX'}
\|f\|_{X'(0,L)}\approx \sup_{\substack{g\in{\Mpl(0,L)}\\\varrho(g)\leq1}}\int_0^Lg(t)R_{u, v, \nu^{-1}}(f^*)(t)\d{t} \quad \text{for every $f\in\Mpl(0,L)$}.
\end{equation}
\end{enumerate}

If, in addition, 
\begin{equation*}
\frac1{v(t)}=\int_0^{\nu^{-1}(t)}\xi(s)\d{s} \quad \text{for every $t\in(0,L)$},
\end{equation*}
where $\xi\colon(0,L)\to(0,\infty)$ is a measurable function satisfying
\begin{equation}\label{thm:char_of_optimal_spaces:integral_cond_on_xi}
\frac{u(t)}{U(t)}\int_0^t \xi(s) \d{s} \lesssim \xi(t) \quad \text{for a.e.~$t\in(0,L)$},
\end{equation}
and the function $\varphi\circ\nu^{-1}$, where $\varphi=\frac{u}{\xi}$, is equivalent to a quasiconcave function, then each of the equivalent statements above implies that
\begin{enumerate}
\item[(iv)] the operator $T_\varphi$, defined by \eqref{opTdef}, is bounded on $X'(0,L)$.
\end{enumerate}
\end{theorem}
\begin{proof}
We start off by observing that each of the three equivalent statements implies that the functional $\varrho$ is actually a rearrangement\hyp{}invariant function norm. Statements (i) and (ii) imply it thanks to \myref{Proposition}{prop:norminducedbyH} and \myref{Proposition}{prop:Roptimalrange}, respectively. If we assume (iii), then, in particular, the set $\{g\in\Mpl(0,L)\colon\varrho(g)\leq1\}$ needs to contain a function $g\in\Mpl(0,L)$ not equal to $0$ a.e. It follows from \myref{Proposition}{prop:norminducedbyH} and its proof that $\varrho(g)=\infty$ for every $g\in\Mpl(0,L)$ not equal to $0$ a.e.~provided that $\varrho$ fails to be a rearrangement\hyp{}invariant function norm. Hence $\varrho$ is a rearrangement\hyp{}invariant function norm if (iii) is assumed. Therefore, in all of the cases, we are entitled to denote the corresponding rearrangement\hyp{}invariant function space over $(0,L)$ by $Z(0,L)$. Moreover, note that \eqref{thm:char_of_optimal_spaces:normonX'} actually reads as
\begin{equation}\label{thm:char_of_optimal_spaces:normonX'notationofnorms}
\|f\|_{X'(0,L)}\approx\|R_{u, v, \nu^{-1}}(f^*)\|_{Z'(0,L)}\quad\text{for every $f\in\Mpl(0,L)$}.
\end{equation}
 
Second, statements (i) and (ii) are clearly equivalent to each other owing to \myref{Remark}{rem:optimal_for_R_iff_for_H}.

Next, the proof of the fact that (ii) implies (iii) is based on the following important observation. If $X'(0,L)$ is the optimal domain space for the operator $R_{u, v, \nu^{-1}}$ and a rearrangement\hyp{}invariant function space $Y(0,L)$, then, in particular, $R_{u, v, \nu^{-1}}\colon X'(0,L)\to Y(0,L)$ is bounded. Consequently, by virtue of \myref{Proposition}{prop:Roptimalrange}, the rearragement\hyp{}invariant function space whose associate function norm is $\varrho$ is the optimal target space for the operator $R_{u, v, \nu^{-1}}$ and $X'(0,L)$. By \eqref{ch1:ri:X''=X}, this optimal target space is actually the space $Z'(0,L)$. Owing to \myref{Proposition}{prop:norminducedbyR}, the optimal domain space for the operator $R_{u, v, \nu^{-1}}$ and $Z'(0,L)$ exists, and we denote it by $W(0,L)$. Moreover,
\begin{equation}\label{thm:char_of_optimal_spaces:eq3}
\|f\|_{W(0,L)}\approx\|R_{u, v, \nu^{-1}}(f^*)\|_{Z'(0,L)}\quad\text{for every $f\in\Mpl(0,L)$}.
\end{equation}
The crucial observation is that we have, in fact, that $X'(0,L)=W(0,L)$. The embedding $X'(0,L)\hookrightarrow W(0,L)$ is valid because $R_{u, v, \nu^{-1}}\colon X'(0,L)\to Z'(0,L)$ is bounded and $W(0,L)$ is the optimal domain space for the operator $R_{u, v, \nu^{-1}}$ and $Z'(0,L)$. The validity of the opposite embedding is slightly more complicated. Since $R_{u, v, \nu^{-1}}\colon X'(0,L)\to Y(0,L)$ is bounded and $Z'(0,L)$ is the optimal target space for the operator $R_{u, v, \nu^{-1}}$ and $X'(0,L)$, we have that $Z'(0,L)\hookrightarrow Y(0,L)$. Consequently, since $R_{u, v, \nu^{-1}}\colon W(0,L)\to Z'(0,L)$ is bounded, so is $R_{u, v, \nu^{-1}}\colon W(0,L)\to Y(0,L)$. Using the fact that $X'(0,L)$ is the optimal domain space for the operator $R_{u, v, \nu^{-1}}$ and $Y(0,L)$, we obtain that $W(0,L)\hookrightarrow X'(0,L)$. Now that we know that $X'(0,L)=W(0,L)$, \eqref{thm:char_of_optimal_spaces:normonX'notationofnorms} follows from \eqref{thm:char_of_optimal_spaces:eq3}.

Next, note that (iii) implies (ii). Indeed, \eqref{thm:char_of_optimal_spaces:normonX'notationofnorms} coupled with \myref{Proposition}{prop:norminducedbyR} tells us that $X'(0,L)$ is the optimal domain space for the operator $R_{u, v, \nu^{-1}}$ and $Z'(0,L)$.

Finally, it only remains to prove that (iii) implies (iv) under the additional assumptions. By combining \eqref{thm:char_of_optimal_spaces:normonX'notationofnorms}  with the fact that $T_\varphi f$ is equivalent to a nonincreasing function for every $f\in\Mpl(0,L)$ (the multiplicative constants in this equivalence are independent of $f$), we have that
\begin{equation}\label{thm:char_of_optimal_spaces:eq7}
\begin{aligned}
\|T_\varphi f\|_{X'(0,L)}&\approx\|R_{u, v, \nu^{-1}}((T_\varphi f)^*)\|_{Z'(0,L)}\approx\|R_{u, v, \nu^{-1}}(T_\varphi f)\|_{Z'(0,L)}\\
&=\Big\|v(t)\int_0^{\nu^{-1}(t)} \xi(s) \esssup_{\tau\in[s,L)}\varphi(\tau)f^*(\tau)\d{s}\Big\|_{Z'(0,L)}\\
&\leq\Big\|v(t)\int_0^{\nu^{-1}(t)} \xi(s) \esssup_{\tau\in[s,\nu^{-1}(t))}\varphi(\tau)f^*(\tau)\d{s}\Big\|_{Z'(0,L)}\\
&\quad+\Big\|v(t) \Big(\esssup_{\tau\in[\nu^{-1}(t),L)}\varphi(\tau)f^*(\tau) \Big) \int_0^{\nu^{-1}(t)}\xi(s)\d{s}\Big\|_{Z'(0,L)}\\
&=\Big\|v(t)\int_0^{\nu^{-1}(t)} \xi(s) \esssup_{\tau\in[s,\nu^{-1}(t))}\varphi(\tau)f^*(\tau)\d{s}\Big\|_{Z'(0,L)}\\
&\quad+\Big\|\esssup_{\tau\in[\nu^{-1}(t),L)}\varphi(\tau)f^*(\tau)\Big\|_{Z'(0,L)}.
\end{aligned}
\end{equation}
Since $\varphi$ is equivalent to a continuous nondecreasing function and $\xi$ satisfies \eqref{thm:char_of_optimal_spaces:integral_cond_on_xi}, it follows from \cite[Theorem~3.2]{GOP:06} that
\begin{equation*}
\int_0^{\nu^{-1}(t)} \xi(s) \esssup_{\tau\in[s,\nu^{-1}(t))}\varphi(\tau)f^*(\tau)\d{s}\lesssim\int_0^{\nu^{-1}(t)}f^*(s) u(s) \d{s}\quad\text{for every $t\in(0,L)$};
\end{equation*}
hence
\begin{equation}\label{thm:char_of_optimal_spaces:eq4}
\Big\|v(t)\int_0^{\nu^{-1}(t)} \xi(s) \esssup_{\tau\in[s,\nu^{-1}(t))}\varphi(\tau)f^*(\tau)\d{s}\Big\|_{Z'(0,L)}\lesssim\|R_{u, v, \nu^{-1}}(f^*)\|_{Z'(0,L)}.
\end{equation}
Since the function $\varphi\circ\nu^{-1}$ is equivalent to a quasiconcave function, it follows from \cite[Lemma~4.10]{EMMP:20} that
\begin{equation*}
\Big\|\esssup_{\tau\in[\nu^{-1}(t),L)}\varphi(\tau)f^*(\tau)\Big\|_{Z'(0,L)}\lesssim\|\varphi(\nu^{-1}(t))f^*(\nu^{-1}(t))\|_{Z'(0,L)}.
\end{equation*}
We note that, although \cite[Lemma~4.10]{EMMP:20} deals only with the case $L=\infty$, its proof translates verbatim to the case of $L\in(0,\infty)$. Furthermore, we have that
\begin{equation}\label{thm:char_of_optimal_spaces:eq6}
\begin{aligned}
\|\varphi(\nu^{-1}(t))f^*(\nu^{-1}(t))\|_{Z'(0,L)}&\lesssim\Bigg\|\frac{U(\nu^{-1}(t))}{\int_0^{\nu^{-1}(t)}\xi(s)\d{s}}f^*(\nu^{-1}(t))\Bigg\|_{Z'(0,L)}\\
&\leq\|R_{u, v, \nu^{-1}}(f^*)\|_{Z'(0,L)},
\end{aligned}
\end{equation}
where we used the fact that $\xi$ satisfies \eqref{thm:char_of_optimal_spaces:integral_cond_on_xi} in the first inequality and the monotonicity of $f^*$ in the second one. By combining \eqref{thm:char_of_optimal_spaces:eq7} with \eqref{thm:char_of_optimal_spaces:eq4} and \eqref{thm:char_of_optimal_spaces:eq6} and using \eqref{thm:char_of_optimal_spaces:normonX'notationofnorms}, we obtain that
\begin{equation*}
\|T_\varphi f\|_{X'(0,L)}\lesssim\|R_{u, v, \nu^{-1}}(f^*)\|_{Z'(0,L)}\approx\|f\|_{X'(0,L)}\quad\text{for every $f\in\Mpl(0,L)$};
\end{equation*}
hence $T_\varphi$ is bounded on $X'(0,L)$.
\end{proof}

\begin{remarks}\ 
\begin{enumerate}[(i)]
\item If $X'(0,L)$ is the optimal domain space for $R_{u, v, \nu^{-1}}$ and some rearrangement\hyp{}invariant function space $Y(0,L)$, then $X'(0,L)$ is actually the optimal domain space for $R_{u, v, \nu^{-1}}$ and its own optimal target space. This follows from the following. Thanks to \myref{Proposition}{prop:Roptimalrange} and \myref{Proposition}{prop:norminducedbyR}, we are entitled to denote by $Z(0,L)$ the optimal target space for $R_{u, v, \nu^{-1}}$ and $X'(0,L)$ and by $W(0,L)$ the optimal domain space for $R_{u, v, \nu^{-1}}$ and $Z(0,L)$. We need to show that $X'(0,L)=W(0,L)$. On the one hand, since $R_{u, v, \nu^{-1}}\colon X'(0,L)\to Z(0,L)$ is bounded and $W(0,L)$ is the optimal domain space for $R_{u, v, \nu^{-1}}$ and $Z(0,L)$, we have that $X'(0,L)\hookrightarrow W(0,L)$. On the other hand, since $R_{u, v, \nu^{-1}}\colon X'(0,L)\to Y(0,L)$ is bounded and $Z(0,L)$ is the optimal target space for $R_{u, v, \nu^{-1}}$ and $X'(0,L)$, we have that $Z(0,L)\hookrightarrow Y(0,L)$; consequently, $R_{u, v, \nu^{-1}}\colon W(0,L)\to Y(0,L)$ is bounded. Finally, since $X'(0,L)$ is the optimal domain space for $R_{u, v, \nu^{-1}}$ and $Y(0,L)$, we obtain that $W(0,L)\hookrightarrow X'(0,L)$. Furthermore, by combining this observation with \myref{Remark}{rem:optimal_for_R_iff_for_H}, we also obtain that, if $X(0,L)$ is the optimal target space for $H_{u, v, \nu}$ and some rearrangement\hyp{}invariant function space $Y(0,L)$, then $X(0,L)$ is actually the optimal target space for $H_{u, v, \nu}$ and its own optimal domain space.

\item If $\xi$ satisfies the averaging condition~\eqref{averaging_condition}, then \eqref{thm:char_of_optimal_spaces:integral_cond_on_xi} is satisfied for every nonincreasing function $u$ inasmuch as $tu(t)\leq U(t)$ for every $t\in(0,L)$.

\item When $u(t)=t^{-1+\alpha}$, $v(t)=t^{-1+\beta}$,  and $\nu(t) = t^\gamma$, $t\in(0,L)$, the additional assumptions of \myref{Theorem}{thm:char_of_optimal_spaces} are satisfied if $\alpha\in(0,1]$, $\beta\in[0,1)$, $\gamma>0$ and $1 \leq \frac{\alpha}{\gamma} + \beta \leq 2$.
\end{enumerate}
\end{remarks}

The remainder of this subsection is devoted to the particular but important case $u\equiv1$, in which we will be able to establish an even stronger connection between the various notions that we have met. First, we need to equip ourselves with the following auxiliary result, which generalizes \cite[Lemma~4.9]{EMMP:20} and whose immediate corollary for $u\equiv1$ (namely \myref{Corollary}{cor:Honsimple:u==1}) is of independent interest.
\begin{proposition}\label{prop:Honsimple}
Let $\|\cdot\|_{X(0,L)}$ be a rearrangement\hyp{}invariant function norm. Let $\nu\colon(0,L)\to(0,L)$ be an increasing bijection. Assume that $\nu^{-1}\in\Dsup{0}$. If $L=\infty$, assume that $\nu^{-1}\in\Dsup{\infty}$. Let $u,v\colon(0,L)\to(0,\infty)$ be nonincreasing. Assume that $v$ satisfies the averaging condition \eqref{averaging_condition}, and denote its averaging constant by $C$.
Set $f=\sum_{i=1}^Nc_i\chi_{(0,a_i)}$, where $c_i\in(0,\infty)$, $i=1,\dots, N$, and $0<a_1<\cdots<a_N < L$. We have that
\begin{equation}\label{prop:Honsimple:equivalence}
\Big\|u(t)\int_{\nu(t)}^Lf(s)v(s)\d{s}\Big\|_{X(0,L)}\approx\Big\|u(t)\sum_{i=1}^Na_ic_iv(a_i)\chi_{(0,\nu^{-1}(a_i))}(t)\Big\|_{X(0,L)},
\end{equation}
in which the multiplicative constants depend only on $\nu$ and $C$.
\end{proposition}
\begin{proof}
First, observe that $\inf_{t\in(0,L)}\frac{\nu^{-1}(\frac{t}{\theta})}{\nu^{-1}(t)}\in(0,1)$, where $\theta>1$ is such that $\nu^{-1}\in\Dsup[\theta]{0}$ and, if $L=\infty$, also $\nu^{-1}\in\Dsup[\theta]{\infty}$. We denote the infimum by $M$.

Second, we have that
\begin{align}
\Big\|u(t)\int_{\nu(t)}^Lf(s)v(s)\d{s}\Big\|_{X(0,L)}&=\Big\|u(t)\sum_{i=1}^Nc_i\chi_{(0,\nu^{-1}(a_i))}(t)\int_{\nu(t)}^{a_i}v(s)\d{s}\Big\|_{X(0,L)}\notag\\
\begin{split}\label{prop:Honsimple:eq1}
&\geq\Big\|u(t)\sum_{i=1}^Nc_i\chi_{(0,\nu^{-1}(\frac{a_i}{\theta}))}(t)v(a_i)(a_i-\nu(t))\Big\|_{X(0,L)}\\
&\geq\frac{\theta-1}{\theta}\Big\|u(t)\sum_{i=1}^Nc_i\chi_{(0,\nu^{-1}(\frac{a_i}{\theta}))}(t)v(a_i)a_i\Big\|_{X(0,L)}\\
&\geq\frac{\theta-1}{\theta}\Big\|u(t)\sum_{i=1}^Nc_i\chi_{(0,M\nu^{-1}(a_i))}(t)v(a_i)a_i\Big\|_{X(0,L)}\\
&\geq M\frac{\theta-1}{\theta}\Big\|u(Mt)\sum_{i=1}^Nc_i\chi_{(0,\nu^{-1}(a_i))}(t)v(a_i)a_i\Big\|_{X(0,L)}\\
&\geq M\frac{\theta-1}{\theta}\Big\|u(t)\sum_{i=1}^Nc_i\chi_{(0,\nu^{-1}(a_i))}(t)v(a_i)a_i\Big\|_{X(0,L)}
\end{split}
\end{align}
thanks to the fact that $u$ and $v$ are nonincreasing and \eqref{ch1:ri:dilation} (the boundedness of the dilation operator $D_{\frac1{M}}$).

Last, using the fact that $v$ satisfies the averaging condition~\eqref{averaging_condition}, we obtain that
\begin{equation}\label{prop:Honsimple:eq2}
\begin{aligned}
\Big\|u(t)\int_{\nu(t)}^Lf(s)v(s)\d{s}\Big\|_{X(0,L)}&=\Big\|u(t)\sum_{i=1}^Nc_i\chi_{(0,\nu^{-1}(a_i))}(t)\int_{\nu(t)}^{a_i}v(s)\d{s}\Big\|_{X(0,L)}\\
&\leq C\Big\|u(t)\sum_{i=1}^Nc_i\chi_{(0,\nu^{-1}(a_i))}(t)a_iv(a_i)\Big\|_{X(0,L)}.
\end{aligned}
\end{equation}
By combining \eqref{prop:Honsimple:eq2} and \eqref{prop:Honsimple:eq1}, we obtain \eqref{prop:Honsimple:equivalence}.
\end{proof}

Since every nonnegative, nonincreasing function on $(0,L)$ is the pointwise limit of a nondecreasing sequence of nonnegative, nonincreasing simple functions, \myref{Proposition}{prop:Honsimple} with $u\equiv1$ has the following important corollary.
\begin{corollary}\label{cor:Honsimple:u==1}
Let $\nu\colon(0,L)\to(0,L)$ be an increasing bijection. Assume that $\nu^{-1}\in\Dsup{0}$. If $L=\infty$, assume that $\nu^{-1}\in\Dsup{\infty}$. Let $v\colon(0,L)\to(0, \infty)$ be a nonincreasing function satisfying the averaging condition~\eqref{averaging_condition}. Let $f\in\Mpl(0,L)$. There is a nondecreasing sequence $\{f_k\}_{k=1}^\infty$ of nonnegative, nonincreasing simple functions on $(0,L)$ such that, for every  rearrangement\hyp{}invariant function norm $\|\cdot\|_{X(0,L)}$,
\begin{equation*}
\lim_{k\to\infty}\|H_{1, v, \nu}(f_k)\|_{X(0,L)}\approx\|f^*\|_{X(0,L)}=\|f\|_{X(0,L)},
\end{equation*}
in which the multiplicative constants depend only on $\nu$ and the averaging constant of $v$.
\end{corollary}

\begin{remark}
The assumption $\nu^{-1}\in\Dsup{0}$ is not overly restrictive. For example, it is satisfied whenever $\nu$ is equivalent to $t\mapsto t^\alpha \ell_1(t)^{\beta_1}\cdots\ell_k(t)^{\beta_k}$ near $0$ for any $\alpha>0$, $k\in\N_0$ and $\beta_j\in\R$, $j=1, 2, \dots, k$, where the functions $\ell_j$ are $j$-times iterated logarithmic functions defined as
\begin{equation*}
\ell_j(t)=
	\begin{cases}
		1 + |\log t| \quad &\text{if $j=1$},\\
		1 + \log\ell_{j-1}(t) \quad &\text{if $j>1$},
	\end{cases}
\end{equation*}
for $t\in(0,L)$. In this case, $\nu^{-1}$ is equivalent to $t\mapsto t^{\frac1{\alpha}} \ell_1(t)^{-\frac{\beta_1}{\alpha}}\cdots\ell_k(t)^{-\frac{\beta_k}{\alpha}}$ near $0$ (cf.~\citep[Appendix~5]{BGT:89}). On the other hand, $\nu(t)=\log^\alpha(\frac1{t})$ near $0$, where $\alpha<0$, is a typical example of a function not satisfying the assumption. The same remark (with the obvious modifications) is true for the assumption $\nu^{-1}\in\Dsup{\infty}$.
\end{remark}

While \myref{Proposition}{prop:norminducedbyT} provides a sufficient condition for simplification of \eqref{prop:norminducedbyH:normdef}, the following proposition provides a necessary one.
\begin{proposition}\label{prop:normH_withu=1_simplified_implies_X_optimal_space}
Let $\|\cdot\|_{X(0,L)}$ be a rearrangement\hyp{}invariant function norm. Let $\nu\colon(0,L)\to(0,L)$ be an increasing bijection. Assume that $\nu^{-1}\in\Dsup{0}$. If $L=\infty$, assume that $\nu^{-1}\in\Dsup{\infty}$ and $\nu\in\Dinf{\infty}$. Let $v\colon (0,L) \to (0,\infty)$ be a nonincreasing function satisfying the averaging condition~\eqref{averaging_condition}. If there is a positive constant $C$ such that
\begin{equation}\label{prop:normH_withu=1_simplified_implies_X_optimal_space:simplifiednorm}
\sup_{h\sim f}\|H_{1, v, \nu}h\|_{X(0,L)}\leq C\|H_{1, v, \nu}(f^*)\|_{X(0,L)}\quad\text{for every $f\in\Mpl(0,L)$},
\end{equation}
then the three equivalent statements from \myref{Theorem}{thm:char_of_optimal_spaces} with $u\equiv1$ are satisfied.
\end{proposition}
\begin{proof}
Let $a$ be defined by \eqref{prop:Roptimaldomain:L_or_1}. Since $v$ is integrable over $(0,a)$, for it satisfies the averaging condition~\eqref{averaging_condition}, we have that
\begin{equation}\label{prop:normH_withu=1_simplified_implies_X_optimal_space:eq3}
\Big\|\chi_{(0,\nu^{-1}(a))}(t)\int_{\nu(t)}^av(s)\d{s}\Big\|_{X(0,L)}\leq\int_0^av(s)\d{s}\|\chi_{(0,\nu^{-1}(a))}\|_{X(0,L)}<\infty.
\end{equation}
Furthermore, if $L=\infty$, then $\limsup_{\tau\to\infty}v(\tau)\|\chi_{(0,\nu^{-1}(\tau))}\|_{X(0,\infty)}<\infty$. Indeed, suppose that $\limsup_{\tau\to\infty}v(\tau)\|\chi_{(0,\nu^{-1}(\tau))}\|_{X(0,\infty)}=\infty$. It follows from the proof of \myref{Proposition}{prop:norminducedbyH} that 
\begin{equation}\label{prop:normH_withu=1_simplified_implies_X_optimal_space:eq4}
\sup_{h\sim \chi_{(0,1)}}\|H_{1, v, \nu}h\|_{X(0,\infty)} = \infty.
\end{equation}
However, since
\begin{equation*}
\sup_{h\sim \chi_{(0,1)}}\|H_{1, v, \nu}h\|_{X(0,\infty)} \approx \|H_{1, v, \nu}\chi_{(0,1)}\|_{X(0,\infty)} = \Big\|\chi_{(0,\nu^{-1}(1))}(t)\int_{\nu(t)}^1v(s)\d{s}\Big\|_{X(0,\infty)}<\infty
\end{equation*}
thanks to \eqref{prop:normH_withu=1_simplified_implies_X_optimal_space:simplifiednorm} and \eqref{prop:normH_withu=1_simplified_implies_X_optimal_space:eq3}, \eqref{prop:normH_withu=1_simplified_implies_X_optimal_space:eq4} is not possible. Hence, \myref{Proposition}{prop:norminducedbyH} guarantees that the optimal domain space for $H_{1, v, \nu}$ and $X(0,L)$ exists; moreover, if we denote it by $Z(0,L)$, then \eqref{prop:normH_withu=1_simplified_implies_X_optimal_space:simplifiednorm} implies that
\begin{equation}\label{prop:normH_withu=1_simplified_implies_X_optimal_space:Znorm}
\|f\|_{Z(0,L)} \approx \| H_{1, v, \nu}(f^*) \|_{X(0,L)} \quad \text{for every $f\in\Mpl(0,L)$}.
\end{equation}

Now, we finally turn our attention to proving that \eqref{prop:normH_withu=1_simplified_implies_X_optimal_space:simplifiednorm} implies statement (iii) from \myref{Theorem}{thm:char_of_optimal_spaces}. Let $Y(0,L)$ be the optimal target space for the operator $H_{1, v, \nu}$ and $Z(0,L)$. Its existence is guaranteed by \myref{Proposition}{prop:Hoptimalrange}, and we have that 
\begin{equation}\label{prop:normH_withu=1_simplified_implies_X_optimal_space:eq1}
\|f\|_{Y'(0,L)}=\|R_{1, v, \nu^{-1}}(f^*)\|_{Z'(0,L)}\quad\text{for every $f\in\Mpl(0,L)$}.
\end{equation}
Using the optimality of $Y(0,L)$ combined with the fact that $H_{1, v, \nu}\colon Z(0,L) \to X(0,L)$ is bounded (and so $Y(0,L)\hookrightarrow X(0,L)$), and \eqref{prop:normH_withu=1_simplified_implies_X_optimal_space:Znorm}, we obtain that
\begin{equation*}
\|H_{1, v, \nu}(f^*)\|_{X(0,L)} \lesssim \|H_{1, v, \nu}(f^*)\|_{Y(0,L)} \lesssim \|f\|_{Z(0,L)} \approx \|H_{1, v, \nu}(f^*)\|_{X(0,L)}
\end{equation*}
for every $f\in\Mpl(0,L)$; hence
\begin{equation*}
\|H_{1, v, \nu}(f^*)\|_{X(0,L)} \approx \|H_{1, v, \nu}(f^*)\|_{Y(0,L)}\quad\text{for every $f\in\Mpl(0,L)$}.
\end{equation*}
In particular, we have that
\begin{equation}\label{prop:normH_withu=1_simplified_implies_X_optimal_space:eq2}
\|H_{1, v, \nu}h\|_{X(0,L)}\approx\|H_{1, v, \nu}h\|_{Y(0,L)}
\end{equation}
for every nonincreasing simple function $h\in\Mpl(0,L)$. By combining \eqref{prop:normH_withu=1_simplified_implies_X_optimal_space:eq2} with \myref{Corollary}{cor:Honsimple:u==1}, we obtain that
\begin{equation*}
\|f^*\|_{X(0,L)}\approx\|f^*\|_{Y(0,L)}\quad\text{for every $f\in\Mpl(0,L)$}.
\end{equation*}
Owing to the rearrangement invariance of both function norms, it follows that $X(0,L)=Y(0,L)$. Hence \eqref{thm:char_of_optimal_spaces:normonX'} with $u\equiv1$ follows from \eqref{prop:normH_withu=1_simplified_implies_X_optimal_space:eq1} combined with \eqref{ch1:ri:normX'}.
\end{proof}

We obtain the final result of this subsection by combining \myref{Theorem}{thm:char_of_optimal_spaces}, \myref{Proposition}{prop:normH_withu=1_simplified_implies_X_optimal_space}, \myref{Proposition}{prop:norminducedbyT} and \myref{Theorem}{thm:Tbounded_iff_Xinterpolation}.
\begin{theorem}\label{thm:charTbounded}
Let $\|\cdot\|_{X(0,L)}$ be a rearrangement\hyp{}invariant function norm. Let $\nu\colon(0,L)\to(0,L)$ be an increasing bijection. Assume that $\nu^{-1}\in\Dsup{0}$. If $L=\infty$, assume that $\nu^{-1}\in\Dsup{\infty}$ and $\nu\in\Dinf{\infty}$. Let $v\colon(0,L)\to(0,\infty)$ be defined by
\begin{equation*}
\frac1{v(t)}=\int_0^{\nu^{-1}(t)}\xi(s)\d{s},\ t\in(0,L),
\end{equation*}
where $\xi\colon(0,L)\to(0,\infty)$ is a measurable function satisfying the averaging condition~\eqref{averaging_condition}. Assume that $v$, too, satisfies the averaging condition~\eqref{averaging_condition}. Set $\varphi=\frac1{\xi}$. Assume that the function $\varphi\circ\nu^{-1}$ is equivalent to a quasiconcave function. Let $\varrho$ be the functional defined by \eqref{prop:norminducedbyH:normdef} with $u\equiv1$. The following five statements are equivalent.
\begin{enumerate}[(i)]
\item The operator $T_\varphi$, defined by \eqref{opTdef}, is bounded on $X'(0,L)$.
\item There is a positive constant $C$ such that
\begin{equation*}
\sup_{h\sim f}\|H_{1, v, \nu}h\|_{X(0,L)}\leq C\|H_{1, v, \nu}(f^*)\|_{X(0,L)}\quad\text{for every $f\in\Mpl(0,L)$}.
\end{equation*}
\item The space $X(0,L)$ is the optimal target space for the operator $H_{1, v, \nu}$ and some rearrangement\hyp{}invariant function space.
\item The space $X'(0,L)$ is the optimal domain space for the operator $R_{1, v, \nu^{-1}}$ and some rearrangement\hyp{}invariant function space.
\item We have that
\begin{equation*}
\|f\|_{X'(0,L)} \approx \sup_{\substack{g\in{\Mpl(0,L)}\\\varrho(g)\leq1}}\int_0^Lg(t)R_{1, v, \nu^{-1}}(f^*)(t)\d{t} \quad \text{for every $f\in\Mpl(0,L)$}.
\end{equation*}
\end{enumerate}
If $L<\infty$, these five statements are also equivalent to
\begin{enumerate}[(i)]
\item[(vi)] $X(0,L) \in \Int\big(\Lambda^1_\xi(0,1), L^\infty(0,1)\big)$.
\end{enumerate}
\end{theorem}

\begin{remarks}\hphantom{}
\begin{itemize}
\item The assumption that $v$ satisfies the averaging condition~\eqref{averaging_condition} is natural because it forbids weights $v$ for which the question of whether $X(0,L)$ (or $X'(0,L)$) is the optimal target (or domain) space for $H_{1, v, \nu}$ (or $R_{1, v, \nu^{-1}}$) and some rearrangement\hyp{}invariant function space cannot be decided by the boundedness of the corresponding supremum operator $T_\varphi$. This can be illustrated by a very simple example. Consider $\nu=\id$ and $\xi\equiv1$. Since $T_\varphi f=f^*$, $T_\varphi$ is bounded on any $X'(0,L)$; however, $R_{1, v,\nu^{-1}}f(t)=\frac1{t}\int_0^t |f(s)|\d{s}$ clearly need not be bounded from $X'(0,L)$ to $(L^1+L^\infty)(0,L)$, which is the largest rearrangement\hyp{}invariant function space. To this end, consider, for example, $X(0,L)=L^\infty(0,\infty)$ (cf.~\cite[Proposition~4.1]{ST:16}).

\item When $v(t)=t^{-1+\beta}$ and $\nu(t) = t^\gamma$, $t\in(0,L)$, the assumptions of \myref{Theorem}{thm:charTbounded} are satisfied if $\beta\in(0,1)$, $\gamma>0$ and $1 \leq \frac1{\gamma} + \beta \leq 2$.
\end{itemize}
\end{remarks}

\subsection{Iteration of optimal r.i.~norms}\label{sec:iteration}
This subsection is devoted to so-called \emph{sharp iteration principles} for the operators $R_{u, v, \nu}$ and $H_{u, v, \nu}$. To illustrate their meaning and importance, suppose that $Y_1(0,L)$ is the optimal target space for $H_{u_1,v_1,\nu_1}$ and a rearrangement\hyp{}invariant function space $X(0,L)$. Let us now go one step further and suppose that $Y_2(0,L)$ is the optimal target space for $H_{u_2,v_2,\nu_2}$ and $Y_1(0,L)$. In the light of \myref{Proposition}{prop:Hoptimalrange}, the associate function norm of $\|\cdot\|_{Y_2(0,L)}$ is equal to $\|f\|_{Y_2'(0,L)}=\|R_{u_1,v_1,\nu_1^{-1}}((R_{u_2,v_2,\nu_2^{-1}}(f^*))^*)\|_{X'(0,L)}$. We immediately see that there is an inevitable difficulty that we face if we wish to understand the iterated norm. This difficulty is caused by the fact that the function $R_{u_2,v_2,\nu_2^{-1}}(f^*)$ is hardly ever (equivalent to) a nonincreasing function (unless $u_2$, $v_2$ and $\nu_2$ are related to each other in a very specific way; see \myref{Proposition}{prop:norminducedbyHwhenRnonincreasing}), and so we cannot just readily ``delete'' the outer star. Nevertheless, with some substantial effort, we shall be able to equivalently express the iterated norm as a noniterated one under suitable assumptions. The suitable assumptions are such that the iteration does not lead to the presence of kernels, which would go beyond the scope of this paper (see~\cite[Section~8]{CPS:15} in that regard). It should be noted that such iteration is not artificial. For example, it is an essential tool for establishing sharp iteration principles for various Sobolev embeddings, which, roughly speaking, ensure that the optimal rearrangement\hyp{}invariant target space in a Sobolev embedding of $(k+l)$th order is the same as that obtained by composing the optimal Sobolev embedding of order $k$ with the optimal Sobolev embedding of order $l$ (see~\cite{CP:16, CPS:20, M:21} and references therein). Another possible application is description of optimal rearrangement\hyp{}invariant function norms for compositions of some operators of harmonic analysis (see~\cite{EMMP:20} and references therein for optimal behavior of some classical operators on rearrangement\hyp{}invariant function spaces). Finally, the motivation behind studying function norms induced by $H_{u_1, v_1, \nu_1}\circ H_{u_1, v_1, \nu_1}$ is similar.

The following proposition is the first step towards the sharp iteration principle for $R_{u, v, \nu}$.
\begin{proposition}\label{prop:Riterationoptimaldomainlowerbound}
Let $\|\cdot\|_{X(0,L)}$ be a rearrangement\hyp{}invariant function norm. Let $\nu_1,\nu_2\colon(0,L)\to(0,L)$ be increasing bijections. Assume that $\nu_2\in\Dsup{0}$. If $L=\infty$, assume that $\nu_2\in\Dsup{\infty}$. Set $\nu=\nu_2\circ\nu_1$. Let $u_1, u_2\colon(0,L)\to(0,\infty)$ be nonincreasing. Let $v_1\colon(0,L)\to(0,\infty)$ be measurable. Let $v_2\colon(0,L)\to(0,\infty)$ be a nonincreasing function satisfying the averaging condition~\eqref{averaging_condition}.  We have that
\begin{equation*}
\Big\|u_1(\nu_1(t))v_1(t)\nu_1(t)v_2(\nu_1(t))\int_0^{\nu(t)}f^*(s) u_2(s)\d{s}\Big\|_{X(0,L)}\lesssim\|R_{u_1, v_1, \nu_1}((R_{u_2, v_2, \nu_2}(f^*))^*)\|_{X(0,L)}
\end{equation*}
for every $f\in\Mpl(0,L)$, in which the multiplicative constant depends only on $\nu_2$ and the averaging constant of $v_2$.
\end{proposition}
\begin{proof}
Note that $\inf_{t\in(0,L)}\frac{\nu_2(\frac{t}{\theta})}{\nu_2(t)}>0$, where $\theta>1$ is such that $\nu_2\in\Dsup[\theta]{0}$ and, if $L=\infty$, also $\nu_2\in\Dsup[\theta]{\infty}$. Consequently, there is $N\in\N$, such that $\nu_2(t)\leq N\nu_2(\frac{t}{\theta})$ for every $t\in(0,L)$. Hence, for every $f\in\Mpl(0,L)$, we have that
\begin{equation}\label{prop:Riterationoptimaldomainlowerbound:eq1}
\int_0^{\nu_2(t)}f^*(s)u_2(s)\d{s}\leq N\int_0^{\nu_2(\frac{t}{\theta})}f^*(s)u_2(s)\d{s} \quad \text{for every $t\in(0,L)$}
\end{equation}
owing to the fact that $f^* u_2$ is nonincreasing. Thanks to the monotonicity of $u_1$ and $v_2$, the fact that $v_2$ satisfies the averaging condition~\eqref{averaging_condition} and the inequality \eqref{prop:Riterationoptimaldomainlowerbound:eq1}, we have that
\begin{align*}
&\Big\|v_1(t)u_1(\nu_1(t))\nu_1(t)v_2(\nu_1(t))\int_0^{\nu(t)}f^*(s)u_2(s)\d{s}\Big\|_{X(0,L)}\\
&\lesssim\Big\|v_1(t)u_1(\nu_1(t))\int_0^{\nu_1(t)}v_2(s)\d{s}\,\int_0^{\nu(t)}f^*(s)u_2(s)\d{s}\Big\|_{X(0,L)}\\
&\approx\Big\|v_1(t)u_1(\nu_1(t))\int_{\frac{\nu_1(t)}{\theta}}^{\nu_1(t)}v_2(s)\d{s}\,\int_0^{\nu(t)}f^*(s)u_2(s)\d{s}\Big\|_{X(0,L)}\\
&\leq\Big\|v_1(t)\int_{\frac{\nu_1(t)}{\theta}}^{\nu_1(t)}v_2(s)u_1(s)\d{s}\,\int_0^{\nu(t)}f^*(s)u_2(s)\d{s}\Big\|_{X(0,L)}\\
&\lesssim\Big\|v_1(t)\int_{\frac{\nu_1(t)}{\theta}}^{\nu_1(t)}v_2(s)u_1(s)\d{s}\,\int_0^{\nu_2(\frac{\nu_1(t)}{\theta})}f^*(s)u_2(s)\d{s}\Big\|_{X(0,L)}\\
&\leq\Big\|v_1(t)\int_{\frac{\nu_1(t)}{\theta}}^{\nu_1(t)}\Big(v_2(s)\int_0^{\nu_2(s)}f^*(\tau)u_2(\tau)\d{\tau}\Big)u_1(s)\d{s}\Big\|_{X(0,L)}\\
&\leq\Big\|v_1(t)\int_0^{\nu_1(t)}\Big(v_2(s)\int_0^{\nu_2(s)}f^*(\tau)u_2(\tau)\d{\tau}\Big)u_1(s)\d{s}\Big\|_{X(0,L)}\\
&\leq\Big\|R_{u_1,v_1,\nu_1}((R_{u_2,v_2,\nu_2}(f^*))^*)\Big\|_{X(0,L)}
\end{align*}
for every $f\in\Mpl(0,L)$, where we used the Hardy--Littlewood inequality \eqref{ch1:ri:HL} in the last inequality.
\end{proof}

We are now in a position to establish the sharp iteration principle for $R_{u,v,\nu}$.
\begin{theorem}\label{thm:Riterationoptimaldomainupperbound}
Let $\|\cdot\|_{X(0,L)}$ be a rearrangement\hyp{}invariant function norm. Let $\nu_1,\nu_2\colon(0,L)\to(0,L)$ be increasing bijections. Assume that $\nu_2\in\Dsup{0}$. If $L=\infty$, assume that $\nu_2\in\Dsup{\infty}$. Let $u_1, u_2\colon(0,L)\to(0,\infty)$ be nonincreasing. Let $v_1\colon(0,L)\to(0,\infty)$ be a continuous function. Let $v_2\colon(0,L)\to(0,\infty)$ be defined by
\begin{equation*}
\frac1{v_2(t)}=\int_0^{\nu_2(t)}\xi(s)\d{s},\ t\in(0,L),
\end{equation*}
where $\xi\colon(0,L)\to(0,\infty)$ is a measurable function. Assume that the function $u_1 v_2$ satisfies the averaging condition~\eqref{averaging_condition}. Set $\nu=\nu_2\circ\nu_1$, 
\begin{align*}
v(t)&=\nu_1(t)u_1(\nu_1(t))v_1(t)v_2(\nu_1(t)),\ t\in(0,L),\\
\intertext{and}
\eta(t)&=\frac1{U_2(t)v(\nu^{-1}(t))},\ t\in(0,L).
\end{align*}
Assume that $\eta$ and $\frac{\eta}{\xi}$ are equivalent to nonincreasing functions. Furthermore, assume that there are positive constants $C_1$ and $C_2$ such that
\begin{align}
\int_0^t \eta(s) u_2(s) \d{s} &\leq C_1 U_2(t) \eta(t) \quad \text{for a.e.~$t\in(0,L)$} \label{thm:Riterationoptimaldomainupperbound:integralassumption_on_eta_times_u2} \\
\intertext{and}
\frac1{t}\int_0^t U_2(\nu(s)) v(s) \d{s} &\geq C_2 U_2(\nu(t)) v(t) \quad \text{for a.e.~$t\in(0,L)$}. \label{thm:Riterationoptimaldomainupperbound:integralassumptionon1overeta}
\end{align}
We have that
\begin{equation*}
\|R_{u_1, v_1, \nu_1}((R_{u_2, v_2, \nu_2}(f^*))^*)\|_{X(0,L)}\approx\|R_{u_2, v, \nu}(f^*)\|_{X(0,L)}
\end{equation*}
for every $f\in\Mpl(0,L)$, in which the multiplicative constants depend only on $\nu_1$, $\nu_2$, $C_1$, $C_2$, the averaging constant of $u_1 v_2$ and the multiplicative constants in the equivalences of $\eta$ and $\frac{\eta}{\xi}$ to nonincreasing functions.
\end{theorem}
\begin{proof}
First, since the fact that $v_2u_1$ satisfies the averaging condition~\eqref{averaging_condition} together with the monotonicity of $u_1$ implies that $v_2$, too, satisfies the averaging condition~\eqref{averaging_condition} (with the same multiplicative constant), we have that
\begin{equation*}
\|R_{u_1, v_1, \nu_1}((R_{u_2, v_2, \nu_2}(f^*))^*)\|_{X(0,L)}\gtrsim\|R_{u_2, v, \nu}(f^*)\|_{X(0,L)}\quad\text{for every $f\in\Mpl(0,L)$}
\end{equation*}
thanks to \myref{Proposition}{prop:Riterationoptimaldomainlowerbound}; consequently we only need to prove the opposite inequality.

We may assume that $u_2$ is nondegenerate and $\psi\in X(0,L)$, where $\psi$ is defined as $\psi(t)=v(t)U_2(\nu(t))\chi_{(0,L)}(t)+v(t)\chi_{(L,\infty)}(t)$, $t\in(0,L)$, for, if it is not the case, then $\|R_{u_2, v, \nu}(f^*)\|_{X(0,L)}=\infty$ for every $f\in\Mpl(0,L)$ that is not equivalent to $0$ a.e. \myref{Proposition}{prop:norminducedbyR} with $u=u_2$ guarantees that there is a rearrangement\hyp{}invariant function space $Z(0,L)$ such that
\begin{equation*}
\|f\|_{Z(0,L)}=\|R_{u_2, v, \nu}(f^*)\|_{X(0,L)}\quad\text{for every $f\in\Mpl(0,L)$}.
\end{equation*}
Furthermore, by \eqref{prop:RXtoYbddiffHY'toX':normRXtoZ=normHZ'toX'} and the Hardy--Littlewood inequality \eqref{ch1:ri:HL}, we have that
\begin{equation}\label{thm:Riterationoptimaldomainupperbound:eq1}
\sup_{\|g\|_{X'(0,L)}\leq1}\|H_{u_2, v, \nu^{-1}}g\|_{Z'(0,L)}=1.
\end{equation}
Note that, for every $f\in\Mpl(0,L)$, the function
\begin{equation*}
(0,L)\ni t\mapsto v_2(t)\int_0^{\nu_2(t)}\xi(s) u_2(s) \sup_{\tau\in[s,L)}\frac1{\xi(\tau)}f^*(\tau)\d{s}
\end{equation*}
is nonincreasing because it is the integral mean of the nonincreasing function $(0,L)\ni s \mapsto u_2(s) \sup_{\tau\in[s,L)}\frac1{\xi(\tau)}f^*(\tau)$ over the interval $(0,\nu_2(t))$ with respect to the measure $\xi(s)\d{s}$. By \eqref{ch1:ri:normX''} and \eqref{RaHdual}, we have that
\begin{align*}
&\|R_{u_1, v_1, \nu_1}((R_{u_2, v_2, \nu_2}(f^*))^*)\|_{X(0,L)}=\sup_{\|g\|_{X'(0,L)}\leq1}\int_0^L(R_{u_2, v_2, \nu_2}(f^*))^*(t)H_{u_1 v_1, \nu_1^{-1}}g(t)\d{t}\\
&=\sup_{\|g\|_{X'(0,L)}\leq1}\int_0^L\Big[v_2(s)\int_0^{\nu_2(s)} u_2(\tau) f^*(\tau)\d{\tau}\Big]^*(t)H_{u_1, v_1, \nu_1^{-1}}g(t)\d{t}\\
&\leq\sup_{\|g\|_{X'(0,L)}\leq1}\int_0^L\Big[v_2(s)\int_0^{\nu_2(s)} \xi(\tau) u_2(\tau) \sup_{x\in[\tau,L)}\frac1{\xi(x)}f^*(x)\d{\tau}\Big]^*(t)H_{u_1, v_1, \nu_1^{-1}}g(t)\d{t}\\
&=\sup_{\|g\|_{X'(0,L)}\leq1}\int_0^Lv_2(t)\int_0^{\nu_2(t)} \xi(s) u_2(s) \sup_{\tau\in[s,L)}\frac1{\xi(\tau)}f^*(\tau)\d{s}\,H_{u_1, v_1, \nu_1^{-1}}g(t)\d{t}\\
&=\sup_{\|g\|_{X'(0,L)}\leq1}\int_0^L\Big(\xi(s)\sup_{\tau\in[s,L)}\frac1{\xi(\tau)}f^*(\tau)\Big)\Big(u_2(s)\int_{\nu_2^{-1}(s)}^Lv_2(t) u_1(t) \int_{\nu_1^{-1}(t)}^Lg(x)v_1(x)\d{x}\d{t}\Big)\d{s}\\
&\leq\Big\|\xi(t)\sup_{s\in[t,L)}\frac1{\xi(s)}f^*(s)\Big\|_{Z(0,L)}\sup_{\|g\|_{X'(0,L)}\leq1}\Big\|u_2(t)\int_{\nu_2^{-1}(t)}^Lv_2(s) u_1(s)\int_{\nu_1^{-1}(s)}^Lg(\tau)v_1(\tau)\d{\tau}\d{s}\Big\|_{Z'(0,L)}\\
&=\Big\|\xi(t)\sup_{s\in[t,L)}\frac1{\xi(s)}f^*(s)\Big\|_{Z(0,L)}\sup_{\|g\|_{X'(0,L)}\leq1}\Big\|u_2(t)\int_{\nu^{-1}(t)}^Lg(\tau)v_1(\tau)\int_{\nu_2^{-1}(t)}^{\nu_1(\tau)}v_2(s)u_1(s)\d{s}\d{\tau}\Big\|_{Z'(0,L)}\\
&\lesssim\Big\|\xi(t)\sup_{s\in[t,L)}\frac1{\xi(s)}f^*(s)\Big\|_{Z(0,L)}\sup_{\|g\|_{X'(0,L)}\leq1}\Big\|u_2(t)\int_{\nu^{-1}(t)}^Lg(\tau)v_1(\tau)\nu_1(\tau)u_1(\nu_1(\tau))v_2(\nu_1(\tau))\d{\tau}\Big\|_{Z'(0,L)}\\
&=\Big\|\xi(t)\sup_{s\in[t,L)}\frac1{\xi(s)}f^*(s)\Big\|_{Z(0,L)}\sup_{\|g\|_{X'(0,L)}\leq1}\|H_{u_2, v, \nu^{-1}}g\|_{Z'(0,L)}\\
&=\Big\|\xi(t)\sup_{s\in[t,L)}\frac1{\xi(s)}f^*(s)\Big\|_{Z(0,L)},
\end{align*}
for every $f\in\Mpl(0,L)$, where we used Fubini's theorem in the fourth and fifth equalities, the H\"older inequality \eqref{ch1:ri:holder} in the second inequality, the fact that $u_1 v_2$ satisfies the averaging condition~\eqref{averaging_condition} in the last inequality, and \eqref{thm:Riterationoptimaldomainupperbound:eq1} in the last equality. Therefore, the proof will be finished once we show that
\begin{equation*}
\Big\|\xi(t)\sup_{s\in[t,L)}\frac1{\xi(s)}f^*(s)\Big\|_{Z(0,L)}\lesssim\|R_{u_2, v, \nu}(f^*)\|_{X(0,L)}\quad\text{for every $f\in\Mpl(0,L)$}.
\end{equation*}
Since the function $\frac{\eta}{\xi}$ is equivalent to a nonincreasing function, we have that
\begin{equation*}
\Big\|\xi(t)\sup_{s\in[t,L)}\frac1{\xi(s)}f^*(s)\Big\|_{Z(0,L)}\lesssim\Big\|\eta(t)\sup_{s\in[t,L)}\frac1{\eta(s)}f^*(s)\Big\|_{Z(0,L)}
\end{equation*}
for every $f\in\Mpl(0,L)$. Hence it is sufficient to show that
\begin{equation}\label{thm:Riterationoptimaldomainupperbound:eq2}
\Big\|\eta(t)\sup_{s\in[t,L)}\frac1{\eta(s)}f^*(s)\Big\|_{Z(0,L)}\lesssim\|R_{u_2, v, \nu}(f^*)\|_{X(0,L)}\quad\text{for every $f\in\Mpl(0,L)$}.
\end{equation}

Note that, for every $f\in\Mpl(0,L)$,
\begin{align}\label{thm:Riterationoptimaldomainupperbound:eq3}
\Big\|\eta(t)\sup_{s\in[t,L)}\frac1{\eta(s)}f^*(s)\Big\|_{Z(0,L)}&\approx\Big\|v(t)\int_0^{\nu(t)} u_2(s) \eta(s)  \sup_{\tau\in[s,L)}\frac1{\eta(\tau)}f^*(\tau)\d{s}\Big\|_{X(0,L)}\notag\\
\begin{split}
&\leq\Big\|v(t)\int_0^{\nu(t)} u_2(s) \eta(s) \sup_{\tau\in[s,\nu(t))}\frac1{\eta(\tau)}f^*(\tau)\d{s}\Big\|_{X(0,L)}\\
&\quad+\Big\|v(t) \Big(\sup_{\tau\in[\nu(t),L)}\frac1{\eta(\tau)}f^*(\tau)\Big) \int_0^{\nu(t)} u_2(s) \eta(s) \d{s}\Big\|_{X(0,L)}
\end{split}
\end{align}
inasmuch as $\eta$ is equivalent to a nonincreasing function. Furthermore, since $\eta$ is equivalent to a nonincreasing function and satisfies \eqref{thm:Riterationoptimaldomainupperbound:integralassumption_on_eta_times_u2}, \cite[Theorem~3.2]{GOP:06} guarantees that
\begin{equation*}
\int_0^{\nu(t)} u_2(s) \eta(s) \sup_{\tau\in[s,\nu(t))}\frac1{\eta(\tau)}f^*(\tau)\d{s}\lesssim\int_0^{\nu(t)}f^*(s) u_2(s) \d{s}
\end{equation*}
for every $t\in(0,L)$ and every $f\in\Mpl(0,L)$, in which the multiplicative constant depends only on $C_2$. Hence
\begin{equation}\label{thm:Riterationoptimaldomainupperbound:eq4}
\begin{aligned}
\Big\|v(t)\int_0^{\nu(t)} u_2(s) \eta(s) \sup_{\tau\in[s,\nu(t))}\frac1{\eta(\tau)}f^*(\tau)\d{s}\Big\|_{X(0,L)}&\lesssim\Big\|v(t)\int_0^{\nu(t)} f^*(s) u_2(s) \d{s}\Big\|_{X(0,L)}\\
&=\|R_{u_2, v, \nu}(f^*)\|_{X(0,L)}
\end{aligned}
\end{equation}
for every $f\in\Mpl(0,L)$. Furthermore, thanks to the fact that $\eta$ satisfies \eqref{thm:Riterationoptimaldomainupperbound:integralassumption_on_eta_times_u2} again, we have that
\begin{align}\label{thm:Riterationoptimaldomainupperbound:eq5}
\Big\|v(t) \Big(\sup_{\tau\in[\nu(t),L)}\frac1{\eta(\tau)}f^*(\tau)\Big) \int_0^{\nu(t)}u_2(s) \eta(s) \d{s}\Big\|_{X(0,L)}&\lesssim\Big\|v(t)U_2(\nu(t))\eta(\nu(t))\sup_{\tau\in[\nu(t),L)}\frac1{\eta(\tau)}f^*(\tau)\Big\|_{X(0,L)}\notag\\
\begin{split}
&=\Big\|\sup_{\tau\in[\nu(t),L)}\frac1{\eta(\tau)}f^*(\tau)\Big\|_{X(0,L)}\\
&=\Big\|\sup_{\tau\in[t,L)}\frac1{\eta(\nu(\tau))}f^*(\nu(\tau))\Big\|_{X(0,L)}
\end{split}
\end{align}
for every $f\in\Mpl(0,L)$. We claim that
\begin{equation}\label{thm:Riterationoptimaldomainupperbound:eq6}
\Big\|\sup_{\tau\in[t,L)}\frac1{\eta(\nu(\tau))}f^*(\nu(\tau))\Big\|_{X(0,L)}\lesssim\|R_{u_2, v, \nu}(f^*)\|_{X(0,L)}.
\end{equation}
Thanks to the Hardy--Littlewood--P\'olya principle \eqref{ch1:ri:HLP}, it is sufficient to show that
\begin{equation}\label{thm:Riterationoptimaldomainupperbound:eq7}
\int_0^t\sup_{\tau\in[s,L)}\frac1{\eta(\nu(\tau))}f^*(\nu(\tau))\d{s}\lesssim\int_0^t(R_{u_2, v, \nu}(f^*))^*(s)\d{s}\quad\text{for every $t\in(0,L)$}.
\end{equation}
To this end, we have that
\begin{align}
\int_0^t\sup_{\tau\in[s,t)}\frac1{\eta(\nu(\tau))}f^*(\nu(\tau))\d{s}&\lesssim\int_0^t\frac1{\eta(\nu(s))}f^*(\nu(s))\d{s}=\int_0^t U_2(\nu(s)) v(s) f^*(\nu(s)) \d{s}\notag\\
&\leq\int_0^tR_{u_2, v, \nu}(f^*)(s)\d{s}\leq\int_0^t(R_{u_2, v, \nu}(f^*))^*(s)\d{s}\label{thm:Riterationoptimaldomainupperbound:eq8}
\end{align}
for every $t\in(0,L)$, where the first inequality follows from \cite[Theorem~3.2]{GOP:06} (the fact that the function $(0,L)\ni s\mapsto \frac1{\eta(\nu(s))}=U_2(\nu(s)) v(s)$ is equivalent to a nondecreasing function and satisfies \eqref{thm:Riterationoptimaldomainupperbound:integralassumptionon1overeta} was used here), the second inequality follows from the monotonicity of $f^*$, and the last one follows from the Hardy--Littlewood inequality \eqref{ch1:ri:HL}. Furthermore, owing to \eqref{thm:Riterationoptimaldomainupperbound:integralassumptionon1overeta} again, we have that
\begin{align}
\sup_{\tau\in[t,L)}\frac1{\eta(\nu(\tau))}f^*(\nu(\tau))&=\sup_{\tau\in[t,L)} U_2(\nu(\tau)) v(\tau) f^*(\nu(\tau))\lesssim\sup_{\tau\in[t,L)}\Big(\frac1{\tau}\int_0^\tau U_2(\nu(s)) v(s)\d{s}\Big)f^*(\nu(\tau))\notag\\
\begin{split}\label{thm:Riterationoptimaldomainupperbound:eq9}
&\leq\sup_{\tau\in[t,L)}\frac1{\tau}\int_0^\tau U_2(\nu(s)) v(s)f^*(\nu(s))\d{s}\leq\sup_{\tau\in[t,L)}\frac1{\tau}\int_0^\tau R_{u_2, v, \nu}(f^*)(s)\d{s}\\
&\leq\sup_{\tau\in[t,L)}\frac1{\tau}\int_0^\tau (R_{u_2, v, \nu}(f^*))^*(s)\d{s}=\frac1{t}\int_0^t (R_{u_2, v, \nu}(f^*))^*(s)\d{s}.
\end{split}
\end{align}
Inequality \eqref{thm:Riterationoptimaldomainupperbound:eq7} now follows from \eqref{thm:Riterationoptimaldomainupperbound:eq8} and \eqref{thm:Riterationoptimaldomainupperbound:eq9} inasmuch as
\begin{equation*}
\int_0^t\sup_{\tau\in[s,L)}\frac1{\eta(\nu(\tau))}f^*(\nu(\tau))\d{s}\leq\int_0^t\sup_{\tau\in[s,t)}\frac1{\eta(\nu(\tau))}f^*(\nu(\tau))\d{s} + t\sup_{\tau\in[t,L)}\frac1{\eta(\nu(\tau))}f^*(\nu(\tau))
\end{equation*}
for every $t\in(0,L)$.

Finally, by combining \eqref{thm:Riterationoptimaldomainupperbound:eq3} with \eqref{thm:Riterationoptimaldomainupperbound:eq4}, \eqref{thm:Riterationoptimaldomainupperbound:eq5} and \eqref{thm:Riterationoptimaldomainupperbound:eq6}, we obtain \eqref{thm:Riterationoptimaldomainupperbound:eq2}.
\end{proof}

\begin{remark}
Since \myref{Theorem}{thm:Riterationoptimaldomainupperbound} has several assumptions, it is instructive to provide a concrete, important example, which is also quite general. Let $\alpha_1, \alpha_2, \beta_1, \beta_2, \gamma_1, \gamma_2 \in(0,\infty)$. Set $\nu_j(t)=t^{\alpha_j}$, $u_j(t)=t^{\beta_j-1}b_j(t)$ and $v_j(t)=t^{\gamma_j-1}c_j(t)$, $t\in(0,L)$, $j=1,2$, where $b_j$, $c_j$ are continuous slowly-varying functions. Set $d=(b_1\circ\nu_1) \cdot c_1 \cdot (c_2\circ\nu_1)$ and $\widetilde{d}=(b_1\circ\nu_1) \cdot c_1$. Assume that $\gamma_2 < 1$, $\beta_1 + \gamma_2 > 1$, and
\begin{equation*}
\alpha_1(\beta_1 + \alpha_2\beta_2 + \gamma_2 - 1) + \gamma_1 \geq 1,\ \alpha_1(\beta_1 +\alpha_2\beta_2 - \alpha_2) + \gamma_1 \geq 1,\ \alpha_1(\beta_1 + \gamma_2 - 1) + \gamma_1 < 1.
\end{equation*}
If $\alpha_1(\beta_1 + \alpha_2\beta_2 + \gamma_2 - 1) + \gamma_1=1$ or $\alpha_1(\beta_1 +\alpha_2\beta_2 - \alpha_2) + \gamma_1 = 1$, also assume that $d$ or $\widetilde{d}$, respectively, is equivalent to a nondecreasing function. Under these assumptions, we can use \myref{Theorem}{thm:Riterationoptimaldomainupperbound} to obtain that
\begin{equation*}
\Big\|v_1(t)\int_0^{t^{\alpha_1}}\big[v_2(\tau)\int_0^{\tau^{\alpha_2}}f^*(\sigma)u_2(\sigma)\d{\sigma}\big]^*(s)u_1(s)\d{s}\Big\|_{X(0,L)}\approx\Big\|t^\delta d(t)\int_0^{t^{\alpha_1\alpha_2}}f^*(s)s^{\beta_2-1}\d{s}\Big\|_{X(0,L)}
\end{equation*}
for every $f\in\Mpl(0,L)$, where $\delta=\alpha_1(\beta_1 + \gamma_2 - 1) + \gamma_1 - 1$.

When $\beta_j=1$ and $b_j=c_j\equiv1$, $j=1,2$, the assumptions are satisfied provided that
\begin{equation}\label{rem:Riterationoptimaldomainupperbound:specialcase}
\alpha_1(\alpha_2 + \gamma_2) + \gamma_1 \geq 1,\ \alpha_1 + \gamma_1 \geq 1,\ \alpha_1\gamma_2 + \gamma_1 < 1.
\end{equation}
In particular, \eqref{rem:Riterationoptimaldomainupperbound:specialcase} is satisfied if (cf.~\cite[Theorem~3.4]{CP:16})
\begin{equation*}
\alpha_2+\gamma_2\geq1,\ \alpha_1+\gamma_1\geq1,\ \alpha_1\gamma_2 + \gamma_1<1.
\end{equation*}
\end{remark}

We conclude this paper with a $H_{u, v, \nu}$ counterpart to \myref{Theorem}{thm:Riterationoptimaldomainupperbound}, whose proof is substantially simpler than that of the theorem.
\begin{proposition}\label{prop:Hiteration}
Let $\|\cdot\|_{X(0,L)}$ be a rearrangement\hyp{}invariant function norm. Let $\nu_1,\nu_2\colon (0,L)\to(0,\infty)$ be increasing bijections. Assume that $\nu_1\in\Dinf{0}$. If $L=\infty$, assume that $\nu_1\in\Dinf{\infty}$. Let $u_1, u_2, v_1, v_2\colon(0,L)\to(0,\infty)$ be measurable. Assume that the function $v_1 u_2$ is equivalent to a nonincreasing function and that it satisfies the averaging condition~\eqref{averaging_condition}.  Set
\begin{equation*}
v(t)=\nu_2^{-1}(t)v_1(\nu_2^{-1}(t))u_2(\nu_2^{-1}(t))v_2(t),\ t\in(0,L),
\end{equation*}
and $\nu=\nu_2\circ\nu_1$. We have that
\begin{equation}\label{prop:Hiteration:H1_circ_H2}
\|H_{u_1, v_1, \nu_1}(H_{u_2, v_2, \nu_2}f)\|_{X(0,L)}\approx\|H_{u_1, v, \nu}f\|_{X(0,L)}\quad\text{for every $f\in\Mpl(0,L)$},
\end{equation}
in which the multiplicative constants depend only on $\nu_1$, the averaging constant of $v_1 u_2$, and the multiplicative constants in the equivalence of $v_1 u_2$ to a nonincreasing function.

If, in addition, $u_1$ and $u_2$ are nonincreasing, if there is a measurable function $\xi\colon (0,L) \to (0,\infty)$ such that
\begin{equation*}
\frac1{v_1(t)} = \int_0^{\nu_1^{-1}(t)} \xi(s) \d{s} \quad \text{for every $t\in(0,L)$},
\end{equation*}
and if the operator $T_\varphi$ defined by \eqref{opTdef} with $\varphi = \frac{u_1}{\xi}$ is bounded on $X'(0,L)$,
then
\begin{equation*}
\sup_{\substack{g\sim f\\g\in\Mpl(0,L)}} \sup_{\substack{h \sim H_{u_2, v_2, \nu_2} g\\h\in\Mpl(0,L)}} \| H_{u_1, v_1, \nu_1}h\|_{X(0,L)} \approx \sup_{\substack{g\sim f\\g\in\Mpl(0,L)}} \|H_{u_1, v, \nu}g\|_{X(0,L)}\quad\text{for every $f\in\Mpl(0,L)$},
\end{equation*}
in which the multiplicative constants depend only on the norm of $T_\varphi$ on $X'(0,L)$ and the multiplicative constant in \eqref{prop:Hiteration:H1_circ_H2}.
\end{proposition}
\begin{proof}
On the one hand, we have that
\begin{align*}
\|H_{u_1, v_1, \nu_1}(H_{u_2, v_2, \nu_2}f)\|_{X(0,L)}&=\Big\|u_1(t) \int_{\nu_1(t)}^L \Big( u_2(s)\int_{\nu_2(s)}^Lf(\tau)v_2(\tau)\d{\tau}\Big)  u_2(s) v_1(s) \d{s}\Big\|_{X(0,L)}\notag\\
\begin{split}
&=\Big\|u_1(t) \int_{\nu(t)}^Lf(\tau) v_2(\tau) \int_{\nu_1(t)}^{\nu_2^{-1}(\tau)}  u_2(s) v_1(s) \d{s}\d{\tau}\Big\|_{X(0,L)}\\
&\leq\Big\|u_1(t) \int_{\nu(t)}^Lf(\tau)v_2(\tau)\int_0^{\nu_2^{-1}(\tau)} u_2(s) v_1(s)\d{s}\d{\tau}\Big\|_{X(0,L)}\\
&\lesssim\Big\|u_1(t) \int_{\nu(t)}^Lf(\tau)v_2(\tau)\nu_2^{-1}(\tau)u_2(\nu_2^{-1}(\tau)) v_1(\nu_2^{-1}(\tau))\d{\tau}\Big\|_{X(0,L)}\\
&=\|H_{u_1, v, \nu}f\|_{X(0,L)}
\end{split}
\end{align*}
for every $f\in\Mpl(0,L)$ thanks to the fact that $v_1 u_2$ satisfies the averaging condition~\eqref{averaging_condition}. 

As for the opposite inequality, observe that $M=\inf_{t\in(0,\frac{L}{\theta})}\frac{\nu_1(\theta t)}{\nu_1(t)}>1$, where $\theta>1$ is such that $\nu_1\in\Dinf[\theta]{0}$ and, if $L=\infty$, also $\nu_1\in\Dinf[\theta]{\infty}$. Set $K=\min\{\frac1{\theta}, \nu_1^{-1}(\frac1{M})\}$. We have that
\begin{align*}
&\|H_{u_1, v_1, \nu_1}(H_{u_2, v_2, \nu_2}f)\|_{X(0,L)}\\
&=\Big\|u_1(t) \int_{\nu(t)}^Lf(\tau)v_2(\tau)\int_{\nu_1(t)}^{\nu_2^{-1}(\tau)} u_2(s) v_1(s)\d{s}\d{\tau}\Big\|_{X(0,L)}\\
&\geq\Big\|\chi_{(0,KL)}(t) u_1(t)\int_{\nu_2(M\nu_1(t))}^Lf(\tau)v_2(\tau)\int_{\nu_1(t)}^{\nu_2^{-1}(\tau)} u_2(s) v_1(s)\d{s}\d{\tau}\Big\|_{X(0,L)}\\
&\gtrsim\Big\|\chi_{(0,KL)}(t) u_1(t)\int_{\nu_2(M\nu_1(t))}^Lf(\tau)v_2(\tau)u_2(\nu_2^{-1}(\tau)) v_1(\nu_2^{-1}(\tau))(\nu_2^{-1}(\tau)-\nu_1(t))\d{\tau}\Big\|_{X(0,L)}\\
&\geq\frac{M-1}{M}\Big\|\chi_{(0,KL)}(t) u_1(t)\int_{\nu_2(M\nu_1(t))}^Lf(\tau)v(\tau)\d{\tau}\Big\|_{X(0,L)}\\
&\geq\frac{M-1}{M}\Big\|\chi_{(0,KL)}(t) u_1(t)\int_{\nu_2(\nu_1(\theta t))}^Lf(\tau)v(\tau)\d{\tau}\Big\|_{X(0,L)}\\
&\geq\frac{M-1}{M}\Big\|\chi_{(0,L)}\Big(\frac{t}{K}\Big) u_1(t)\int_{\nu_2(\nu_1(\frac{t}{K}))}^Lf(\tau)v(\tau)\d{\tau}\Big\|_{X(0,L)}\\
&\geq \frac{M-1}{M}K\Big\|u_1(t) \int_{\nu(t)}^Lf(\tau)v(\tau)\d{\tau}\Big\|_{X(0,L)}
\end{align*}
for every $f\in\Mpl(0,L)$, where we used the fact that $v_1u_2$ is equivalent to a nonincreasing function and the boundedness of the dilation operator $D_{\frac1{K}}$ (see~\eqref{ch1:ri:dilation}).

Finally, under the additional assumptions, we have that
\begin{align*}
\sup_{\substack{g\sim f\\g\in\Mpl(0,L)}} \sup_{\substack{h \sim H_{u_2, v_2, \nu_2} g\\h\in\Mpl(0,L)}} \| H_{u_1, v_1, \nu_1}h\|_{X(0,L)} &\approx \sup_{\substack{g\sim f\\g\in\Mpl(0,L)}} \| H_{u_1, v_1, \nu_1}(H_{u_2, v_2, \nu_2} g)\|_{X(0,L)}\\
&\approx \sup_{\substack{g\sim f\\g\in\Mpl(0,L)}} \|H_{u_1, v, \nu}g\|_{X(0,L)}
\end{align*}
for every $f\in\Mpl(0,L)$ thanks to \eqref{prop:norminducedbyT:equivalencewithHf*} combined with \eqref{prop:Hiteration:H1_circ_H2}.
\end{proof}

\begin{remark}
If $T_{\varphi}$ is not bounded on $X'(0,\infty)$, then, while we still have that
\begin{equation*}
\sup_{\substack{g\sim f\\g\in\Mpl(0,L)}} \sup_{\substack{h \sim H_{u_2, v_2, \nu_2} g\\h\in\Mpl(0,L)}} \| H_{u_1, v_1, \nu_1}h\|_{X(0,L)} \gtrsim \sup_{\substack{g\sim f\\g\in\Mpl(0,L)}} \|H_{u_1, v, \nu}g\|_{X(0,L)}\quad\text{for every $f\in\Mpl(0,L)$},
\end{equation*}
it remains an open question whether the opposite inequality (is)/(can be) valid.
\end{remark}

\bibliography{bibliography}
\end{document}